\documentclass[authoryear]{elsarticle}

\usepackage{multirow}

\usepackage{graphicx}  
\usepackage[font=small,format=plain,labelfont=bf,up,textfont=it,up]{caption}

\usepackage{subfigure}
\usepackage{geometry}
\usepackage{amsmath}
\usepackage{amssymb}
\usepackage{amsfonts}
\usepackage{amsthm}

\usepackage{fancyhdr}
\usepackage{float}

\usepackage{pgf,tikz}
\usetikzlibrary{shapes,arrows,calc,automata}

\usepackage{ctable}

\let\oldref\ref

\renewcommand{\ref}[1]{(\oldref{#1})}

\title{\textbf{Concurrent Implicit Spectral Deferred Correction Scheme \\ for Low-Mach Number Combustion \\ with Detailed Chemistry}}
\date{}

\author[lbl]{Fran\c cois P. Hamon\corref{cor1}}
\ead{fhamon@lbl.gov}
\author[lbl]{Marcus S. Day}

\author[lbl]{Michael L. Minion}

\cortext[cor1]{Corresponding author}
\address[lbl]{Center for Computational Sciences and Engineering, Lawrence Berkeley National Laboratory, Berkeley, USA}

\begin{document}

\begin{abstract} 

We present a parallel multi-implicit time integration scheme for the advection-diffusion-reaction systems arising from 
the equations governing low-Mach number combustion with complex chemistry. Our strategy employs parallelization across 
the method to accelerate the serial Multi-Implicit Spectral Deferred Correction (MISDC) scheme used to couple the advection, 
diffusion, and reaction processes. In our approach, the diffusion solves and the reaction solves are performed concurrently 
by different processors. Our analysis shows that the proposed parallel scheme is stable for stiff problems and that the sweeps converge to 
the fixed-point solution at a faster rate than with serial MISDC. We present numerical examples to demonstrate that the new algorithm 
is high-order accurate in time, and achieves a parallel speedup compared to serial MISDC.

\end{abstract}


\maketitle

\section{\label{section_introduction}Introduction} 

Many reacting flow problems are modeled by advection-diffusion-reaction partial differential equations. These systems are often characterized by a 
large disparity in the time scales of these three processes, therefore making the design of accurate and efficient time integration schemes 
particularly challenging. This is the case for the low-Mach number combustion approaches used to simulate laboratory-scale flames. Low-Mach number 
combustion models are obtained from the fully compressible equations by an asymptotic analysis which eliminates the fast acoustic waves 
while still accounting for compressibility effects caused by diffusion and reaction processes \citep{majda1985derivation}. This methodology is 
advantageous for the computational efficiency of the time integration scheme because low-Mach number combustion models can therefore be numerically 
evolved on the time scale of the relatively slow advection process. The resulting stability constraint is much less restrictive than in fully 
compressible models and significantly larger stable time steps can be used. However, this approach requires a sophisticated integration method to 
couple advection with the diffusion and reaction processes which operate on much faster time scales. The coupling strategy becomes a 
key determinant of the largest time step size that preserves the accuracy, stability, and computational efficiency of the overall numerical 
simulation framework.

A popular approach to advance low-Mach number reactive systems is operator splitting, 
in which the different physical processes are decoupled and solved sequentially. Operator splitting 
is attractive because it allows specialized schemes and solvers to be used for advection, diffusion, 
and reaction, with different time step sizes depending on the stiffness of each process. Applications
to low-Mach number combustion include the second-order Strang-split schemes of \cite{knio1999semi,day2000numerical}.
The former is based on an explicit evaluation of the advection and diffusion terms and an implicit treatment of the 
reaction term. The latter couples an explicit scheme for advection with a Crank-Nicolson treatment of diffusion and a 
fully implicit method for reactions; the implicit diffusion treatment accommodates that diffusive time scales in 
combustion can be much shorter than those of advective transport.
Several authors present alternatives to address the stiffness of the diffusion terms relative to advection 
via extended-stability Runge-Kutta-Chebyshev schemes \citep{najm2005modeling,safta2010high,motheau2016high}. 
However, in all such attempts, the underlying splitting approach introduces significant time-splitting errors,
and therefore requires a dynamic adaptation of the time step sizes 
to preserve the accuracy of the approximation, as in \cite{duarte2013time}. 
In addition to this drawback, a significant limitation of operator-splitting approaches in general
is that they cannot be easily extended beyond 
second-order accuracy in time.

Recently, the Multi-Implicit Spectral Deferred Correction (MISDC) scheme of \cite{bourlioux2003high,layton2004conservative} has been successfully 
employed to couple advection, diffusion, and reaction processes in the low-Mach number equation set \citep{nonaka2012deferred,pazner2016high}. 
The methodology builds on the Semi-Implicit Spectral Deferred Correction (SISDC) scheme of \cite{minion2003semi} and aims at handling
a right-hand side with one or more explicit terms coupled with multiple implicit terms. Specifically, MISDC relies on a temporal splitting in which advection 
is treated explicitly while diffusion and reaction are treated implicitly but decoupled as in operator splitting. 
In MISDC, the integration is based on a sequence of low-order corrections applied to the variables to iteratively achieve high-order accuracy. 
At a given Spectral Deferred Correction (SDC) temporal node, the diffusion and reaction corrections are computed separately in two successive implicit solves. Numerical results 
showed that MISDC is able to accurately resolve stiff kinetics while capturing the coupling between chemical reactions and transport processes. 
The finite-volume schemes based on MISDC avoid the often catastrophic splitting errors of Strang splitting, while also reducing the cost 
of the implicit reaction solves
and achieving second-order space-time accuracy in \cite{nonaka2012deferred}, and fourth-order in \cite{pazner2016high}.

The implicit solves corresponding to the diffusion and reaction steps still represent most of the computational cost incurred by the multi-implicit time 
integration scheme for these stiff advection--diffusion--reaction systems. In the numerical methodology for low-Mach number combustion presented in 
\cite{pazner2016high}, the treatment of the diffusion 
and reaction steps is sequential to avoid solving a global diffusion-reaction nonlinear system involving all the degrees of freedom. Instead, the 
diffusion step requires solving a linearized system  for each mass conservation equation and for the energy equation to update the species mass fractions 
and the enthalpy. Then, the diffusion step is followed by the reaction step to update the species mass 
fractions. The latter entails solving complex, but local, nonlinear Newton-based backward-Euler systems to evolve the multistep kinetic network 
of the combustion process under consideration. These sequential solves become very expensive for detailed combustion problems with a large number of 
species involved in stiff and highly nonlinear chemical reaction paths.

In this paper, we design a multi-implicit spectral deferred correction-based integration method in which the implicit diffusion and reaction steps 
are performed concurrently to overcome the sequential limitation of the standard MISDC method. Our algorithm, referred to as Concurrent Implicit 
SDC (CISDC), relies on \textit{parallelization across the method} according to the classification of \cite{burrage1997parallel}.
Examples of parallel-across-the-method approaches include Runge-Kutta schemes in which intermediate stage values 
are computed in parallel \citep{iserles1990theory,butcher1997order}, and SDC schemes in which the corrections are 
applied at multiple temporal nodes concurrently \citep{christlieb2010parallel,speck2018parallelizing}. This strategy 
can be contrasted with \textit{parallelization across the steps} 
\citep{nievergelt1964parallel,miranker1967parallel,lions2001resolution,emmett2012toward,falgout2014parallel} in which 
the parallel work is done on multiple time steps concurrently, and \textit{parallelization across the problem} (e.g. 
wave-form relaxation methods \citep{1270004, gander1999waveform}), in which the problem is divided into temporal subproblems 
solved simultaneously. 

In the proposed CISDC algorithm, we formally decouple the diffusion step at a given SDC temporal node from the reaction step at the 
previous node. These two steps are therefore independent and can be executed in parallel by multiple processors. In this work, 
we focus primarily on the convergence of the correction iterations to show that the parallel scheme retains robust 
stability properties for stiff problems, which is key for combustion problems with complex chemistry. We also analyze the theoretical 
computational cost of the proposed CISDC scheme. Then, we modify the numerical methodology of \cite{pazner2016high} to perform the temporal 
integration of the low-Mach number equation set with our parallel-across-the-method scheme. Finally, we demonstrate the robustness and 
efficiency of the new approach using synthetic test cases and a challenging one-dimensional unsteady flame simulation.

We proceed to the presentation of the low-Mach number equation set in Section \oldref{section_governing_equations}. In Section 
\oldref{section_parallel_spectral_deferred_correction_scheme}, we propose our parallel Concurrent Implicit SDC scheme in which the 
implicit solves can be performed in parallel to reduce the computational cost. The stability and convergence properties of the 
scheme are analyzed in Section \oldref{subsection_convergence_analysis}. Numerical examples are presented in Section \oldref{section_numerical_examples}.

\section{\label{section_governing_equations}Low-Mach number governing equations}

Following previous work by \cite{day2000numerical,nonaka2012deferred,pazner2016high}, we consider a low-Mach number model presented by 
\cite{baum1978equations} and derived from an asymptotic analysis in \cite{majda1985derivation}. It consists of a system of partial differential 
equations governing the evolution of a gaseous mixture in an open container in a non-gravitationally stratified environment.
The system describes coupled advection, differential/preferential diffusion, and chemical reaction processes, and relies on a mixture model for species 
diffusion \citep{kee1983fortran,warnatz1982numerical}.  The Soret and Dufour transport effects are ignored. Sound waves are analytically eliminated from the system, but the local compressibility 
effects caused by reactions and diffusion are included in the model. The system is closed by an equation of state (EOS), written as a constraint
on the divergence of the velocity, which acts to ensure that the thermodynamic pressure remains constant over space and time. Physically, this approach amounts to enforcing an instantaneous
acoustic equilibration of the entire system to the constant ambient pressure.

Considering a mixture with $N$ species, conservation of mass for species $i \in \{1, \dots, N\}$, and conservation of energy are expressed by the evolution
equations
\begin{align}
\frac{\partial (\rho Y_i)}{\partial t} &= - \nabla \cdot ( U \rho Y_i ) + \nabla \cdot \rho \mathcal{D}_i \nabla Y_i + \dot{\omega}_i, \qquad i \in \{ 1, \dots, N \}, \label{species_equation} \\
\frac{\partial ( \rho h )}{\partial t} &= - \nabla \cdot ( U \rho h )  + \nabla \cdot \frac{\lambda}{c_p} \nabla h + \sum_i \nabla \cdot h_i \bigg( \rho \mathcal{D}_i - \frac{\lambda}{c_p} \bigg) \nabla Y_i,
\label{enthalpy_equation}
\end{align}
where $\rho$ is the density, $\boldsymbol{Y} = [ Y_1, \dots, Y_N ]$ are the species mass fractions, $\mathcal{D}_i(\boldsymbol{Y},T)$ is the 
mixture-averaged diffusion coefficient of species $i$, $T$ is the temperature, $\dot{\omega}_i(\boldsymbol{Y},T)$ is the production rate of 
species $i$ due to chemical reactions, $h(\boldsymbol{Y},T) = \sum_i Y_i h_i(T)$ is the enthalpy with $h_i(T)$ denoting the enthalpy of 
species $i$, $\lambda(\boldsymbol{Y},T)$ is the thermal conductivity, and $c_p = \sum_i Y_i d h_i / d T$ is the specific heat at constant 
pressure. In this formulation, there is no exchange in enthalpy due to chemical reactions since $h$ includes the standard enthalpy of formation. 
These evolution equations are closed by an EOS, which states that the thermodynamic pressure, $p_{\textit{EOS}}$, must remain constant and equal to the 
ambient pressure, $p_0$,
\begin{equation}
p_0 = p_{\textit{EOS}}, \label{thermodynamic_constraint_eos}
\end{equation}
where $p_{\textit{EOS}}$ is computed as
\begin{equation}
p_{\textit{EOS}} := \rho \mathcal{R} T \sum_i \frac{Y_i}{W_i}.
\label{p_eos_definition}
\end{equation}
$\mathcal{R}$ is the universal gas constant and $W_i$ is the molecular weight of species $i$. In the mixture model considered here, the 
diffusion flux, $\Gamma_i$, of species $i$, is defined as $\Gamma_i := - \rho \mathcal{D}_i \nabla Y_i$. To enforce that the sum of the fluxes 
is equal to zero and therefore to guarantee mass conservation, we define a dominant species, $i_0$, whose diffusion flux is set to 
$\Gamma_{i_0} := - \sum_{i \neq i_0} \Gamma_i$. Using this approach, one can sum the species equation \ref{species_equation} to obtain the 
continuity equation, 
\begin{equation}
\frac{\partial \rho}{\partial t} = - \nabla \cdot ( U \rho ), \label{continuity_equation}
\end{equation}
where we used the constraints
\begin{equation}
\sum_i Y_i = 1,
\qquad 
\text{and}
\qquad
\sum_i \dot{\omega}_i = 0.
\end{equation}
A possible approach consists in evolving all the thermodynamic variables but one with \ref{species_equation}-\ref{enthalpy_equation},
and then use the EOS \ref{thermodynamic_constraint_eos} to compute the last thermodynamic variable \citep{najm1998semi,knio1999semi}. 
However, this approach fails to strictly conserve energy. Instead, we follow the volume discrepancy method of \cite{pember1998adaptive,day2000numerical}. That is,
the EOS \ref{thermodynamic_constraint_eos} is written in the form of a 
constraint on the divergence of the velocity. This is achieved by taking the derivative of \ref{thermodynamic_constraint_eos} in the 
Lagrangian frame while enforcing that the thermodynamic pressure remains constant, and then substituting the evolution 
equations \ref{species_equation}-\ref{enthalpy_equation} for $\rho$, $\boldsymbol{Y}$, and $T$ to obtain
\begin{equation}
\nabla \cdot U = S. \label{divergence_constraint}
\end{equation}
The quantity $S$ is defined as
\begin{equation}
S := \frac{1}{\rho c_p T} \bigg( \nabla \cdot \lambda \nabla T + \sum_i \Gamma_i \cdot \nabla h_i \bigg) 
             + \frac{1}{\rho} \sum_i \frac{W}{W_i} \nabla \cdot \Gamma_i + \frac{1}{\rho} \sum_i \bigg( \frac{W}{W_i} - \frac{h_i}{c_p T} \bigg) \dot{\omega}_i, \label{S_def}
\end{equation}
where $W = (\sum_i Y_i/W_i )^{-1}$ denotes the mixture-averaged molecular weight. Equations \ref{divergence_constraint}-\ref{S_def} 
represent a linearized approximation to the velocity field required to hold the thermodynamic pressure equal to $p_0$ in the presence 
of local compressibility effects due to reaction heating, compositional changes, and thermal diffusion.  In summary, the system of partial 
differential equations solved in this work is given by \ref{species_equation}, \ref{enthalpy_equation}, and \ref{continuity_equation}, 
coupled with the velocity constraint, \ref{divergence_constraint}-\ref{S_def}. As an aside, since we only consider the one-dimensional 
case in the present work, the numerical scheme does not require a velocity projection; the velocity is fully determined via the divergence constraint 
\ref{divergence_constraint} and the inflow Dirichlet boundary condition specified in the numerical examples. Next, we proceed to the 
description of the temporal discretization applied to these PDEs, with an emphasis on the proposed parallel spectral deferred correction scheme.

\section{\label{section_parallel_spectral_deferred_correction_scheme}Parallel spectral deferred correction scheme}

\subsection{\label{subsection_multi_implicit_spectral_deferred_correction_schemes}Multi-implicit spectral deferred correction schemes}

We start with a review of the fundamentals of the Spectral Deferred Correction (SDC) and Multi-Implicit Spectral Deferred Correction (MISDC) schemes. 
Consider the ODE on a generic time step
\begin{align}
\phi_t(t) &= F_A(\phi(t)) + F_D(\phi(t)) +F_R(\phi(t)), \qquad t \in [t^n,t^{n} + \Delta t], \\
\phi(t^n) &= \phi^n,
\end{align}
and its solution in integral form given by
\begin{equation}
  \phi(t) = \phi^n + \int^t_{t^n} F( \phi(\tau) )d \tau,
\end{equation}
where $F_A$, $F_D$, and $F_R$ represent the advection, diffusion, and reaction terms, respectively, with $F = F_A + F_D + F_R$. 
We denote by $\tilde{\phi}(t)$ the approximation of $\phi(t)$, and we define the correction $\delta(t) := \phi(t) - \tilde{\phi}(t)$.
The SDC scheme of \cite{dutt2000spectral} iteratively improves the accuracy of the approximation with the update equation
\begin{align}
\tilde{\phi}(t) + \delta(t) & = \phi^n + \int_{t^{n}}^{t} \big[ F\big( \tilde{\phi}(\tau) + \delta(\tau) \big) - F\big( \tilde{\phi}(\tau) \big) \big] d \tau \nonumber \\
               & \qquad \, \, \,  + \int_{t^n}^{t} F\big( \tilde{\phi}(\tau) \big)  d \tau.
              \label{standard_sdc_sweep}
\end{align}
In \ref{standard_sdc_sweep}, the first integral is approximated with a low-order discretization such as backward- or forward-Euler, 
whereas the second integral is approximated with a high-order quadrature rule. The resulting discrete update is applied iteratively in sweeps to 
increase the order of accuracy of the approximation. Specifically, each sweep increases the formal order of accuracy by one until the order of 
accuracy of the quadrature applied to the second integral is reached \citep{hagstrom2007spectral,xia2007efficient,christlieb2009comments}.

MISDC, proposed in \cite{bourlioux2003high,layton2004conservative}, and used in \cite{pazner2016high}, is based on a variant of the original SDC update equation \ref{standard_sdc_sweep}. It is 
well suited for advection-diffusion-reaction problems in which the time scales corresponding to these three processes are very different, which is the 
case for low-Mach number combustion with complex chemistry. MISDC provides a methodology to treat these processes sequentially while accounting for the 
physical coupling between them to minimize the splitting error. This approach is based on the update equations
\begin{align}
\tilde{\phi}(t) + \delta_A(t) &= \phi^n + \int_{t^{n}}^{t} \big[ F_A\big( \tilde{\phi}(\tau) + \delta_A(\tau) \big) - F_A\big(  \tilde{\phi}(\tau) \big) \big] d \tau \nonumber \\ 
                   &  \qquad \, \, \, + \int_{t^n}^{t} F\big( \tilde{\phi}(\tau \big) d \tau, \label{update_equation_misdc_integral_form_A} \\[10pt]
\tilde{\phi}(t) + \delta_{AD}(t) &= \phi^n + \int_{t^{n}}^{t} \big[ F_A\big( \tilde{\phi}(\tau) + \delta_A(\tau) \big) - F_A\big(  \tilde{\phi}(\tau) \big) + F_D\big( \tilde{\phi}(\tau) + \delta_{AD}(\tau) \big) - F_D\big(  \tilde{\phi}(\tau) \big) \big] d \tau \nonumber \\ 
                 &  \qquad \, \, \, + \int_{t^n}^{t} F\big( \tilde{\phi}(\tau) \big) d \tau, \label{update_equation_misdc_integral_form_D} \\[10pt]
\tilde{\phi}(t) + \delta(t) &= \phi^n + \int_{t^{n}}^{t} \big[ F_A\big( \tilde{\phi}(\tau) + \delta_A(\tau) \big) - F_A\big(  \tilde{\phi}(\tau) \big) + F_D\big( \tilde{\phi}(\tau) + \delta_{AD}(\tau) \big) - F_D\big(  \tilde{\phi}(\tau) \big)  \nonumber \\ 
&  \qquad \, \, \, + F_R\big(  \tilde{\phi}(\tau) + \delta(\tau) \big) - F_R\big(  \tilde{\phi}(t) \big) \big] d \tau + \int_{t^n}^{t} F\big( \tilde{\phi}(\tau) \big) d \tau. \label{update_equation_misdc_integral_form_R} 
\end{align}
To discretize the update equations, SDC-based methods rely on a decomposition of the time interval $[t^n, t^{n+1}]$ into $M$ subintervals using $M+1$ temporal nodes, such that
\begin{equation}
  t^n = t^{n,0} < t^{n,1} < \dots < t^{n,M} = t^n + \Delta t = t^{n+1}.
\end{equation}
In this work, we consider Gauss-Lobatto nodes for the definition of the subintervals.
For brevity, we use the shorthand notations $t^m = t^{n,m}$ and $\Delta t^m = t^{m+1} - t^m$ in the remainder of the paper. We denote by $\phi^{m,(k)}$ the 
approximation of $\phi(t^m)$ at sweep $(k)$. 
In the application considered here, the diffusion and reaction terms are very stiff compared to the advection term which operates on a much slower 
time scale. Therefore, in MISDC the diffusion and reaction terms are treated implicitly and discretized with a backward-Euler method, whereas the advection term is 
treated explicitly and discretized with a forward-Euler method. Based on this temporal splitting, computing $\phi^{m+1,(k+1)}_A$ is not necessary, and the 
updates \ref{update_equation_misdc_integral_form_A}-\ref{update_equation_misdc_integral_form_D}-\ref{update_equation_misdc_integral_form_R} simplify to
the discrete update equations
\begin{align}
\qquad \quad \phi_{AD}^{m+1,(k+1)} = \phi^{m,(k+1)} + \Delta t^m &\big[ F_A(  \phi_{}^{m,(k+1)} ) - F_A( \phi^{m,(k)} ) \nonumber \\ 
+ & \, \, F_D( \phi_{AD}^{m+1,(k+1)} ) - F_D(  \phi^{m+1,(k)} ) \big] + \Delta t S^{m:m+1}\big( F(\phi^{(k)}) \big), \label{update_equation_misdc_discrete_form_D} \\[10pt]
\phi^{m+1,(k+1)} = \phi^{m,(k+1)} + \Delta t^m &\big[ F_A(  \phi_{}^{m,(k+1)} ) - F_A( \phi^{m,(k)} ) \nonumber \\ 
+ & \, \, F_D( \phi_{AD}^{m+1,(k+1)} ) - F_D(  \phi^{m+1,(k)} ) \big] \nonumber \\
+ & \, \, F_R( \phi^{m+1,(k+1)} ) - F_R( \phi^{m+1,(k)} ) \big] + \Delta t S^{m:m+1}\big( F(\phi^{(k)}) \big), \label{update_equation_misdc_discrete_form_R}
\end{align}
where $S^{m:m+1}\big( F(\phi^{(k)}) \big)$ is a high-order numerical quadrature approximating the last integral in \ref{standard_sdc_sweep} over 
the interval between two consecutive SDC nodes using the Lagrange polynomials $L_j$, $j \in \{0, \dots, M \}$,
\begin{equation}
S^{m:m+1}\big( F( \phi^{(k)} ) \big) := \sum_{j = 0}^{M} s_{m+1,j} F( \phi^{(k)} ),
\end{equation}
with
\begin{equation}
 s_{m+1,j} := \frac{1}{\Delta t} \int_{t^m}^{t^{m+1}} L_j( \tau ) d \tau.
\end{equation}

\cite{weiser2015faster} proposed a new class of SDC schemes in which the SDC sweep is cast as a stage of a diagonally implicit Runge-Kutta method. 
The choice of quadrature weights, based on LU decomposition, leads to a faster convergence of the iterative correction process than simpler backward-Euler. 
This new approach is also advantageous because in the resulting scheme, the sweeps remain convergent when the underlying problem is very stiff. Here, 
we adapt the scheme of \cite{weiser2015faster} to our multi-implicit framework based on Gauss-Lobatto nodes by writing the update equations 
for $m \in \{ 0, \dots, M-1 \}$ as
\begin{align}
\phi_{AD}^{m+1,(k+1)} = \phi^{0}  + \Delta t    & \sum_{j = 1}^m \tilde{q}^E_{m+1,j}  \big[ F_A(  \phi^{j,(k+1)} ) - F_A( \phi^{j,(k)} ) \big]   \nonumber \\ 
                             + \, \Delta t &  \sum_{j = 1}^m \tilde{q}^I_{m+1,j} \big[ F_D(  \phi^{j,(k+1)} ) - F_D( \phi^{j,(k)} ) \big]  \nonumber \\
                             + \, \Delta t & \tilde{q}^I_{m+1,m+1}              \big[ F_D(  \phi_{AD}^{m+1,(k+1)} ) - F_D(  \phi^{m+1,(k)} ) \big] \nonumber \\
                             + \, \Delta t & Q^{0:m+1}\big( F(\phi^{(k)}) \big), \label{update_equation_misdcq_discrete_form_D} \\[5pt]
\phi^{m+1,(k+1)}     =\phi^{0} + \Delta t     & \sum_{j = 1}^m \tilde{q}^E_{m+1,j} \big[ F_A(  \phi^{j,(k+1)} ) - F_A( \phi^{j,(k)} ) \big]   \nonumber \\ 
                             + \, \Delta t &  \sum_{j = 1}^m \tilde{q}^I_{m+1,j} \big[ F_D(  \phi^{j,(k+1)} ) - F_D( \phi^{j,(k)} ) + F_R(  \phi^{j,(k+1)} ) - F_R( \phi^{j,(k)} ) \big]  \nonumber \\
                             + \, \Delta t & \tilde{q}^I_{m+1,m+1} \big[ F_D( \phi_{AD}^{m+1,(k+1)} ) - F_D(  \phi^{m+1,(k)} ) + F_R( \phi^{m+1,(k+1)} ) - F_R(  \phi^{m+1,(k)} )\big] \nonumber \\
                             + \, \Delta t & Q^{0:m+1}\big( F(\phi^{(k)}) \big), \label{update_equation_misdcq_discrete_form_R}
\end{align}
where $Q^{0:m+1}\big( F(\phi^{(k)}) \big)$ approximates the integral of $F$ over the temporal interval $[t^0, t^{m+1}]$,
\begin{equation}
Q^{0:m+1}\big( F(\phi^{(k)}) \big) := \sum_{j = 0}^{M} q_{m+1,j} F( \phi^{(k)} ), 
\end{equation}
with
\begin{equation}
q_{m+1,j} := \frac{1}{\Delta t} \int_{t^0}^{t^{m+1}} L_j( \tau ) d \tau.
\label{weights}
\end{equation}
We denote by $\boldsymbol{Q} = \{ q_{i,j} \} \in \mathbb{R}^{M \times (M+1)}$ the matrix containing the weights defined in \ref{weights}. We decompose $\boldsymbol{Q}$ into its first column 
$\boldsymbol{q} \in \mathbb{R}^M$, and the matrix containing the remaining $M$ columns, denoted by $\boldsymbol{\tilde{Q}} \in \mathbb{R}^{M \times M}$.
The coefficients $\{ \tilde{q}^I_{i,j} \}$ in $\boldsymbol{\tilde{Q}}^I \in \mathbb{R}^{M \times M}$ are obtained by setting 
\begin{equation}
\boldsymbol{\tilde{Q}}^I := \boldsymbol{\tilde{U}}^T,
\end{equation}
where $\boldsymbol{\tilde{U}}^T$ is the transpose of the upper triangular matrix in the LU decomposition of $\boldsymbol{\tilde{Q}}^T$.
We refer to \cite{weiser2015faster} for a detailed discussion of the implications of this choice of coefficients on the convergence of 
the SDC iterations in the case of Gauss-Radau nodes. Finally, the coefficients $\{ \tilde{q}^E_{i,j} \}$ in 
$\boldsymbol{\tilde{Q}}^E \in \mathbb{R}^{M \times M}$ are set to $\tilde{q}^E_{m+1,m} := \Delta t^m / \Delta t$, and $\tilde{q}^E_{m+1,j \neq m+1} := 0$, 
corresponding to forward-Euler.

In the remainder of this paper, the modified MISDC scheme of \ref{update_equation_misdcq_discrete_form_D}-\ref{update_equation_misdcq_discrete_form_R} 
will be referred to as MISDCQ. In the convergence analysis of Section \oldref{subsection_convergence_analysis} and in numerical examples of 
Section \oldref{section_numerical_examples}, we will show that for stiff problems, MISDCQ converges to the fixed-point solution in fewer sweeps 
than the standard MISDC of \ref{update_equation_misdc_discrete_form_D}-\ref{update_equation_misdc_discrete_form_R}. Next, we use MISDCQ as a 
basis for the construction of a parallel-across-the-method Concurrent Implicit SDC scheme.

\subsection{
\label{subsection_parallel_multi_implicit_spectral_deferred_correction_scheme}
Concurrent implicit spectral deferred correction scheme
}

In the MISDC and MISDCQ schemes, the implicit solves corresponding to the diffusion and reaction steps represent most of the computational cost 
of the temporal integration scheme. At each SDC node, these solves are performed sequentially, with the reaction step following the diffusion step. These 
solves can become expensive for highly resolved combustion problems with a large number of species involved in a stiff and highly nonlinear chemical 
reaction path. In this section, we design a stable parallel-across-the-method Concurrent Implicit Spectral Deferred Correction (CISDC) scheme in which the 
implicit diffusion and reaction solves can be performed concurrently.

\subsubsection{Decoupling strategy}

In the MISDCQ scheme \ref{update_equation_misdcq_discrete_form_D}-\ref{update_equation_misdcq_discrete_form_R}, the diffusion step at node $m+1$ depends 
on the output of the reaction step at node $m$, denoted by $\phi^{m, (k+1)}$. Then, the output of the diffusion step at node $m+1$, denoted by 
$\phi^{m+1,(k+1)}_{AD}$, is used to compute the reaction step at node $m+1$, and obtain $\phi^{m+1,(k+1)}$. The MISDCQ algorithm is 
therefore inherently serial across the sweeps and cannot be parallelized without modification. 
We overcome this limitation by formally decoupling the diffusion and reaction steps to compute them concurrently.

There are multiple ways to achieve this goal. A simple strategy to decouple the diffusion and reaction steps consists in 
lagging the index $(k)$ in \ref{update_equation_misdcq_discrete_form_D}, i.e., to use a lagged diffusion term, $\phi_{AD}^{m+1,(k)}$, in the computation of
the reaction term, $\phi^{m+1,(k+1)}$. This decoupling allows a concurrent update of the diffusion and reaction terms at the same node, $m+1$. With this approach,
$M$ steps are needed to update the state variables at the $M$ SDC nodes. But, we found that his approach led to unstable schemes for stiff 
problems and therefore we do not further discuss it here.

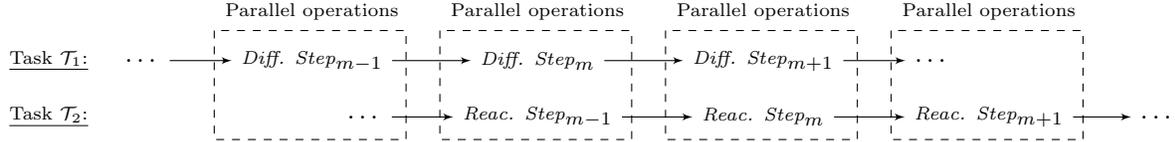
\begin{figure}[ht]
\centering
\tikzstyle{int}=[draw, minimum size=2em]
\tikzstyle{init} = [pin edge={to-,thick,black}]

\begin{tikzpicture}[node distance=3cm,auto,>=latex']
    \node (a_diff) {$\scriptsize \textit{Diff. Step}_{m-1}$};
    \node (b_diff) [left of=a_diff,node distance=2.25cm] {$\dots$};
    \node (c_diff) [right of=a_diff] {$\scriptsize \textit{Diff. Step}_{m}$};
    \node (d_diff) [left of=c_diff,node distance=1.5cm,coordinate] {};
    \node (f_diff) [right of=c_diff]{$\scriptsize \textit{Diff. Step}_{m+1}$};
    \node (e_diff) [left of=f_diff,node distance=1.5cm,coordinate] {};
    \node (g_diff) [right of=f_diff,node distance=2.25cm] {$\dots$};
    \node (label_diff) [left of=d_diff,node distance=5cm] {\scriptsize \underline{Task $\mathcal{T}_1$:}};
    \path[->] (b_diff) edge node {} (a_diff);  
    \path[->] (a_diff) edge node {} (c_diff);
    \draw[->] (c_diff) edge node {} (f_diff);
    \draw[->] (f_diff) edge node {} (g_diff);

    \node (a_reac) [below of=a_diff,node distance=0.75cm] {$ \quad \quad \qquad \dots$};
    \node (c_reac) [right of=a_reac] {$\scriptsize \textit{Reac. Step}_{m-1}$};
    \node (d_reac) [left of=c_reac,node distance=1.5cm,coordinate] {};
    \node (f_reac) [right of=c_reac]{$\scriptsize \textit{Reac. Step}_{m}$};
    \node (e_reac) [left of=f_reac,node distance=1.5cm,coordinate] {};
    \node (g_reac) [right of=f_reac] {$\scriptsize \textit{Reac. Step}_{m+1}$};
    \node (i_reac) [right of=g_reac,node distance=2.25cm] {$\dots$};
    \node (label_reac) [left of=d_reac,node distance=5cm] {\scriptsize \underline{Task $\mathcal{T}_2$:}};
    \path[->] (a_reac) edge node {} (c_reac);
    \draw[->] (c_reac) edge node {} (f_reac) ;
    \draw[->] (f_reac) edge node {} (g_reac);
    \draw[->] (g_reac) edge node {} (i_reac);

    \path (-1.3,0.4) node (d) {};
    \path (1.25,-1.05) node (e) {};
    \path [draw=black,dashed] (d) rectangle (e); 
    \node (ib1) at (0,0.65) {\scriptsize Parallel operations};

    \path (1.7,0.4) node (d) {};
    \path (4.25,-1.05) node (e) {};
    \path [draw=black,dashed] (d) rectangle (e); 
    \node (ib1) at (3,0.65) {\scriptsize Parallel operations};

    \path (4.7,0.4) node (d) {};
    \path (7.25,-1.05) node (e) {};
    \path [draw=black,dashed] (d) rectangle (e); 
    \node (ib1) at (6,0.65) {\scriptsize Parallel operations};

    \path (7.7,0.4) node (d) {};
    \path (10.25,-1.05) node (e) {};
    \path [draw=black,dashed] (d) rectangle (e); 
    \node (ib1) at (9,0.65) {\scriptsize Parallel operations};




\end{tikzpicture}
\vspace{-0.4cm}
\caption{\label{fig:sketch_amisdcQ}Sketch of the sweep parallelization strategy based on a decoupling of the diffusion step at $m+1$ from the reaction step at $m$ to solve them concurrently.}
\end{figure}

Instead, we adopt a decoupling strategy that allows a concurrent update of the diffusion and reaction terms at two consecutive SDC nodes.
Specifically, in the proposed marching scheme, we solve in parallel the diffusion step at node $m+1$ and the reaction step at the previous node, $m$. A sketch of this 
parallel marching scheme is in Fig.~\oldref{fig:sketch_amisdcQ}. We highlight that this algorithm requires $M+1$ steps to traverse the $M$ SDC nodes -- that is, 
one additional step compared to the approach outlined in the previous paragraph --,
but remains stable for the stiff problems considered in this work, as shown in Section \oldref{subsection_convergence_analysis}.


The parallel marching scheme described in Fig.~\oldref{fig:sketch_amisdcQ} is achieved by lagging in SDC sweep the terms that introduce a dependence on the reaction step at 
node $m$ in the diffusion step at node $m+1$. 
Specifically, in the diffusion solve at node $m+1$ of sweep $(k+1)$, we replace $F_A(  \phi^{m,(k+1)} )$, $F_D(  \phi^{m,(k+1)} )$, 
and $F_R(  \phi^{m,(k+1)} )$ with lagged values $F_A(  \phi^{m,(k+1,0)} )$, $F_D(  \phi^{m,(k+1,0)} )$, and $F_R(  \phi^{m,(k+1,0)} )$. 
Using these lagged values, the temporal interval is traversed to compute the corrected values $\phi_{AD}^{m+1,(k+1, 1)}$ and $\phi^{m,(k+1,1)}$ in 
parallel at all SDC nodes using the parallelization strategy sketched in Fig.~\oldref{fig:sketch_amisdcQ}. The accuracy of the approximation 
is improved iteratively by repeating this procedure, i.e., at iteration $\ell+1$, we correct $\phi_{AD}^{m+1,(k+1, \ell+1)}$ and $\phi^{m,(k+1,\ell+1)}$ 
in parallel at all nodes using lagged values computed at iteration $\ell$. This results in a loop on $\ell$ nested in sweep $(k+1)$. At node 
$m+1 \in \{ 1, \dots, M+1\}$, each correction entails solving the update equations
\begin{align}
\phi_{AD}^{m+1,(k+1, \ell+1)} = \phi^{0} + \Delta t    & \sum_{j = 1}^{m-1} \tilde{q}^E_{m+1,j}  \big[ F_A(  \phi^{j,(k+1,\ell+1)} ) - F_A( \phi^{j,(k)} ) \big]   \nonumber \\ 
                             + \, \Delta t & \sum_{j = 1}^{m-1} \tilde{q}^I_{m+1,j} \big[ F_D(  \phi^{j,(k+1,\ell+1)} ) - F_D( \phi^{j,(k)} )  + F_R(  \phi^{j,(k+1,\ell+1)} ) - F_R( \phi^{j,(k)} ) \big]  \nonumber \\
                             + \, \Delta t & \tilde{q}^E_{m+1,m} \big[ F_A(  \phi^{m,(k+1,\ell)} ) - F_A( \phi^{m,(k)} ) \big]   \nonumber \\ 
                             + \, \Delta t & \tilde{q}^I_{m+1,m} \big[ F_D(  \phi^{m,(k+1,\ell)} ) - F_D( \phi^{m,(k)} )  + F_R(  \phi^{m,(k+1,\ell)} ) - F_R( \phi^{m,(k)} ) \big]  \nonumber \\
                             + \, \Delta t & \tilde{q}^I_{m+1,m+1} \big[ F_D(  \phi_{AD}^{m+1,(k+1,\ell+1)} ) - F_D(  \phi^{m+1,(k)} ) + F_R( \phi^{m+1,(k+1,\ell)} ) - F_R(  \phi^{m+1,(k)} ) \big] \nonumber \\
                             + \, \Delta t & Q^{0:m+1}\big( F(\phi^{(k)}) \big), \label{update_equation_amisdcq_discrete_form_D} 
\end{align}
\begin{align}
\phi^{m+1,(k+1,\ell+1)}     = \phi^{0} + \Delta t    & \sum_{j = 1}^{m-1} \tilde{q}^E_{m+1,j}  \big[ F_A(  \phi^{j,(k+1,\ell+1)} ) - F_A( \phi^{j,(k)} ) \big]   \nonumber \\ 
                             + \, \Delta t & \sum_{j = 1}^{m-1} \tilde{q}^I_{m+1,j} \big[ F_D(  \phi^{j,(k+1,\ell+1)} ) - F_D( \phi^{j,(k)} )  + F_R(  \phi^{j,(k+1,\ell+1)} ) - F_R( \phi^{j,(k)} ) \big]  \nonumber \\
                             + \, \Delta t & \tilde{q}^E_{m+1,m} \big[ F_A(  \phi^{m,(k+1,\ell)} ) - F_A( \phi^{m,(k)} ) \big]   \nonumber \\ 
                             + \, \Delta t & \tilde{q}^I_{m+1,m} \big[ F_D(  \phi^{m,(k+1,\ell)} ) - F_D( \phi^{m,(k)} )  + F_R(  \phi^{m,(k+1,\ell+1)} ) - F_R( \phi^{m,(k)} ) \big]  \nonumber \\
                             + \, \Delta t & \tilde{q}^I_{m+1,m+1}              \big[ F_D(  \phi_{AD}^{m+1,(k+1,\ell+1)} ) - F_D(  \phi^{m+1,(k)} ) + F_R( \phi^{m+1,(k+1,\ell+1)} ) - F_R(  \phi^{m+1,(k)} ) \big] \nonumber \\
                             + \, \Delta t & Q^{0:m+1}\big( F(\phi^{(k)}) \big), \label{update_equation_amisdcq_discrete_form_R}
\end{align}
where the matrices $\boldsymbol{\tilde{Q}} = \{ q_{ij} \}$, $\boldsymbol{\tilde{Q}}^I = \{ \tilde{q}^I_{ij} \}$, and 
$\boldsymbol{\tilde{Q}}^E = \{ \tilde{q}^E_{ij} \}$, are the same as in the MISDCQ method of Section 
\oldref{subsection_multi_implicit_spectral_deferred_correction_schemes}. 
In \ref{update_equation_amisdcq_discrete_form_D} (respectively, \ref{update_equation_amisdcq_discrete_form_R}), the lagged terms -- i.e., the
terms evaluated at the previous iteration $\ell$ --, are in the third, fourth, and fifth lines (respectively, third and fourth lines). A procedure 
is required to initialize these lagged terms at the first iteration. The lagged advection and diffusion terms are 
initialized using the most recent advection-diffusion update, whereas the lagged reaction term is initialized using the solution at the previous 
sweep, that is,
\begin{align}
F_A(  \phi^{m,(k+1,0)} ) &:= F_A(  \phi_{AD}^{m,(k+1,0)} ), \label{initial_lagged_values_A} \\
F_D(  \phi^{m,(k+1,0)} ) &:= F_D(  \phi_{AD}^{m,(k+1,0)} ), \label{initial_lagged_values_D} \\
F_R(  \phi^{j,(k+1,0)} ) &:= F_R(  \phi^{j,(k)} ) \qquad j \in \{m, m+1\}.        \label{initial_lagged_values_R}
\end{align}
In the remainder of this paper, the scheme defined by 
\ref{update_equation_amisdcq_discrete_form_D} to \ref{initial_lagged_values_R} will be referred to as the Concurrent Implicit SDCQ
(CISDCQ) scheme, since the implicit diffusion and reaction steps can be performed in parallel.

We highlight that unlike in MISDCQ, the CISDCQ diffusion step \ref{update_equation_amisdcq_discrete_form_D} at node $m+1$ contains a reaction 
correction. This correction is based on known values computed at the same iteration $\phi^{j,(k+1,\ell+1)}$ ($j \in \{1, \dots, m-1\}$) in the second 
line of \ref{update_equation_amisdcq_discrete_form_D}, and on lagged values $\phi^{j, (k+1,\ell)}$ ($j \in \{m, m+1\}$) in the fourth and fifth lines of 
\ref{update_equation_amisdcq_discrete_form_D}. Therefore, if the nested iteration on $\ell$ converges, computing $\phi^{m+1, (k+1, \ell+1)}$ in 
\ref{update_equation_amisdcq_discrete_form_D} yields an approximation of the solution obtained from an implicit fully coupled diffusion-reaction solve. 
The reaction step can then be seen as a correction of the error incurred by the lagged terms in the diffusion step. This becomes apparent when 
\ref{update_equation_amisdcq_discrete_form_D} is used to write \ref{update_equation_amisdcq_discrete_form_R} in the equivalent form
\begin{align}
\phi^{m+1,(k+1,\ell+1)}     & = \phi_{AD}^{m+1,(k+1,\ell+1)} \nonumber \\
                        &+ \, \Delta t  \tilde{q}^I_{m+1,m}   \big[ F_R(  \phi^{m,(k+1,\ell+1)} ) - F_R( \phi^{m,(k+1,\ell)} ) \big] \nonumber \\
                        &+ \, \Delta t  \tilde{q}^I_{m+1,m+1} \big[ F_R( \phi^{m+1,(k+1,\ell+1)} ) - F_R(  \phi^{m+1,(k+1,\ell)} ) \big].
\label{update_equation_amisdcq_discrete_form_R_alt}
\end{align}
Using a linear convergence analysis in Section \oldref{subsection_convergence_analysis}, we will demonstrate that this modification of the MISDCQ 
algorithm yields a stable integration scheme whose sweeps can converge to the fixed-point solution faster than the original MISDCQ.

\subsubsection{\label{subsection_parallelization_strategy}Pipelining}

The nested iteration scheme on $\ell$ used to improve the accuracy of the lagged terms in 
\ref{update_equation_amisdcq_discrete_form_D}-\ref{update_equation_amisdcq_discrete_form_R} increases the computational cost of a sweep
since multiple iterations may be used at each sweep. However, this nested loop on $\ell$ can be performed efficiently, in parallel, when 
multiple iterations on $\ell$ are employed. This additional degree of parallelism, similar to the technique proposed in 
\cite{christlieb2010parallel,christlieb2012parallel}, exploits the structure of the update equations 
\ref{update_equation_amisdcq_discrete_form_D}-\ref{update_equation_amisdcq_discrete_form_R}. 

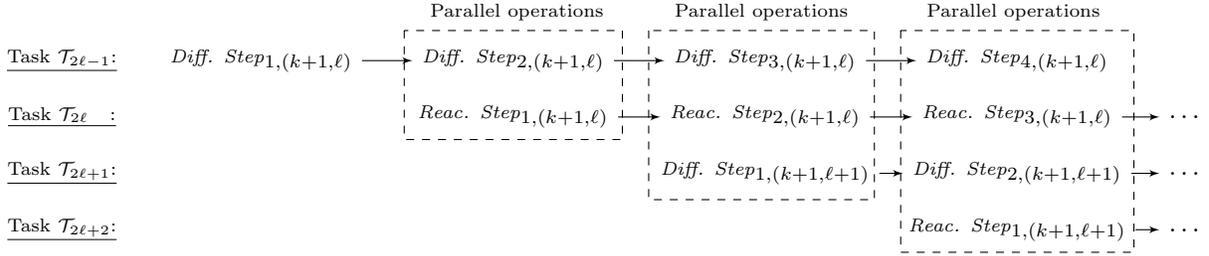
\begin{figure}[ht]
\centering
\tikzstyle{int}=[draw, minimum size=2em]
\tikzstyle{init} = [pin edge={to-,thick,black}]

\begin{tikzpicture}[node distance=3.35cm,auto,>=latex']
    \node (a_diff) {$\scriptsize \textit{Diff. Step}_{1,(k+1,\ell)}$};
    \node (c_diff) [right of=a_diff] {$\scriptsize \textit{Diff. Step}_{2,(k+1,\ell)}$};
    \node (d_diff) [left of=c_diff,node distance=1.5cm,coordinate] {};
    \node (f_diff) [right of=c_diff]{$\scriptsize \textit{Diff. Step}_{3,(k+1,\ell)}$};
    \node (e_diff) [left of=f_diff,node distance=1.5cm,coordinate] {};
    \node (g_diff) [right of=f_diff] {$\scriptsize \textit{Diff. Step}_{4,(k+1,\ell)}$};
    \node (i_diff) [right of=g_diff,node distance=2.25cm] {};
    \node (label_diff) [left of=d_diff,node distance=4.5cm] {\scriptsize \underline{Task $\mathcal{T}_{2 \ell -1}$:}};
    \path[->] (a_diff) edge node {} (c_diff);
    \draw[->] (c_diff) edge node {} (f_diff);
    \draw[->] (f_diff) edge node {} (g_diff);

    \node (a_reac) [below of=a_diff,node distance=0.75cm] {};
    \node (c_reac) [right of=a_reac] {$\scriptsize \textit{Reac. Step}_{1,(k+1,\ell)}$};
    \node (d_reac) [left of=c_reac,node distance=1.5cm,coordinate] {};
    \node (f_reac) [right of=c_reac]{$\scriptsize \textit{Reac. Step}_{2,(k+1,\ell)}$};
    \node (e_reac) [left of=f_reac,node distance=1.5cm,coordinate] {};
    \node (g_reac) [right of=f_reac] {$\scriptsize \textit{Reac. Step}_{3,(k+1,\ell)}$};
    \node (i_reac) [right of=g_reac,node distance=2.25cm] {$\dots$};
    \node (label_reac) [left of=d_reac,node distance=4.5cm] {\scriptsize \underline{Task $\mathcal{T}_{2 \ell~~~}$:}};
    \draw[->] (c_reac) edge node {} (f_reac) ;
    \draw[->] (f_reac) edge node {} (g_reac);
    \draw[->] (g_reac) edge node {} (i_reac);

    \node (a_diff_2) [below of=a_reac,node distance=0.75cm] {};
    \node (c_diff_2) [right of=a_diff_2] {\textcolor{white}{$\scriptsize \textcolor{white}{\textit{Reac. Step}_{1,(k+1,\ell)}}$}};
    \node (d_diff_2) [left of=c_diff_2,node distance=1.5cm,coordinate] {};
    \node (f_diff_2) [right of=c_diff_2]{$\scriptsize \textit{Diff. Step}_{1,(k+1,\ell+1)}$};
    \node (e_diff_2) [left of=f_diff_2,node distance=1.5cm,coordinate] {};
    \node (g_diff_2) [right of=f_diff_2] {$\scriptsize \textit{Diff. Step}_{2,(k+1,\ell+1)}$};
    \node (i_diff_2) [right of=g_diff_2,node distance=2.25cm] {$\dots$};
    \node (label_diff_2) [left of=d_diff_2,node distance=4.5cm] {\scriptsize \underline{Task $\mathcal{T}_{2 \ell +1 }$:}};
    \draw[->] (f_diff_2) edge node {} (g_diff_2);
    \draw[->] (g_diff_2) edge node {} (i_diff_2);

    \node (a_reac_2) [below of=a_diff_2,node distance=0.75cm] {};
    \node (c_reac_2) [right of=a_reac_2] {$\scriptsize \textcolor{white}{\textit{Reac. Step}_{1,(k+1,\ell)}}$};
    \node (d_reac_2) [left of=c_reac_2,node distance=1.5cm,coordinate] {};
    \node (f_reac_2) [right of=c_reac_2]{$\scriptsize \textcolor{white}{\textit{Reac. Step}_{2,(k+1,\ell+1)}}$};
    \node (e_reac_2) [left of=f_reac_2,node distance=1.5cm,coordinate] {};
    \node (g_reac_2) [right of=f_reac_2] {$\scriptsize \textit{Reac. Step}_{1,(k+1,\ell+1)}$};
    \node (i_reac_2) [right of=g_reac_2,node distance=2.25cm] {$\dots$};
    \node (label_reac_2) [left of=d_reac_2,node distance=4.5cm] {\scriptsize \underline{Task $\mathcal{T}_{2 \ell + 2}$:}};
    \draw[->] (g_reac_2) edge node {} (i_reac_2);

    \path (1.9,0.4) node (d) {};
    \path (4.8,-1.05) node (e) {};
    \path [draw=black,dashed] (d) rectangle (e); 
    \node (ib1) at (3.4,0.65) {\scriptsize Parallel operations};

    \path (5.15,0.4) node (d) {};
    \path (8.15,-1.85) node (e) {};
    \path [draw=black,dashed] (d) rectangle (e); 
    \node (ib1) at (6.65,0.65) {\scriptsize Parallel operations};

    \path (8.5,0.4) node (d) {};
    \path (11.6,-2.55) node (e) {};
    \path [draw=black,dashed] (d) rectangle (e); 
    \node (ib1) at (10,0.65) {\scriptsize Parallel operations};




\end{tikzpicture}
\vspace{-0.4cm}
\caption{\label{fig:sketch_amisdcQ_pipelined}Sketch of the pipelined CISDCQ sweep with $2 \nu$ tasks for five Gauss-Lobatto nodes. 
The diffusion step $m+1,(k+1,\ell+1)$ entails solving \ref{update_equation_amisdcq_discrete_form_D} for $\phi^{m+1,(k+1,\ell+1)}_{AD}$, 
and the reaction step $m+1,(k+1,\ell+1)$ requires solving \ref{update_equation_amisdcq_discrete_form_R} for $\phi^{m+1,(k+1,\ell+1)}$.}
\end{figure}

We assume that at least two nested iterations on $\ell$ are used. As in Fig.~\oldref{fig:sketch_amisdcQ}, task $\mathcal{T}_1$, which 
corresponds to the diffusion steps necessary to compute $\phi^{m, (k+1,1)}_{AD}$ ($m \in \{1, \dots, M\}$), and task $\mathcal{T}_2$, 
which corresponds to the reaction steps needed to compute $\phi^{m, (k+1,1)}$ ($m \in \{1, \dots, M\}$), can be executed in parallel. 
But, as soon as $\phi^{1,(k+1,1)}$ has been computed, task $\mathcal{T}_3$ can be launched to perform the diffusion steps at the next 
iteration ($\ell = 2$) and obtain $\phi_{AD}^{m, (k+1, 2)}$ ($m \in \{1, \dots, M\}$) using the lagged values updated by task $\mathcal{T}_2$. 
Then, once $\phi^{1, (k+1, 2)}_{AD}$ is available, task $\mathcal{T}_4$ can be executed in parallel to compute the reaction updates 
$\phi^{m, (k+1, 2)}$. This approach can be generalized to the case of $2 \nu$ tasks with $\nu \geq 2$. A sketch of this pipelined 
nested loop on $\ell$ is in Fig.~\oldref{fig:sketch_amisdcQ_pipelined}. 

%

The parallelism inherent in pipelining the substeps is assumed to be in addition to any spatial parallelization.  
Hence the term {\it processor} used here denotes a group of cores performing the parallel tasks in the CISDCQ algorithm.
The number of tasks does not necessarily correspond to the number of processors used for the implementation of the algorithm
as a given processor can be reused for multiple tasks.
If three Gauss-Lobatto nodes are used ($M = 2$), only two diffusion steps and two reaction steps have to be performed at each nested 
iteration on $\ell$. This means that when task $\mathcal{T}_{3}$ is launched, the two diffusion steps involved in task $\mathcal{T}_1$ 
have already been performed by the first processor which is now idling. Therefore, the first processor can be reused 
to perform task $\mathcal{T}_{3}$. For the same reason, the second processor is idling after completing task $\mathcal{T}_2$ and
can be reused to perform task $\mathcal{T}_4$. With this strategy, the first processor performs tasks $\mathcal{T}_1$, $\mathcal{T}_3$, 
$\dots$, $\mathcal{T}_{2(\nu-1)-1}$ while the second processor performs tasks $\mathcal{T}_2$, $\mathcal{T}_4$, $\dots$, $\mathcal{T}_{2(\nu-1)}$. 
Similarly, if five Gauss-Lobatto nodes are used ($M = 4$), then at most four processors can be used. In the general case, the $2 \nu$
pipelined tasks involved in CISDCQ-$\nu$ use up to $\max( 2\nu, M)$ processors working in parallel.

To assess the computational cost reduction generated by the CISDCQ algorithm, we assume that the cost of a diffusion step 
(respectively, a reaction step) is a constant $\Upsilon_{AD}$ (respectively, $\Upsilon_R$). $N_{\textit{M}}$ (respectively, $N_{\textit{C}}$) 
denotes the number of sweeps required to achieve convergence with MISDCQ (respectively, CISDCQ). For a time interval $[t^n, t^{n+1}]$ 
decomposed into $M$ subintervals using $M+1$ temporal nodes, the computational cost of MISDCQ, $C_{\textit{M}}$, is given by
\begin{equation}
C_{\textit{M}} = N_{\textit{M}} M ( \Upsilon_{AD} + \Upsilon_R ).
\end{equation}
In the evaluation of the computational cost of CISDCQ, $C_{\textit{C}}$, we neglect the communication costs and assume that $\nu$ iterations on $\ell$ 
are nested in each sweep to improve the accuracy of the lagged terms. First, we compute the cost of iteration $\ell +1 $ in sweep $(k+1)$, i.e., 
the cost of the operations performed in tasks $\mathcal{T}_{2 \ell +1}$ and $\mathcal{T}_{2 \ell +2}$ (see Fig.~\oldref{fig:sketch_amisdcQ_pipelined}). 
These operations include the initial diffusion step to compute $\phi_{AD}^{1,(k+1,\ell+1)}$, followed by $(M-1)$ parallel diffusion and reaction steps, 
and by the final reaction step to compute $\phi^{M, (k+1, \ell+1)}$. Therefore, the cost of an iteration on $\ell$ nested in sweep $(k+1)$ is 
\begin{equation}
\Upsilon_{AD} + (M-1) \max(\Upsilon_{AD},\Upsilon_R) + \Upsilon_R.
\label{cost_of_amisdcq_nested_iteration}
\end{equation}
Task $\mathcal{T}_{2 \ell+1}$ can only start after task $\mathcal{T}_{2 \ell - 1}$ has computed $\phi_{AD}^{2,(k+1,\ell)}$ and after task 
$\mathcal{T}_{2 \ell}$ has used this updated value to calculate $\phi^{2,(k+1,\ell)}$. The total cost of these two operations is 
$\Upsilon_{AD} + \Upsilon_R$. Using this result along with \ref{cost_of_amisdcq_nested_iteration}, we find by induction on $\ell$ that the cost of 
the CISDCQ sweep with $\nu$ nested iterations on $\ell$ is 
\begin{equation}
\nu (\Upsilon_{AD} + \Upsilon_R) + (M-1) \max(\Upsilon_{AD},\Upsilon_R).
\end{equation}
The ratio of the computational cost of MISDCQ over that of CISDCQ-$\nu$ -- i.e., CISDCQ with $\nu$ nested iterations
on $\ell$ --, denoted by $R_{\nu}$, is therefore given by
\begin{equation}
R_{\nu} := \frac{C_{\textit{M}}}{C_{\textit{C}}} 
= \frac{N_\textit{M} M ( \Upsilon_{AD} + \Upsilon_R ) }{N_\textit{A} \big( \nu (\Upsilon_{AD} + \Upsilon_R) + (M-1) \max(\Upsilon_{AD},\Upsilon_R) \big)}.
\end{equation}
The CISDCQ-$\nu$ algorithm reduces the computational cost compared to MISDCQ whenever $R_{\nu} > 1$. With the introduction 
of the ratio $\alpha$ as
\begin{equation}
\alpha := \frac{\Upsilon_{AD} + \Upsilon_R}{\max(\Upsilon_{AD}, \Upsilon_R)},
\end{equation}
the ratio $R_{\nu}$ comparing the respective computational costs of CISDCQ and MISDCQ simplifies to
\begin{equation}
R_{\nu} := \frac{N_{\textit{M}}}{N_{\textit{C}}} \times \frac{\alpha M}{\alpha \nu + M -1}.
\label{simple_parallel_speedup}
\end{equation}
We note that $\alpha = 2$ whenever $\Upsilon_{AD} = \Upsilon_R$. 
In \ref{simple_parallel_speedup}, the first ratio compares the convergence of the MISDCQ sweeps with that of the CISDCQ-$\nu$ sweeps, while the second 
ratio compares the computational cost of a single sweep for the two schemes. This decomposition illustrates the two mechanisms leading to a parallel 
speedup. The first mechanism results from the modification of the update equations 
\ref{update_equation_amisdcq_discrete_form_D}-\ref{update_equation_amisdcq_discrete_form_R} with the introduction of a reaction correction 
in the diffusion step. It consists in achieving convergence to the fixed-point solution in fewer sweeps with CISDCQ-$\nu$ to increase the ratio $N_{\textit{M}}/N_{\textit{C}}$.
This enhanced convergence rate can be obtained by increasing the number of iterations on $\ell$ in the nested loop to improve the accuracy of the lagged 
terms of \ref{update_equation_amisdcq_discrete_form_D}. But, performing more nested iterations is in contradiction to the second mechanism that can be 
leveraged to reduce the computational cost, which is the reduction of the cost of a single sweep thanks to parallelization
to increase the ratio $\alpha M / (\alpha \nu + M -1)$. This 
trade-off will be explored with numerical examples in Section \oldref{section_numerical_examples}.

For completeness, we also compare the computational cost of CISDCQ-$\nu$ to that of an SDC scheme in which advection is treated explicitly, 
while diffusion and reaction are treated implicitly in a fully coupled fashion. This scheme, referred to as IMEXQ, is a variant of the IMEX scheme introduced 
in \cite{minion2003semi}, and is based on the LU factorization of the integration matrix explained in Section 
\oldref{subsection_multi_implicit_spectral_deferred_correction_schemes}. 
The ratio of the computational cost of IMEXQ over that of CISDCQ-$\nu$ reads
\begin{equation}
R'_{\nu} := \frac{N_{\textit{I}}}{N_{\textit{C}}} \times \frac{\beta M}{\beta \nu + M -1},
\label{simple_parallel_speedup_2}
\end{equation}
where the coefficient $\beta$ is defined as
\begin{equation}
\beta := \frac{\Upsilon_{ADR}}{\max(\Upsilon_{AD}, \Upsilon_R)}.
\end{equation}
where $N_{\textit{I}}$ is the number of sweeps needed by IMEXQ to converge to the fixed-point solution, and $\Upsilon_{ADR}$ is the cost 
of solving the fully coupled advection-diffusion-reaction system. In Section \oldref{subsection_convergence_analysis}, we show that the IMEXQ sweeps 
converge faster than the CISDCQ-$\nu$ and MISDCQ sweeps for stiff problems. But, one expects $\beta$ to be very large for multidimensional problems with 
stiff chemistry, since solving the global nonlinear system is expensive and requires efficient physics-based preconditioners for the linear systems.

\subsection{\label{section_numerical_methodology}Application to low-Mach number combustion}

We now describe the application of the CISDCQ time integration scheme to the low-Mach number equation set presented in Section 
\oldref{section_governing_equations}. The governing equations are discretized in space using a finite-volume 
formulation with uniform grid spacing. The fourth-order discretization in space is the same as in \cite{pazner2016high}, and relies on the 
operators found in the finite-volume literature \citep{mccorquodale2011high,zhang2012fourth}. 
We also use the same volume discrepancy method based on a correction to the divergence constraint. This correction is needed 
because even if the initial state satisfies the EOS, the variables are updated with fluxes that vary linearly over the
time step. Due to the nonlinearity of the EOS, there is no guarantee that each component will evolve in a way that the 
new state will also satisfy the EOS. The purpose of the $\delta_{\chi}$-correction to the divergence constraint is to 
adjust the face velocities so that a conservative mass and enthalpy update will end up satisfying the EOS.


For the description of the CISDCQ integration method, we assume that the integration is based on $M+1$ Gauss-Lobatto nodes, $N_C$ sweeps, and $\nu$ 
iterations on $\ell$ nested in each sweep. The algorithm to advance the solution from $t^n$ to $t^{n+1}$ is described below. \\

\vspace{0.2cm}

\noindent \textbf{Initialization}

\begin{enumerate} \itemsep 8pt 

\item[\textbf{I1.}] Set $(\rho h, \, \rho \boldsymbol{Y})^{0, (k)} := (\rho h, \rho \boldsymbol{Y})^n$ for all $k \in \{0, \dots, N_C\}$, i.e. the solution at temporal 
node $m = 0$ is a copy of the solution at $t^n$ for all CISDCQ sweeps.

\item[\textbf{I2.}] Set $(\rho h, \, \rho \boldsymbol{Y})^{m, (0)} := (\rho h, \rho \boldsymbol{Y})^n$ for all $m \in \{1, \dots, M\}$, i.e. the solution for 
$k = 0$ is a copy of the solution at $t^n$ for all temporal nodes.

\item[\textbf{I3.}] Define a divergence constraint correction for each temporal interval, $\delta^{m-1:m,(k)}_{\chi}$, as in \cite{pazner2016high}. This 
correction will be applied to the divergence constraint at each node to ensure consistency with the EOS. Initialize the correction
to $\delta^{m-1:m,(0)}_{\chi} := 0$ for all $m \in \{1, \dots, M\}$. 

\item[\textbf{I4.}] Compute face-averaged velocities at $t^n$ by solving the divergence constraint \ref{divergence_constraint}
\begin{equation}
\nabla \cdot U^n = S^n,
\end{equation}
which in one dimension can be done by writing $U^n_{i+1} = U^n_i + \int^{x_{i+1}}_{x_i} S$. Then, set $U^{0,(k)} := U^n$ for all $k \in \{ 0, \dots, N_C \}$.
In addition, evaluate the right-hand side of the discretized species and enthalpy equations obtained from \ref{species_equation} and \ref{enthalpy_equation}, respectively. 
These terms will be used to evaluate $Q^{0:m+1}[ F(\phi^{(0)}) ]$ in the first CISDCQ sweep. 

\end{enumerate}

\noindent \textbf{Sweeps}

\begin{enumerate} \itemsep 8pt

\item[\textbf{for}] $k = 0$ \textbf{to} $N_C-1$ \textbf{do} \\

\begin{enumerate} \itemsep 8pt

\item[\textbf{S1.}] Set $\boldsymbol{\dot{\omega}}^{m, (k+1, 0)} := \boldsymbol{\dot{\omega}}^{m,(k)}$ for all $m \in \{0, \dots, M\}$, 
i.e. the reaction term  at iteration $\ell = 0$ at sweep $(k+1)$ is a copy of the state of the system at the end of sweep $(k)$. 
The vector of production terms is
$\boldsymbol{\dot{\omega}}^{m,(k+1,0)} = [\dot{\omega}^{m,(k+1,0)}_1, \dots, \dot{\omega}^{m,(k+1,0)}_N]^T$.

\item[\textbf{for}] $\ell = 0$ \textbf{to} $\nu-1$ \textbf{do}

\begin{enumerate} \itemsep 8pt

\item[\textbf{for}] $m = 0$ \textbf{to} $M - 1$ \textbf{do} \\

\begin{enumerate} \itemsep 8pt 


\item[\textbf{S2.}] Update the density $\rho^{m+1,(k+1,\ell+1)}$ explicitly by applying the CISDCQ correction \ref{update_equation_amisdcq_discrete_form_D} 
to the discretized continuity equation 
obtained from \ref{continuity_equation}. Since \ref{continuity_equation} only contains advection terms, \ref{update_equation_amisdcq_discrete_form_D} 
simplifies to
\begin{align}
\rho^{m+1,(k+1,\ell+1)} &= \rho^0 \nonumber \\ &+ \Delta t \sum^{m-1}_{j=1} \tilde{q}^E_{m+1,j} \big[ - \nabla \cdot (U \rho)^{j,(k+1,\ell+1)} + \nabla \cdot (U \rho)^{j,(k)} \big] \nonumber \\
                    &+ \Delta t                \tilde{q}^E_{m+1,m} \big[ - \nabla \cdot (U \rho)^{m,(k+1,\ell)} + \nabla \cdot (U \rho)^{m,(k)} \big] \nonumber \\
                    &+ \Delta t                 Q^{0:m+1}  \big[ - \nabla \cdot (U \rho)^{(k)} \big],
\label{mass_fraction_update}
\end{align}
where the advection flux at substep $m$ and sweep $(k+1)$ in the third line is lagged in iteration on $\ell$.

\item[\textbf{S3.}] Compute the new mass fractions $Y^{m+1,(k+1,\ell+1)}_{i,AD}$ by applying the CISDCQ correction 
\ref{update_equation_amisdcq_discrete_form_D} to the discretized species equation
obtained from \ref{species_equation}. We write the implicit banded linear system solved with a direct solver during this step as
\begin{align}
(\rho Y_{i,AD} )^{m+1,(k+1,\ell+1)} &= (\rho Y_i)^0 + A_{\textit{mass}} + D_{\textit{mass}} + R_{\textit{mass}} \nonumber \\ 
                              & + \Delta t                Q^{0:m+1}  \big[ - \nabla \cdot ( U \rho Y_i )^{(k)} + \nabla \cdot \Gamma^{(k)}_i + \dot{\omega}^{(k)}_i \big].
\end{align}
In the explicit advection correction, $A_{\textit{mass}}$, the advection mass flux at substep $m$ and sweep $(k+1)$ is lagged in iteration on $\ell$ as follows:
\begin{align}
A_{\textit{mass}} &=  \Delta t \sum^{m-1}_{j=1} \tilde{q}^E_{m+1,j}  \big[ - \nabla \cdot (U \rho Y_i)^{j,(k+1,\ell+1)} + \nabla \cdot (U \rho Y_i)^{j,(k)} \big]\nonumber \\
               &+  \Delta t                \tilde{q}^E_{m+1,m} \quad \, \,  \, \, \, \big[ - \nabla \cdot (U \rho Y_i)^{m,(k+1,\ell)} + \nabla \cdot (U \rho Y_i)^{m,(k)} \big].
\end{align}
Similarly, in the implicit diffusion correction, $D_{\textit{mass}}$, the diffusion mass flux at substep $m$ and sweep $(k+1)$ is lagged in iteration on $\ell$,
\begin{align}
D_{\textit{mass}} &=  \Delta t \sum^{m-1}_{j=1} \tilde{q}^I_{m+1,j}  \big[  \nabla \cdot \Gamma^{j,(k+1,\ell+1)}_i - \nabla \cdot \Gamma^{j,(k)}_i
                                                                    \big] \nonumber \\
               &+ \Delta t                \tilde{q}^I_{m+1,m} \quad \, \,  \, \, \,  \big[ \nabla \cdot \Gamma^{m,(k+1,\ell)}_i - \nabla \cdot \Gamma^{m,(k)}_i 
                                                                                \big] \nonumber \\
               &+ \Delta t                \tilde{q}^I_{m+1,m+1} \, \, \, \, \, \big[ \nabla \cdot \Gamma^{m+1,(k+1,\ell+1)}_{i,AD} - \nabla \cdot \Gamma^{m+1,(k)}_i 
                                                                               \big], 
\label{diffusion_correction_mass}
\end{align}
where the lagged discrete diffusion flux at substep $m+1$ is
\begin{equation}
\Gamma^{m+1,(k+1,\ell+1)}_{i,AD} := 
\left\{
\begin{array}{l l}
    \rho^{m+1,(k)} \mathcal{D}^{m+1,(k)}_i \nabla Y^{m+1,(k+1,\ell+1)}_{i,AD} & \quad \text{if } \ell = 0 \\[3pt]
    \rho^{m+1,(k+1,\ell)} \mathcal{D}^{m+1,(k+1,\ell)}_i \nabla Y^{m+1,(k+1,\ell+1)}_{i,AD}& \quad \text{otherwise}.
\end{array} \right. \label{lagged_diffusion_flux}
\end{equation}
Finally, in the implicit reaction correction, $R_{\textit{mass}}$, the reaction terms at substeps $m$ and $m+1$ at sweep $(k+1)$ are lagged in iteration on $\ell$. The 
reaction correction is therefore defined as
\begin{align}
R_{\textit{mass}} &= \Delta t \sum^{m-1}_{j=1} \tilde{q}^I_{m+1,j}  \big[   \dot{\omega}^{j, (k+1,\ell+1)}_i - \dot{\omega}^{j,(k)}_i  \big] \nonumber \\
               &+ \Delta t                \tilde{q}^I_{m+1,m} \quad \, \,  \, \, \,  \big[  \dot{\omega}^{m, (k+1,\ell)}_i    - \dot{\omega}^{m,(k)}_i \big] \nonumber \\
               & + \Delta t                \tilde{q}^I_{m+1,m+1} \, \, \, \, \, \big[  \dot{\omega}^{m+1, (k+1,\ell)}_i    - \dot{\omega}^{m+1,(k)}_i \big].
\end{align}



\item[\textbf{S4.}] Compute the new enthalpy $h^{m+1,(k+1,\ell+1)}_{AD}$ by applying the CISDCQ correction \ref{update_equation_amisdcq_discrete_form_D} 
to the discretized energy equation 
obtained from \ref{enthalpy_equation}. 
This step involves using a direct solver to solve the implicit system
\begin{align}
&(\rho h_{AD} )^{m+1,(k+1,\ell+1)} = (\rho h)^0 + A_{\textit{energy}} + D^{\textit{diffdiff}}_{\textit{energy}} + D^{\textit{diff}}_{\textit{energy}}  \label{energy_update} \\
                               & + \Delta t  Q^{0:m+1}    \big[ - \nabla \cdot ( U \rho h )^{(k)} + \nabla \cdot \frac{\lambda^{(k)}}{c^{(k)}_p} \nabla h^{(k)} 
                                                                                                        + \sum_i \nabla \cdot h^{(k)}_i \big( \Gamma^{(k)}_i - \frac{\lambda^{(k)}}{c^{(k)}_p} \nabla Y^{(k)}_i \big) \big]. \nonumber
\end{align}
As in the mass fraction update \ref{mass_fraction_update}, the explicit advection piece in the previous equation, denoted by $A_{\textit{energy}}$, is based on
a lagged advection flux at substep $m$ and sweep $(k+1)$. Specifically, we write
\begin{align}
A_{\textit{energy}} &=  \Delta t \sum^{m-1}_{j=1} \tilde{q}^E_{m+1,j}  \big[ - \nabla \cdot (U \rho h)^{j,(k+1,\ell+1)} + \nabla \cdot (U \rho h)^{j,(k)} \big] \nonumber \\
& + \Delta t               \tilde{q}^E_{m+1,m} \quad \, \,  \, \, \, \, \big[ - \nabla \cdot (U \rho h)^{m,(k+1,\ell)} + \nabla \cdot (U \rho h)^{m,(k)} \big].
\end{align}
As in \cite{pazner2016high}, the differential diffusion terms corresponding to the sum in the right-hand side of 
\ref{enthalpy_equation} are evaluated 
explicitly to simplify the linear system that needs to be solved. Therefore, the differential diffusion correction, $D^{\textit{diffdiff}}_{\textit{energy}}$, is
\begin{align}
D^{\textit{diffdiff}}_{\textit{energy}} &= \Delta t \sum^{m-1}_{j=1} \tilde{q}^E_{m+1,j}  \big[ \sum_i \nabla \cdot h^{j,(k+1,\ell+1)}_i \big( \Gamma^{j,(k+1,\ell+1)}_i - \frac{\lambda^{j,(k+1,\ell+1)}}{c^{j,(k+1,\ell+1)}_p} \nabla Y^{j,(k+1,\ell+1)}_i \big) \nonumber \\ 
          & \qquad \qquad \qquad \qquad \quad \,  - \nabla \cdot h^{j,(k)}_i \big( \Gamma^{j,(k)}_i - \frac{\lambda^{j,(k)}}{c^{(k)}_p} \nabla Y^{j,(k)}_i \big) \big] \nonumber \\
&+ \Delta t               \tilde{q}^E_{m+1,m}   \quad \, \,  \, \, \, \big[ \sum_i \nabla \cdot h^{m,(k+1,\ell)}_i \big( \Gamma^{m,(k+1,\ell)}_i - \frac{\lambda^{m,(k+1,\ell)}}{c^{m,(k+1,\ell)}_p} \nabla Y^{m,(k+1,\ell)}_i \big) \nonumber \\
& \qquad \qquad \qquad \qquad \quad \,  - \nabla \cdot h^{m,(k)}_i \big( \Gamma^{m,(k)}_i - \frac{\lambda^{m,(k)}}{c^{(k)}_p} \nabla Y^{m,(k)}_i \big) \big]. 
\end{align}
Following the approach of \ref{diffusion_correction_mass}, the diffusion correction, $D^{\textit{diff}}_{\textit{energy}}$, contains a lagged diffusion flux at substep $m$ and sweep $(k+1)$
\begin{align}
D^{\textit{diff}}_{\textit{energy}} &= \Delta t \sum^{m-1}_{j=1} \tilde{q}^I_{m+1,j}  \big[ \nabla \cdot \frac{\lambda^{j,(k+1,\ell+1)}}{c^{j,(k+1,\ell+1)}_p} \nabla h^{j,(k+1,\ell+1)}
                                                                                                 - \nabla \cdot \frac{\lambda^{j,(k)}}{c^{j,(k)}_p} \nabla h^{j,(k)} \big] \\     
                             &\quad + \Delta t               \tilde{q}^I_{m+1,m}   \quad   \big[ \nabla \cdot \frac{\lambda^{m,(k+1,\ell)}}{c^{m,(k+1,\ell)}_p} \nabla h^{m,(k+1,\ell)}
                                                                                                   - \nabla \cdot \frac{\lambda^{m,(k)}}{c^{m,(k)}_p} \nabla h^{m,(k)} \big] \nonumber \\
                             &\quad + \Delta t               \tilde{q}^I_{m+1,m+1}   \big[ \nabla \cdot \big( \frac{\lambda}{c_p} \big)^{m+1,(k+1,\ell+1)}_{AD} \nabla h^{m+1,(k+1,\ell+1)}_{AD}
                                                                                                    - \nabla \cdot \frac{\lambda^{m+1,(k)}}{c^{m+1,(k)}_p} \nabla h^{m+1,(k)} \big], \nonumber
\end{align}
where the diffusion coefficient at substep $m+1$, denoted by $\big( \lambda / c_p  \big)^{m+1,(k+1,\ell+1)}_{AD}$, is lagged using the same method as in \ref{lagged_diffusion_flux}.
After the enthalpy update, we can set 
\begin{equation}
(\rho h)^{m+1,(k+1,\ell+1)} = \rho^{m+1,(k+1,\ell+1)} h^{m+1,(k+1,\ell+1)}_{AD},
\end{equation}
since there is no contribution due to reactions in the enthalpy update.
If $\ell = 0$, recompute the diffusion terms in the right-hand side of the discrete species and energy equations. Following \ref{initial_lagged_values_A}-\ref{initial_lagged_values_D}, 
these terms will be used to define the quantities evaluated at $m+1, (k+1,0)$ at the next iteration $\ell = 1$.

\item[\textbf{S5.}] Compute the new mass fractions $Y^{m+1,(k+1,\ell+1)}_{i}$ by applying the CISDCQ correction \ref{update_equation_amisdcq_discrete_form_R_alt} to the discretized species equation
\begin{align}
(\rho Y_i)^{m+1, (k+1,\ell+1)} &= \rho^{m+1,(k+1,\ell+1)} Y^{m+1,(k+1,\ell+1)}_{i,AD} \nonumber \\
                          &+ \Delta t \tilde{q}^I_{m+1,m} \, \, \, \, \, \, \big[ \dot{\omega}^{m, (k+1,\ell+1)}_i - \dot{\omega}^{m,(k+1,\ell)}_i \big] \nonumber \\
                          &+ \Delta t \tilde{q}^I_{m+1,m+1} \big[ \dot{\omega}^{m+1, (k+1,\ell+1)}_i - \dot{\omega}^{m+1,(k+1,\ell)}_i \big]. 
\label{reaction_step_numerical_methodology}
\end{align}
To obtain the new mass fractions from \ref{reaction_step_numerical_methodology}, we form the backward-Euler type nonlinear system
\begin{equation}
\rho^{m+1,(k+1,\ell+1)} \boldsymbol{Y} - \Delta t \tilde{q}^I_{m+1,m+1} \boldsymbol{\dot{\omega}}^{m+1,(k+1,\ell+1)}( \boldsymbol{Y} ) = \boldsymbol{b},
\end{equation}
where we defined a right-hand side $\boldsymbol{b}$ obtained from \ref{reaction_step_numerical_methodology}. 
The density and the 
enthalpy have been computed in the previous steps of the algorithm and are known. The system is solved with Newton's method in which the initial 
guess is the solution at the previous nested iteration on $\ell$. The Jacobian matrix is computed analytically \citep{perini2012analytical}.
\item[\textbf{S6.}] Increment the divergence constraint correction by setting
\begin{equation}
\delta_{\chi}^{m:m+1,(k+1,\ell+1)} := \delta^{m:m+1,(k)}_{\chi} + \frac{2}{p_0} \big( \frac{p^{m+1,(k+1,\ell+1)}_{\textit{EOS}} - p_0}{\Delta t^{m}} \big),
\end{equation}
where the thermodynamic pressure, $p_{\textit{EOS}}$, is defined in \ref{p_eos_definition}.
Then, solve the divergence constraint \ref{divergence_constraint}
\begin{equation}
\nabla \cdot U^{m+1,(k+1,\ell+1)} = S^{m+1,(k+1,\ell+1)} + \delta^{m:m+1,(k+1,\ell+1)}_{\chi}.
\end{equation}
Finally, recompute the diffusion terms in the right-hand side of the species and enthalpy equations \ref{species_equation}-\ref{enthalpy_equation}.

\end{enumerate}

\item[\textbf{end}] \textbf{for} (end loop over temporal nodes $m$)

\end{enumerate}

\item[\textbf{end}] \textbf{for} (end loop over nested iterations $\ell$)

\item[\textbf{S7.}] Set $(\rho \boldsymbol{Y}, \rho h)^{m,(k+1)} := (\rho \boldsymbol{Y}, \rho h)^{m,(k+1,\nu)}$, i.e., we save the solution at the last iteration $\nu$ of sweep $(k+1)$.

\end{enumerate}

\item[\textbf{end}] \textbf{for} (end loop over CISDCQ sweeps $k$)

\end{enumerate}

We highlight that by construction, steps \textbf{S2}, \textbf{S3}, and \textbf{S4} at node $m+1, (k+1,\ell+1)$ can be executed without 
knowing the state of the system at $m, (k+1,\ell+1)$ thanks to the decoupling method based on lagging (see Section 
\oldref{subsection_parallel_multi_implicit_spectral_deferred_correction_scheme}). Therefore, if enough processors are available, a given 
processor executes task $\mathcal{T}_1$, i.e., it performs the advection and diffusion steps, \textbf{S2}, \textbf{S3}, and \textbf{S4}, 
at $m+1, (k+1,1)$ while another processor computes task $\mathcal{T}_2$ corresponding to the reaction step and right-hand side update 
(\textbf{S5} and \textbf{S6}) at $m, (k+1,1)$. If more iterations on $\ell$ are used, the pipelining strategy sketched in Fig. 
\oldref{fig:sketch_amisdcQ_pipelined} is used.  The accuracy and efficiency of the numerical methodology for low-Mach number combustion 
will be assessed with a flame simulation in Section \oldref{section_numerical_examples}.







\section{\label{subsection_convergence_analysis}Convergence of the SDC sweeps}

A key property of all SDC-based methods is their convergence with iteration to the collocation solution used to advance each time step. 
This is to be distinguished from the {\it convergence}\ (in the more traditional sense) of the temporal discretization with decreasing 
time step size.  In this section, we explore the former using a parameterized linear model problem, which allows us to demonstrate 
the behavior of the time stepping scheme in the various limits of competing stiffnesses between the advection, diffusion and reaction 
components of the system.  A simple model problem for this purpose is given by the ODE
\begin{equation}
\left\{
\begin{array}{l l}
  \phi'(t) & = a \phi(t) + d \phi(t) + r \phi(t), \qquad a, d, r \in \mathbb{R} \\[2pt]
  \phi(0)  & = \phi^0,
  \end{array} \right. \label{linear_model_problem}
\end{equation}
where $a$, $d$, and $r$  represent advection, diffusion, and reaction processes, respectively. In order to quantify the convergence of the sweeps, 
we recast the generalized update equations \ref{update_equation_amisdcq_discrete_form_D}-\ref{update_equation_amisdcq_discrete_form_R} into an 
iteration matrix for CISDCQ (and similarly for MISDC and MISDCQ), and we analyze the spectrum of this matrix.

\begin{figure}[ht!]
\centering
\subfigure[]{
\begin{tikzpicture}
\node[anchor=south west,inner sep=0] at (0,0){\includegraphics[scale=0.4]{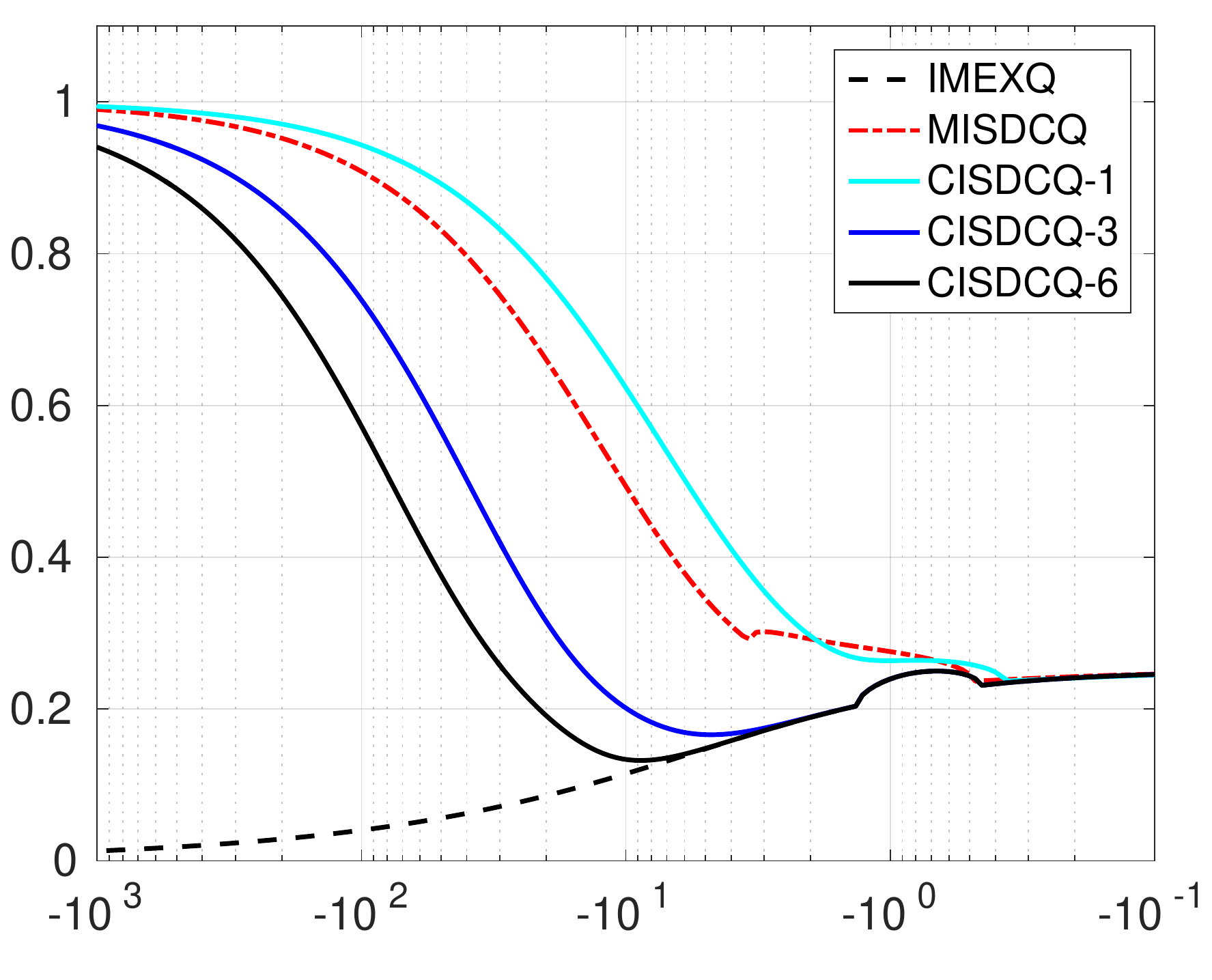}};
\node (ib_1) at (3.7,0.) {\large $r$};
\node[rotate=90] (ib_1) at (-0.3,3.) {\large $\gamma(G)$};
\end{tikzpicture}
\label{fig:convergence_analysis_linear_model_problem_fixed_ratio_d_r_1}
} 
\subfigure[]{
\begin{tikzpicture}
\node[anchor=south west,inner sep=0] at (0,0){\includegraphics[scale=0.4]{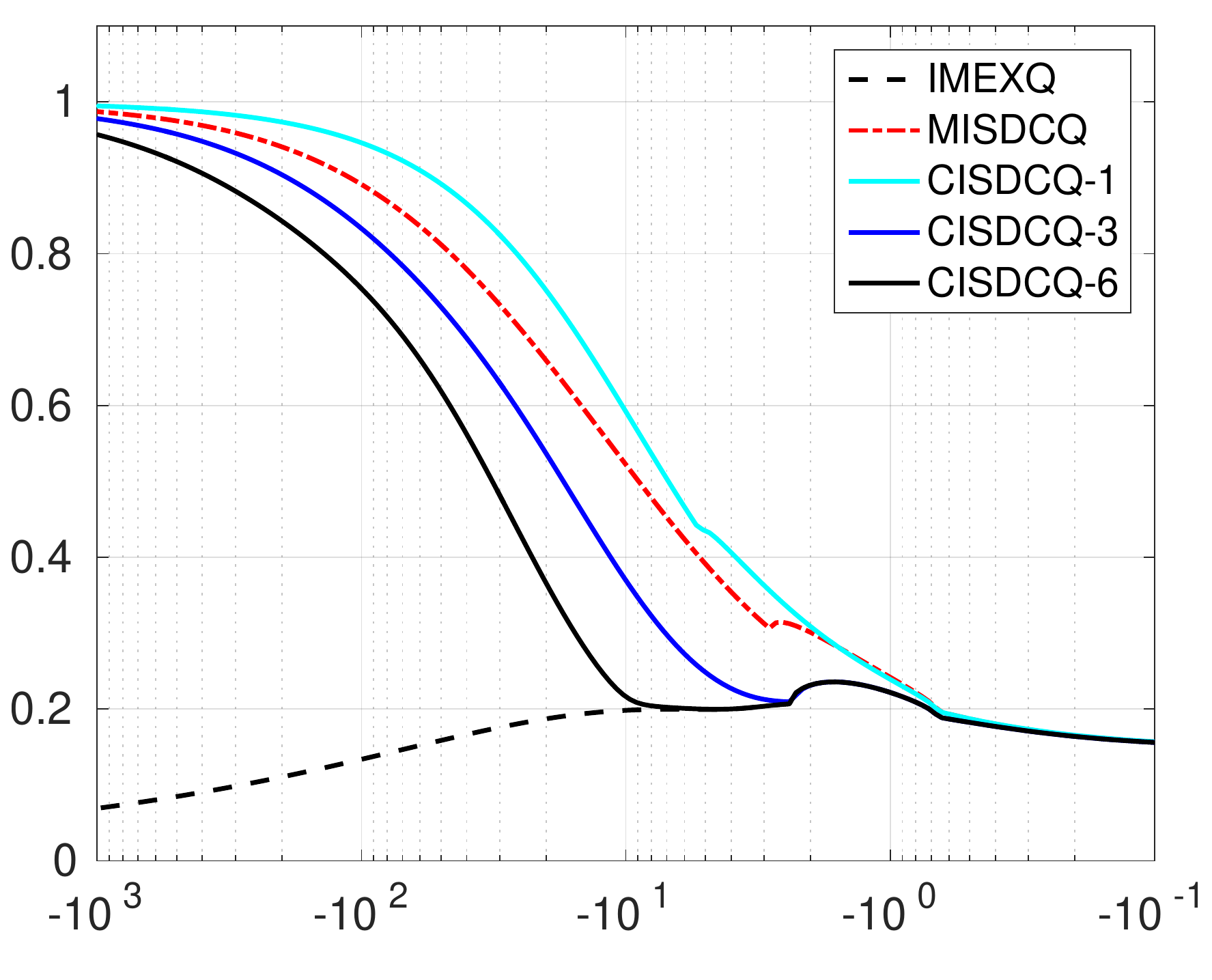}};
\node (ib_1) at (3.7,0.) {\large $r$};
\end{tikzpicture}
\label{fig:convergence_analysis_linear_model_problem_fixed_ratio_d_r_2}
} 
\vspace{-0.4cm}
\caption{\label{fig:convergence_analysis_linear_model_problem_fixed_ratio_d_r}
Spectral radius of  the iteration matrix, $\boldsymbol{G}$, as a function of the reaction coefficient, $r$. The advection and diffusion coefficients are set to $a = 1$
and $d = r/2$, respectively. In \oldref{fig:convergence_analysis_linear_model_problem_fixed_ratio_d_r_1}, three Gauss-Lobatto nodes are used, whereas in 
\oldref{fig:convergence_analysis_linear_model_problem_fixed_ratio_d_r_2}, five Gauss-Lobatto nodes are used. 
}
\end{figure}

In matrix notation, the approximate solutions obtained with CISDCQ at the end of sweep $(k+1)$ are denoted by 
$\boldsymbol{\Phi}^{(k+1)} = [ \phi^{1,(k+1)}, \dots, \phi^{M,(k+1)} ]^T \in \mathbb{R}^M$. This vector is computed iteratively using $\nu$ 
successive nested iterations on $\ell$ applied to the $M$ substeps of the temporal interval. For $m \in \{ 0, \dots, M-1 \}$ and 
$\ell \in \{ 0, \dots, \nu-1 \}$, the intermediate solutions at nested iteration $\ell + 1$, after substep $m+1$, are stored in the vector
\begin{equation}
\boldsymbol{\Phi}^{(k+1, \ell +1)}_{m+1} := [ \phi^{1, (k+1, \ell+1)}, \dots, \phi^{m+1, (k+1, \ell+1)}, 0, \dots, 0 ]^T \in \mathbb{R}^M.
\end{equation}
The vector $\boldsymbol{\Phi}^{(k+1, \ell+1)}_{AD,m+1}$ storing the intermediate solutions after the diffusion update is defined analogously. 
In \oldref{section_derivation_of_the_amisdcq_iteration_matrix_g}, we use the matrix form of the update equations 
\ref{update_equation_amisdcq_discrete_form_D}-\ref{update_equation_amisdcq_discrete_form_R} to show by induction on $m$ and $\ell$ 
that the following relationship holds for all $m \in \{ 0, \dots, M-1 \}$ and $\ell \in \{ 0, \dots, \nu - 1\}$:
\begin{align}
\boldsymbol{\Phi}^{(k+1,\ell+1)}_{m+1} &=  
                                     \boldsymbol{M}^{(\ell+1)}_{1,m+1} \big[ \phi^0 \boldsymbol{1} + s \Delta t \phi^0 \boldsymbol{q} 
                                     + \Delta t \big( s \boldsymbol{\tilde{Q}} - a \boldsymbol{\tilde{Q}}^E - (d + r) \boldsymbol{\tilde{Q}}^I \big) \boldsymbol{\Phi}^{(k)}
                                                                       \big] \nonumber \\
                                   &+ \boldsymbol{M}^{(\ell+1)}_{2,m+1}   \boldsymbol{\Phi}^{(k)},
\label{phi_kp1_as_fct_phi_k_phi_m_phi_l}
\end{align}
where $s = a + d + r$. The matrices $\boldsymbol{M}^{(\ell +1)}_{1,m+1}$, $\boldsymbol{M}^{(\ell +1)}_{2,m+1} \in \mathbb{R}^{M \times M}$ 
depend on the matrices $\boldsymbol{\tilde{Q}}$, $\boldsymbol{\tilde{Q}}^I$, and $\boldsymbol{\tilde{Q}}^E$, on the scalars 
$m$, $\ell$, $\Delta t$, $a$, $d$, and $r$, and can be computed iteratively as explained in the Appendix. The vector $\boldsymbol{q}$ 
is the first column of the matrix $\boldsymbol{Q}$ defined by \ref{weights}. Assuming that the scheme is based on $M+1$ temporal 
nodes and $\nu$ nested iterations on $\ell$, \ref{phi_kp1_as_fct_phi_k_phi_m_phi_l} yields an expression of $\boldsymbol{\Phi}^{(k+1)}$ 
as a function of the previous iterate, $\boldsymbol{\Phi}^{(k)}$, given by
\begin{align}
  \boldsymbol{\Phi}^{(k+1)} = \boldsymbol{\Phi}^{(k+1,\nu)}_{M} &= 
                                     \boldsymbol{M}^{(\nu)}_{1,M} \big[ \phi^0 \boldsymbol{1} + s \Delta t \phi^0 \boldsymbol{q} 
                                     + \Delta t \big( s \boldsymbol{\tilde{Q}} - a \boldsymbol{\tilde{Q}}^E - (d + r) \boldsymbol{\tilde{Q}}^I \big) \boldsymbol{\Phi}^{(k)}
                                                                          \big] \nonumber \\
                                     &+ \boldsymbol{M}^{(\nu)}_{2,M} \boldsymbol{\Phi}^{(k)}.
\label{phi_kp1_as_fct_phi_k}
\end{align}
Taking the difference between two consecutive iterates in \ref{phi_kp1_as_fct_phi_k}, we obtain
\begin{equation}
\boldsymbol{\Phi}^{(k+1)} - \boldsymbol{\Phi}^{(k)} = \boldsymbol{G} ( \boldsymbol{\Phi}^{(k)} - \boldsymbol{\Phi}^{(k-1)} ).
\end{equation}
The iteration matrix, $\boldsymbol{G} \in \mathbb{R}^{M \times M}$, is defined as
\begin{equation}
\boldsymbol{G} := 
 \boldsymbol{M}^{(\nu)}_{1,M} \Delta t \big( s \boldsymbol{\tilde{Q}} - a \boldsymbol{\tilde{Q}}^E - (d + r) \boldsymbol{\tilde{Q}}^I \big) 
                           + \boldsymbol{M}^{(\nu)}_{2,M}.
\end{equation}
We compare the properties of the CISDCQ iteration matrix with those of MISDCQ and IMEXQ. We reiterate that the IMEXQ scheme
is used to illustrate the convergence of the spectral correction process when advection is treated explicitly while 
diffusion and reaction are treated implicitly in a fully coupled fashion. Using these results, we analyze the correction process, which is convergent if 
and only if the spectral radius of the iteration matrix $\boldsymbol{G}$ is strictly smaller than one, that is
\begin{equation}
\gamma( \boldsymbol{G} ) < 1.
\end{equation}
In Fig.~\oldref{fig:convergence_analysis_linear_model_problem_fixed_ratio_d_r}, we first compute the spectral radius of the iteration matrix while 
keeping the ratio $d/r$ fixed. With this approach, we evaluate the asymptotic convergence rate of the schemes when the diffusion and reaction terms 
become stiff simultaneously. We assume a fixed unit time step size $\Delta t = 1$, and we use Gauss-Lobatto nodes. When the problem 
is not stiff ($|r| < 1$), the advection term dominates and all the schemes, including IMEXQ, converge in sweeps to the fixed-point solution at the 
same asymptotic rate. When the problem becomes very stiff ($|r| > 10^3$), the IMEXQ sweeps retain a fast convergence rate and one can show that 
$\lim_{(d+r) \rightarrow - \infty} \gamma(\boldsymbol{G}) = 0$ (see \cite{weiser2015faster}). But, this property does not hold when the diffusion 
and reaction updates are computed in two separate steps. Therefore, the asymptotic convergence rate of the MISDCQ and CISDCQ sweeps
deteriorates significantly for large negative values of $d$ and $r$. Finally, for intermediate values of $r$ ($1 \leq |r| \leq 10^3$), we observe that 
the number of nested iterations on $\ell$ has a strong impact on the asymptotic convergence rate of CISDCQ. If only one nested iteration on $\ell$ 
is used, the CISDCQ-1 sweeps converge slightly slower than with MISDCQ. However, with respectively three and six nested iterations, the CISDCQ-3 
and CISDCQ-6 sweeps converge faster than with MISDCQ. In particular, the CISDCQ-6 sweeps achieve the same convergence rate as IMEXQ for a relatively 
larger fraction of the parameter space. 

\begin{figure}[ht!]
\centering
\subfigure[]{
\begin{tikzpicture}
\node[anchor=south west,inner sep=0] at (0,0){\includegraphics[scale=0.4]{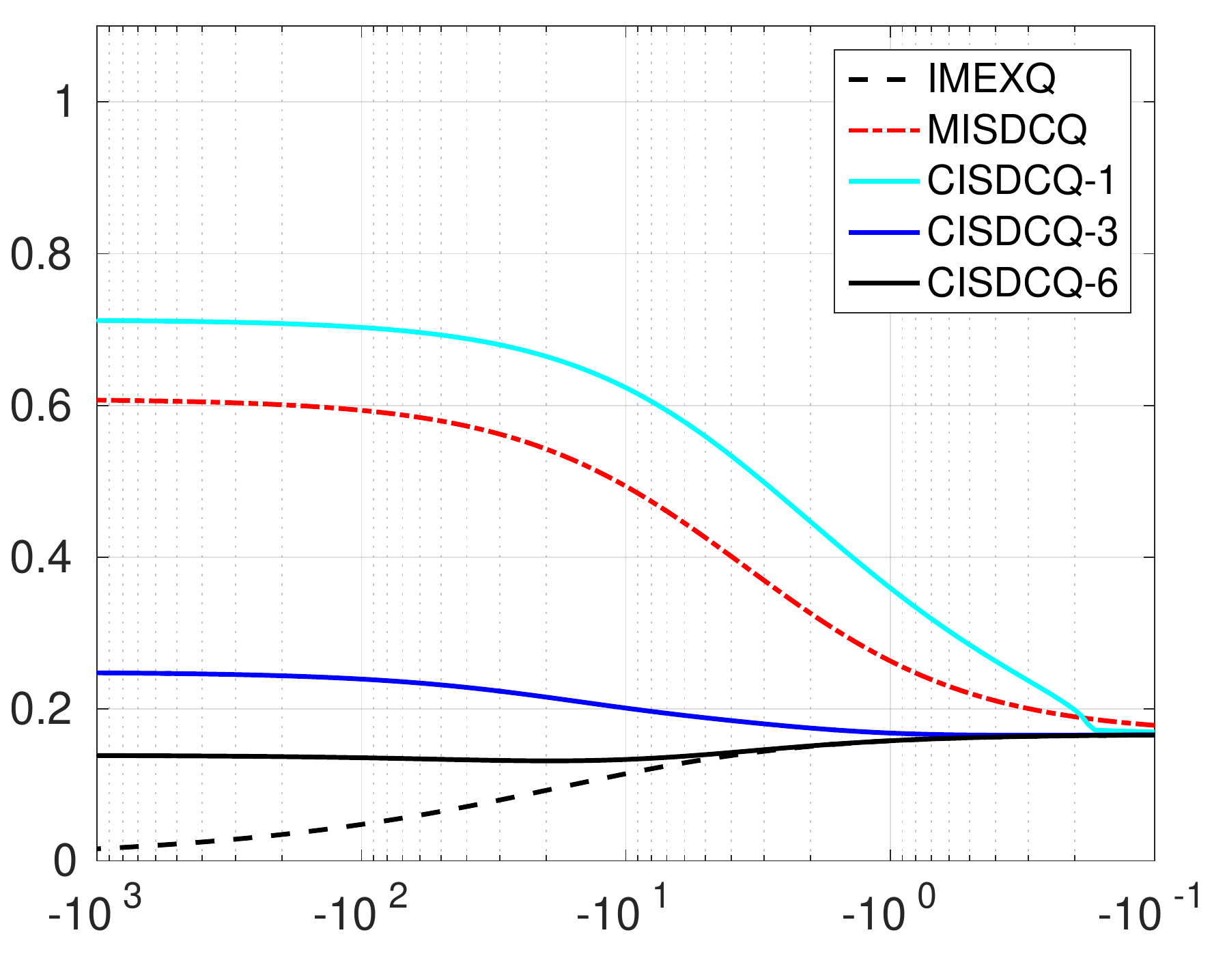}};
\node (ib_1) at (3.7,0.) {\large $r$};
\node[rotate=90] (ib_1) at (-0.3,3.) {\large $\gamma(G)$};
\end{tikzpicture}
\label{fig:convergence_analysis_linear_model_problem_varyin_ratio_d_r_1}
} 
\subfigure[]{
\begin{tikzpicture}
\node[anchor=south west,inner sep=0] at (0,0){\includegraphics[scale=0.4]{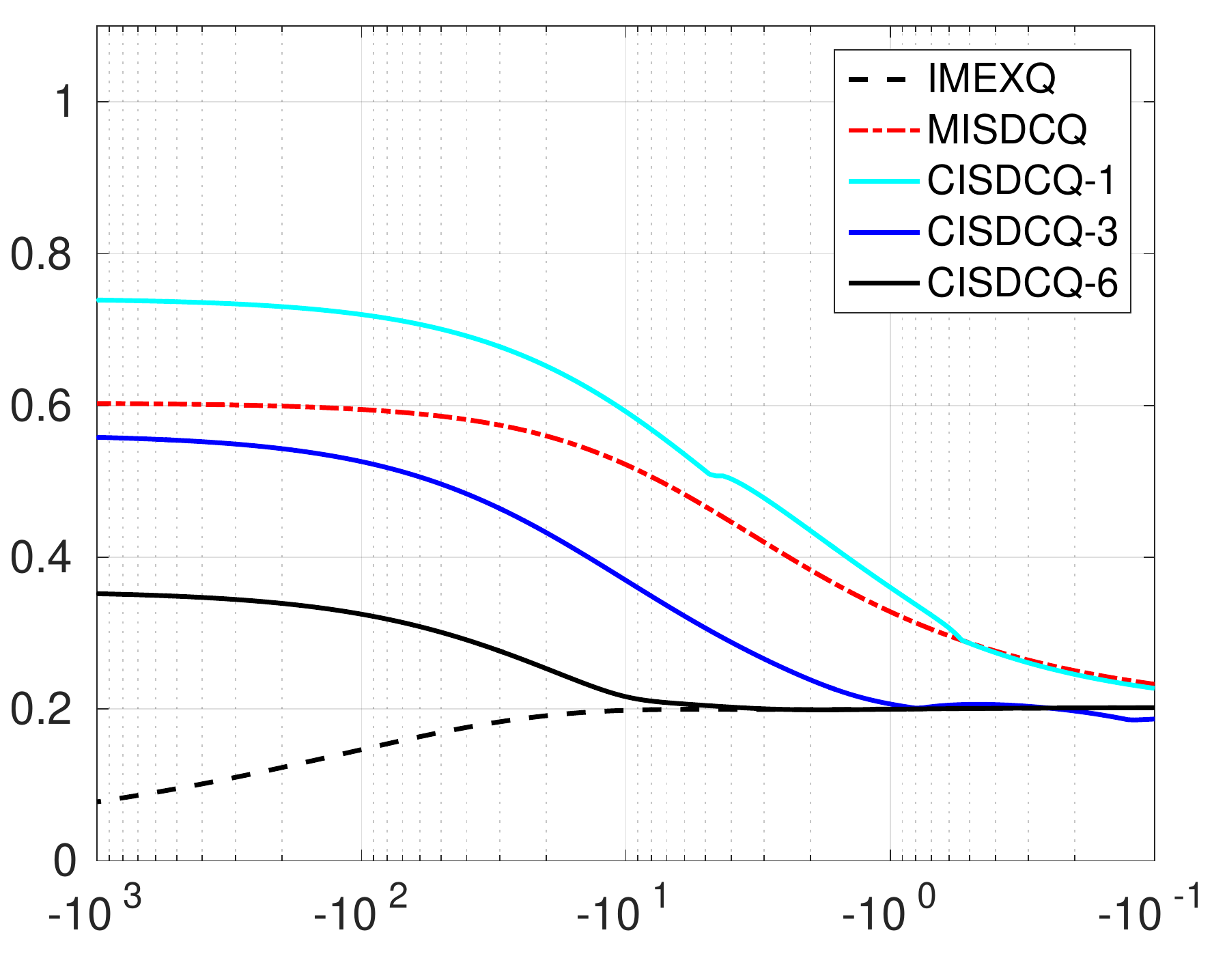}};
\node (ib_1) at (3.7,0.) {\large $r$};
\end{tikzpicture}
\label{fig:convergence_analysis_linear_model_problem_varyin_ratio_d_r_2}
} 

\vspace{-0.4cm}
\caption{\label{fig:convergence_analysis_linear_model_problem_varying_ratio_d_r}
Spectral radius of  the iteration matrix, $\boldsymbol{G}$, as a function of the reaction coefficient, $r$. The advection and diffusion coefficients 
are set to $a = 1$
and $d = -5$, respectively. In \oldref{fig:convergence_analysis_linear_model_problem_varyin_ratio_d_r_1}, three Gauss-Lobatto nodes are used, whereas in 
\oldref{fig:convergence_analysis_linear_model_problem_varyin_ratio_d_r_2}, five Gauss-Lobatto nodes are used. 
}
\end{figure}

In Fig.~\oldref{fig:convergence_analysis_linear_model_problem_varying_ratio_d_r}, we explore another dimension of the parameter space. We compute the 
spectral radius of $\boldsymbol{G}$ while keeping the diffusion coefficient $d$ constant and varying the ratio $d/r$ to evaluate the asymptotic convergence 
rate of the sweeps when the reaction term strongly dominates the diffusion term, or conversely vanishes. We still assume that the time step size is constant, equal to 
one. We observe that with all the schemes, the sweeps converge at the same rate when the reaction term vanishes ($|r| \leq 1$), in which case 
\ref{linear_model_problem} reduces to an advection-diffusion problem. But, if $|r|$ is increased with a constant diffusion strength, increasing the number 
of nested iterations on $\ell$ improves the convergence rate of the sweeps, even when the problem is very stiff ($|r| \geq 10^3$). 
In this configuration, the CISDCQ-1 sweeps still do not converge as fast as the MISDCQ sweeps. But, CISDCQ-3 and CISDCQ-6 converge in sweeps faster 
than the reference scheme.

\section{\label{section_numerical_examples}Numerical examples}
 


For simplicity, we assess the reduction in computational cost obtained with CISDCQ using \ref{simple_parallel_speedup}. This provides a fair comparison
between the time integration methods while using a serial code to simulate parallelism in CISDCQ. The implementation of a truly
parallel version of CISDCQ on a shared-memory platform is a non-trivial task left for future work.

\subsection{\label{subsection_linear_model_problem}Validation of  convergence analysis -- Linear model}

Next, we confirm the analysis above by evolving the linear model problem \ref{linear_model_problem} with the initial condition $\phi^0 = 1$, for several 
representative values of the stiffness parameters, $a$, $d$, and $r$.  For these cases, we discretize the temporal interval $\Delta t = 1$ 
with five Gauss-Lobatto nodes ($M = 4$). We keep the advection coefficient constant equal to $ a=1$ and we vary the diffusion and reaction 
coefficients $d, r \in \mathbb{R}^-$ to assess the stability and convergence rate of the CISDCQ scheme when the problem becomes very stiff. 
We reiterate here that performing the $2 \nu$ pipelined tasks involved in CISDCQ-$\nu$ requires up to $\max(2 \nu, M)$ processors 
working in parallel. Given that $M = 4$ in this section, the CISDCQ-1, CISDCQ-3, and CISDCQ-6 schemes use respectively up to 2, 4, and still 4 processors
to compute the solution.

\begin{figure}[ht!]
\centering
\subfigure[]{
\begin{tikzpicture}
\node[anchor=south west,inner sep=0] at (0,0){\includegraphics[scale=0.26]{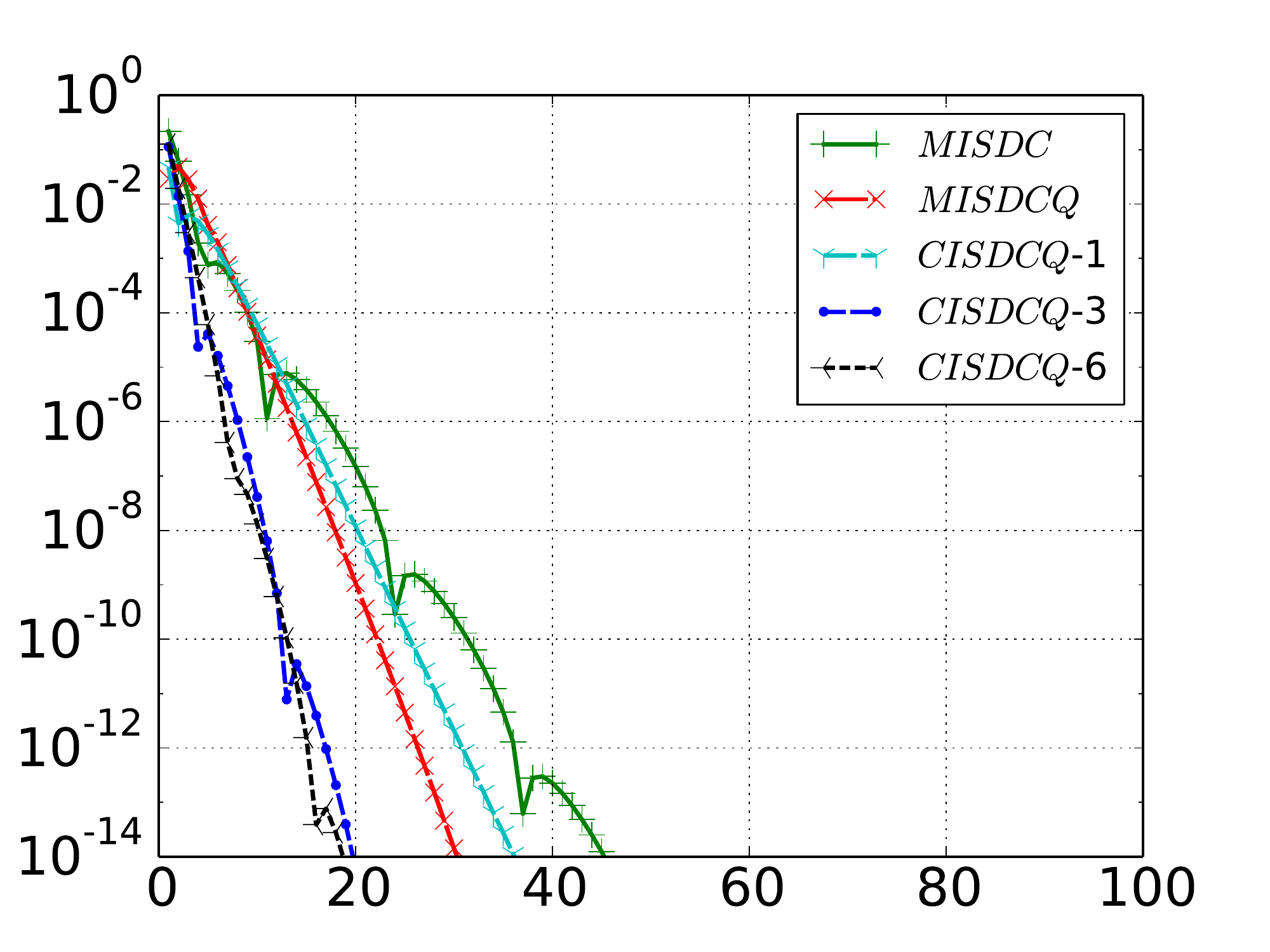}};
\node (ib_1) at (2.7,0.) {\scriptsize Number of sweeps};
\node[rotate=90] (ib_1) at (-0.16,2.) {\scriptsize $| \phi^{M,(k+1)} - \phi^{M,(k)} |$};
\end{tikzpicture}
\label{fig:linear_model_problem_fixed_ratio_d_r_1}
}
\hspace{-0.7cm}
\subfigure[]{
\begin{tikzpicture}
\node[anchor=south west,inner sep=0] at (0,0){\includegraphics[scale=0.26]{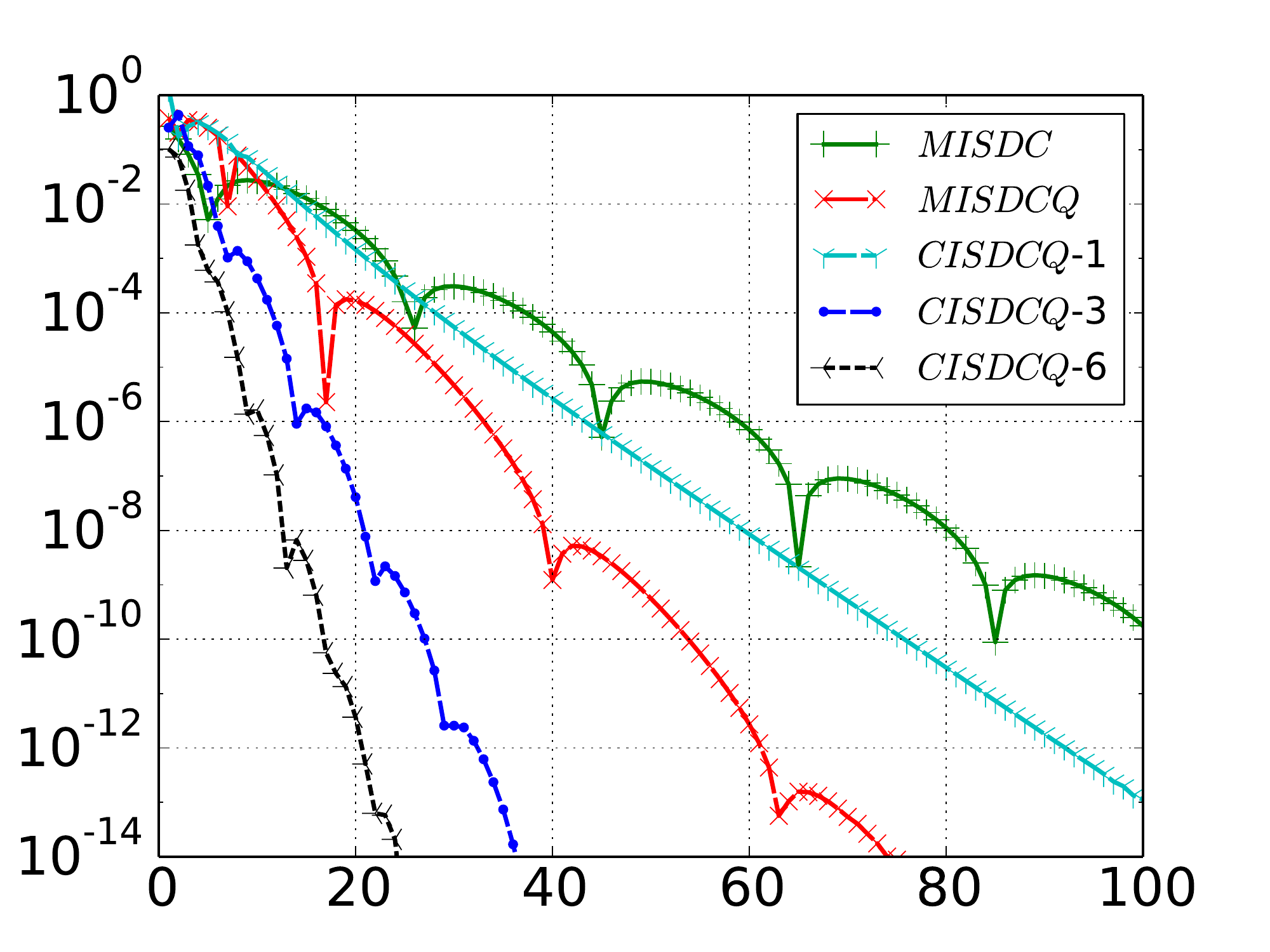}};
\node (ib_1) at (2.7,0.) {\scriptsize Number of sweeps};
\end{tikzpicture}
\label{fig:linear_model_problem_fixed_ratio_d_r_2}
} 
\hspace{-0.7cm}
\subfigure[]{
\begin{tikzpicture}
\node[anchor=south west,inner sep=0] at (0,0){\includegraphics[scale=0.26]{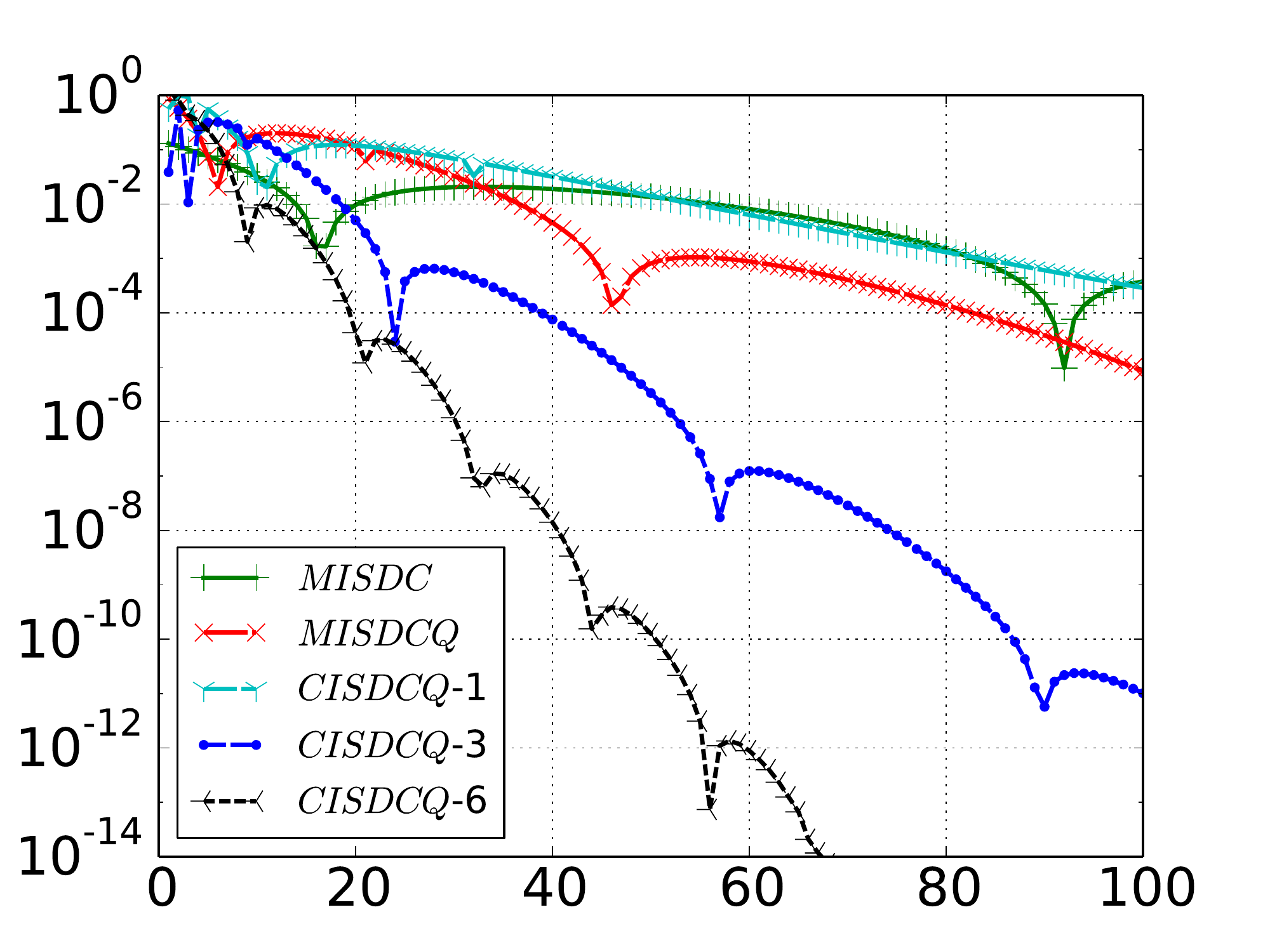}};
\node (ib_1) at (2.7,0.) {\scriptsize Number of sweeps};
\end{tikzpicture}
\label{fig:linear_model_problem_fixed_ratio_d_r_3}
} 

\vspace{-0.4cm}
\caption{\label{fig:linear_model_problem_fixed_ratio_d_r}
Convergence of $| \phi^{M,(k+1)} - \phi^{M,(k)} |$ as a function of the number of sweeps for the linear model problem \ref{linear_model_problem} with $\phi^0 = 1$ and $a = 1$. The diffusion
and reaction coefficients are set to $d = -2$, $r = -4$ in \oldref{fig:linear_model_problem_fixed_ratio_d_r_1}, $d = -10$, $r = -20$ in \oldref{fig:linear_model_problem_fixed_ratio_d_r_2}, and
$d = -50$, $r = -100$ in \oldref{fig:linear_model_problem_fixed_ratio_d_r_3}.
}
\end{figure}

Following the approach of Section \oldref{subsection_convergence_analysis}, we first fix the ratio $r/d = 2$ and increase the absolute value of 
the coefficients $d$ and $r$ at the same rate. Fig.~\oldref{fig:linear_model_problem_fixed_ratio_d_r} shows the convergence rate of the difference 
between two iterates, $| \phi^{M,(k+1)} - \phi^{M,(k)} |$, as a function of the number of sweeps for $r = -4, \, -20, \, -100$. In each case, the 
convergence rate of this quantity is consistent with the analysis of Section \oldref{subsection_convergence_analysis}. Specifically, for the non-stiff
problem of Fig.~\oldref{fig:linear_model_problem_fixed_ratio_d_r_1}, the CISDCQ-1 sweeps converge at a slower rate than those of MISDCQ, whereas the 
CISDCQ-3 and CISDCQ-6 sweeps converge slightly faster than those of MISDC and MISDCQ. As the problem becomes stiffer in 
Figs.~\oldref{fig:linear_model_problem_fixed_ratio_d_r_2} and \oldref{fig:linear_model_problem_fixed_ratio_d_r_3}, CISDCQ-3 and CISDCQ-6 converge in 
sweeps significantly faster than MISDC and MISDCQ. We also note that for the stiff cases, performing six nested iterations on $\ell$ instead of three 
produces a faster convergence of the CISDCQ sweeps. However, increasing the number of nested iterations beyond three does not necessarily improve the 
efficiency of CISDCQ, as shown by the ratio $R_{\nu}$ computed with \ref{simple_parallel_speedup} and the assumption that 
$\Upsilon_{AD} = \Upsilon_{R}$. Considering the number of sweeps necessary to achieve $| \phi^{M,(k+1)} - \phi^{M,(k)} | \leq 10^{-14}$, the values 
of $R_{\nu}$ for each configuration are in Table \oldref{parallel_speedup_linear_model_problem_fixed_ratio_d_r_2}.

\begin{table}[!ht]
\centering
 \scalebox {1}{
         \begin{tabular}{cccc}
           \\ \toprule 
           $r$   & $R_{\nu = 1}$ & $R_{\nu = 3}$ & $R_{\nu = 6}$   \\\toprule
           -4    &    1.4    &    1.5    &   0.9       \\ 
           -20   &    1.1    &    2.6    &   1.6       \\ 
           -100  &    0.9    &    1.8    &   2.0       \\ 
           \bottomrule 
         \end{tabular}}
\caption{\label{parallel_speedup_linear_model_problem_fixed_ratio_d_r_2}
Ratio of the computational cost of MISDCQ over that of CISDCQ-$\nu$ for the linear model problem \ref{linear_model_problem} with $a = 1$ and a fixed ratio $r/d = 2$.
}
\end{table}




\begin{figure}[ht!]
\centering
\subfigure[]{
\begin{tikzpicture}
\node[anchor=south west,inner sep=0] at (0,0){\includegraphics[scale=0.26]{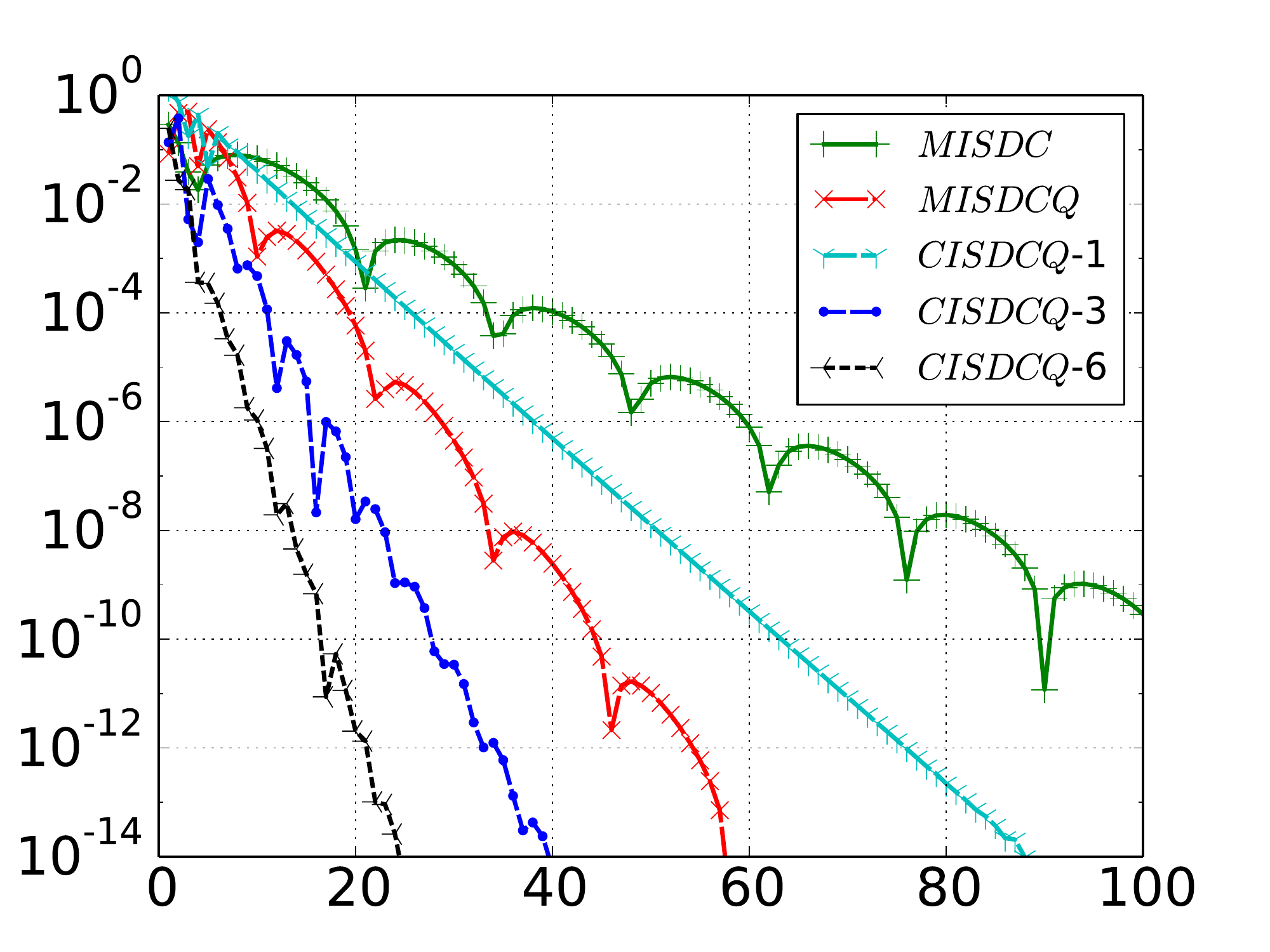}};
\node (ib_1) at (2.7,0.) {\scriptsize Number of sweeps};
 \node[rotate=90] (ib_1) at (-0.16,2.) {\scriptsize $| \phi^{M,(k+1)} - \phi^{M,(k)} |$};
\end{tikzpicture}
\label{fig:linear_model_problem_varyin_ratio_d_r_1}
}
\hspace{-0.7cm}
\subfigure[]{
\begin{tikzpicture}
\node[anchor=south west,inner sep=0] at (0,0){\includegraphics[scale=0.26]{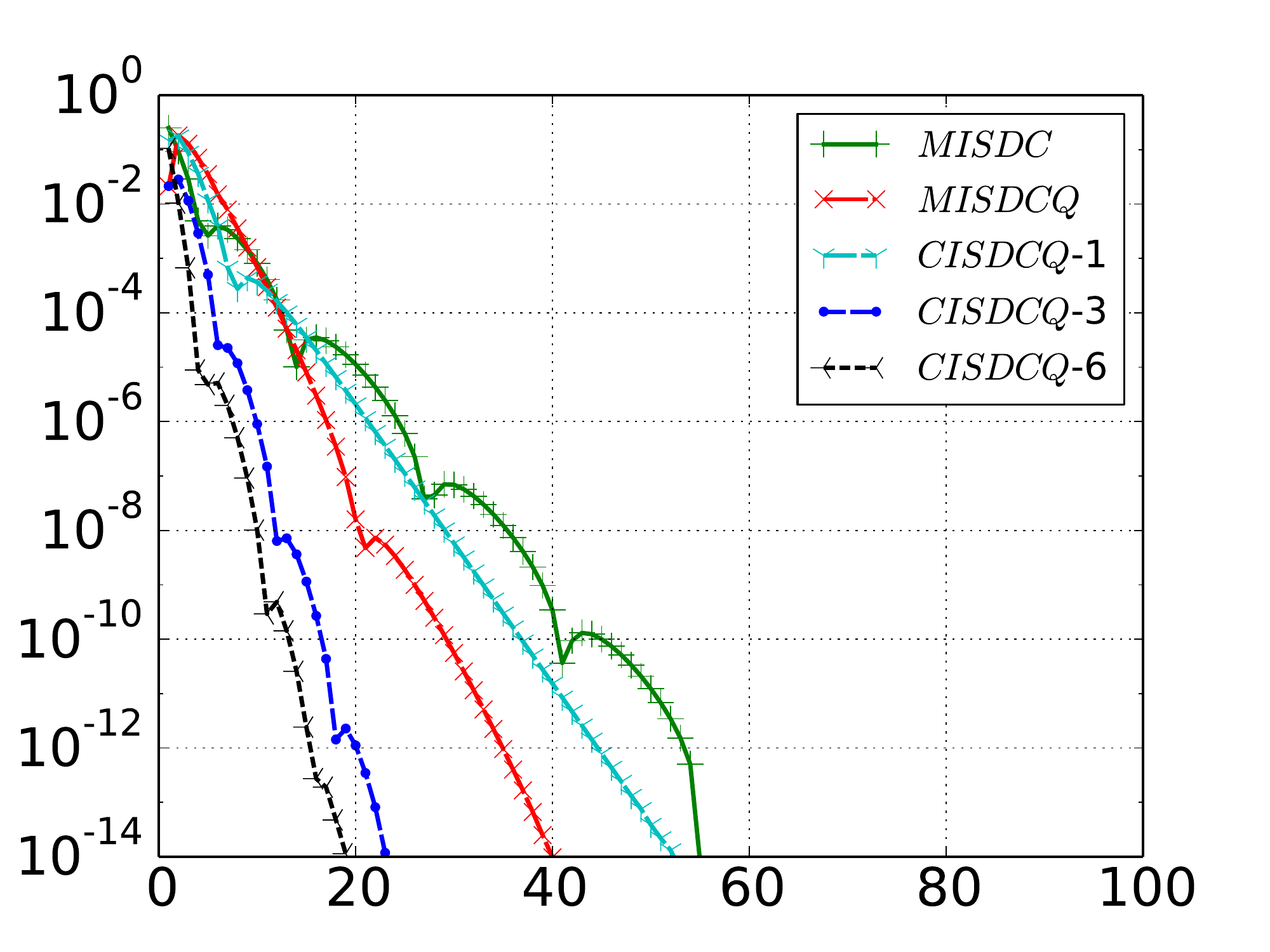}};
\node (ib_1) at (2.7,0.) {\scriptsize Number of sweeps};
\end{tikzpicture}
\label{fig:linear_model_problem_varyin_ratio_d_r_2}
} 
\hspace{-0.7cm}
\subfigure[]{
\begin{tikzpicture}
\node[anchor=south west,inner sep=0] at (0,0){\includegraphics[scale=0.26]{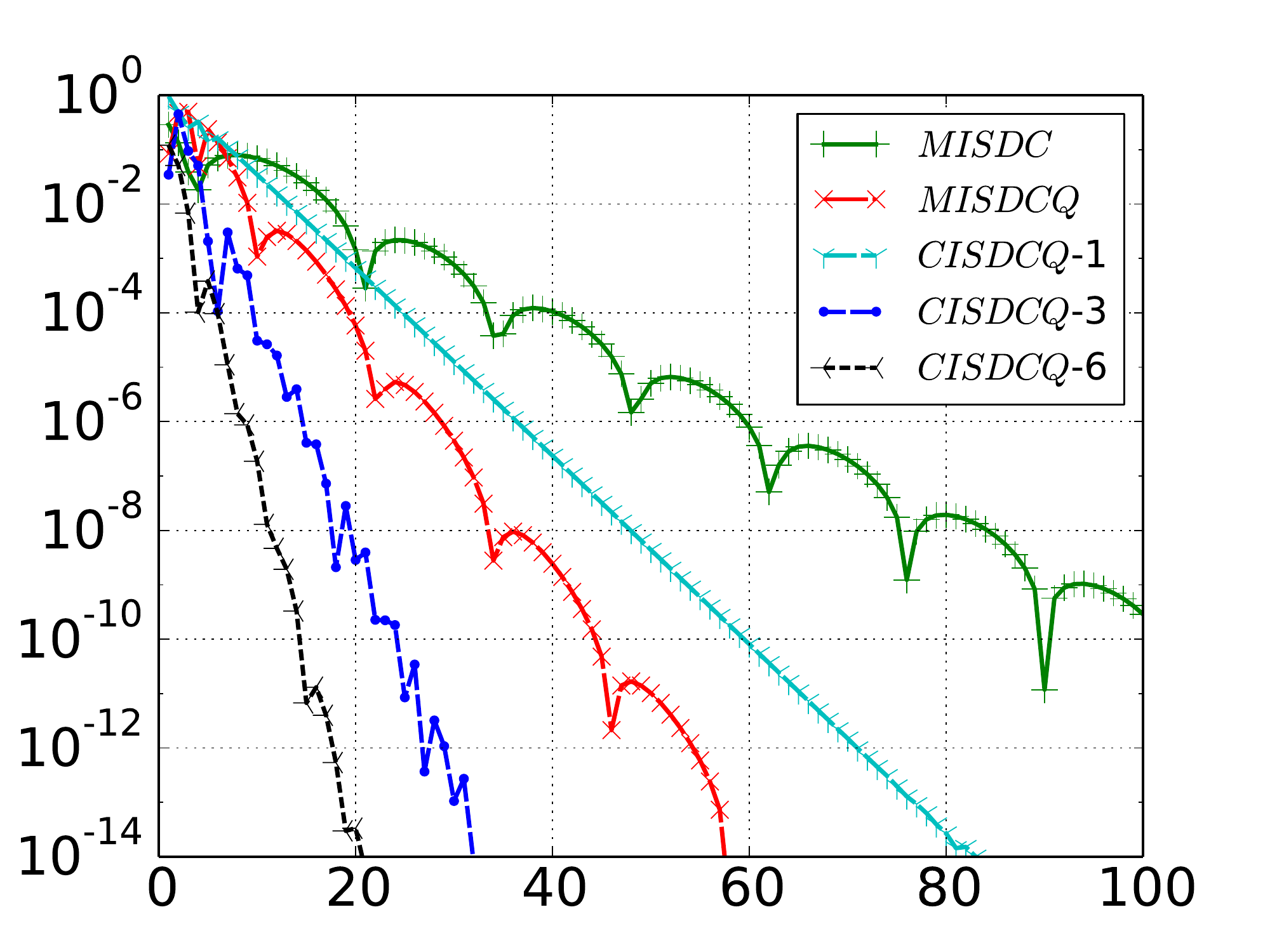}};
\node (ib_1) at (2.7,0.) {\scriptsize Number of sweeps};
\end{tikzpicture}
\label{fig:linear_model_problem_varyin_ratio_d_r_3}
} 

\vspace{-0.4cm}
\caption{\label{fig:linear_model_problem_varying_ratio_d_r}
Convergence of $| \phi^{M,(k+1)} - \phi^{M,(k)} |$ as a function of the number of sweeps for the linear model problem \ref{linear_model_problem} with 
$\phi^0 = 1$ and $a = 1$. The diffusion and reaction coefficients are set to $d = -100$, $r = -5$ in \oldref{fig:linear_model_problem_varyin_ratio_d_r_1}, 
$d = -5$, $r = -5$ in \oldref{fig:linear_model_problem_varyin_ratio_d_r_2}, and $d = -5$, $r = -100$ in \oldref{fig:linear_model_problem_varyin_ratio_d_r_3}.
}
\end{figure}

Next, we study the computational cost of CISDCQ for multiple ratios $r/d$ to make sure that the scheme remains stable when the 
diffusion term dominates the reaction term, and vice-versa. The convergence of $| \phi^{M,(k+1)} - \phi^{M,(k)} |$ as a function of the number of sweeps 
for these cases is in Fig.~\oldref{fig:linear_model_problem_varying_ratio_d_r}. The CISDCQ-3 and CISDCQ-6 sweeps converge faster than those of MISDCQ 
when the absolute magnitude of the diffusion term is 20 times larger than that of the reaction term ($d = -100, \, r = -5$) in Fig. 
\oldref{fig:linear_model_problem_varyin_ratio_d_r_1}, when the two terms have the same magnitude ($d = r = -5$) in Fig. 
\oldref{fig:linear_model_problem_varyin_ratio_d_r_2}, and when the absolute magnitude of the diffusion term is 20 times smaller than that of the reaction 
term ($d = -5, r = -100$) in Fig.~\oldref{fig:linear_model_problem_varyin_ratio_d_r_3}. We observe again that increasing the number of nested iterations 
on $\ell$ beyond three accelerates the convergence of the sweeps but can degrade the reduction in computational cost obtained with CISDCQ. 
The ratio $R_{\nu}$ for a tolerance of $10^{-14}$ is in Table \oldref{parallel_speedup_linear_model_problem_varyin_ratio_d_r}.
This highlights the need for an adaptive mechanism to select the optimal number of nested iterations on $\ell$ at each sweep based 
on the stiffness of the problem, the number of SDC nodes, and the number of available processors. Such a strategy would in particular stop the nested iterations 
before the CISDCQ performance starts degrading. This will be explored in future work.




\begin{table}[!ht]
\centering
 \scalebox {1}{
         \begin{tabular}{ccccc}
           \\ \toprule 
           $d$  & $r$      & $R_{\nu = 1}$ & $R_{\nu = 3}$ & $R_{\nu = 6}$   \\\toprule
           -100 & -5       &    1.0    &    1.2    &   1.2       \\ 
           -5   & -5       &    2.1    &    1.5    &   1.1       \\ 
           -5   & -100     &    1.1    &    1.6    &   1.4       \\ 
           \bottomrule 
         \end{tabular}}
\caption{\label{parallel_speedup_linear_model_problem_varyin_ratio_d_r}
Ratio of the computational cost of MISDCQ over that of CISDCQ-$\nu$ for the linear model problem \ref{linear_model_problem} with $a = 1$ and 
multiple values of the ratio $r/d$.
}
\end{table}

\subsection{\label{subsubsection_nonlinear_test_pde}Nonlinear test PDE}
It is instructive to verify that the results of the previous section can predict reasonably well the performance
of CISDCQ on a more complex nonlinear coupled problem.  The one-dimensional advection-diffusion-reaction model,
\begin{equation}
\left\{
\begin{array}{lll}
\phi_t(x,t) & = & a \phi_x + d \phi_{xx} + r \phi (\phi - 1) ( \phi - 1/2 ) \qquad \text{for} \, (x,t) \in [0,20] \times [0,T], \\
\phi(0,t)   & = & 1, \\
\phi(20,t)  & = & 0, \\
\phi(x,0)   & = & \phi^{0}(x),
\end{array}
\right. \label{nonlinear_test_pde}
\end{equation}
captures many of the generic features and coupling of the low-Mach combustion system, yet avoids complexities associated with parameterized
equations of state, transport properties, and reaction physics.  A suitable initial condition is given by
\begin{equation}
\phi^0(x) = \frac{1}{2} \big( 1 + \tanh(20 - 2x) \big).
\end{equation}
This PDE is solved with the method of lines based on a finite-volume approach. The advection term is discretized using a fourth-order 
operator, $A(\phi)$, and the diffusion term is approximated with a fourth-order Laplacian operator, $D(\phi)$. This discretization scheme 
leads to the nonlinear system of ODEs
\begin{equation}
\phi_t = A(\phi) + D(\phi) + R(\phi), \end{equation}
where $R(\phi) = r \phi ( \phi - 1) ( \phi - 1/2 )$ denotes the reaction term. This system is solved with the MISDC, MISDCQ, and CISDCQ schemes based on 
five Gauss-Lobatto nodes to compare their accuracy and number of sweeps needed to achieve convergence to the fixed-point solution. As 
in the previous section, CISDCQ-1, CISDCQ-3, and CISDCQ-6 allow to use respectively up to 2, 4, and 4 parallel processors.

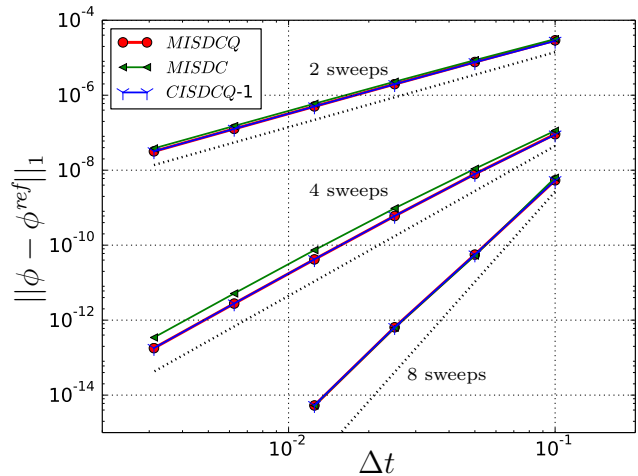
\begin{figure}[ht!]
\centering
\begin{tikzpicture}
\node[anchor=south west,inner sep=0] at (0,0){\includegraphics[scale=0.45]{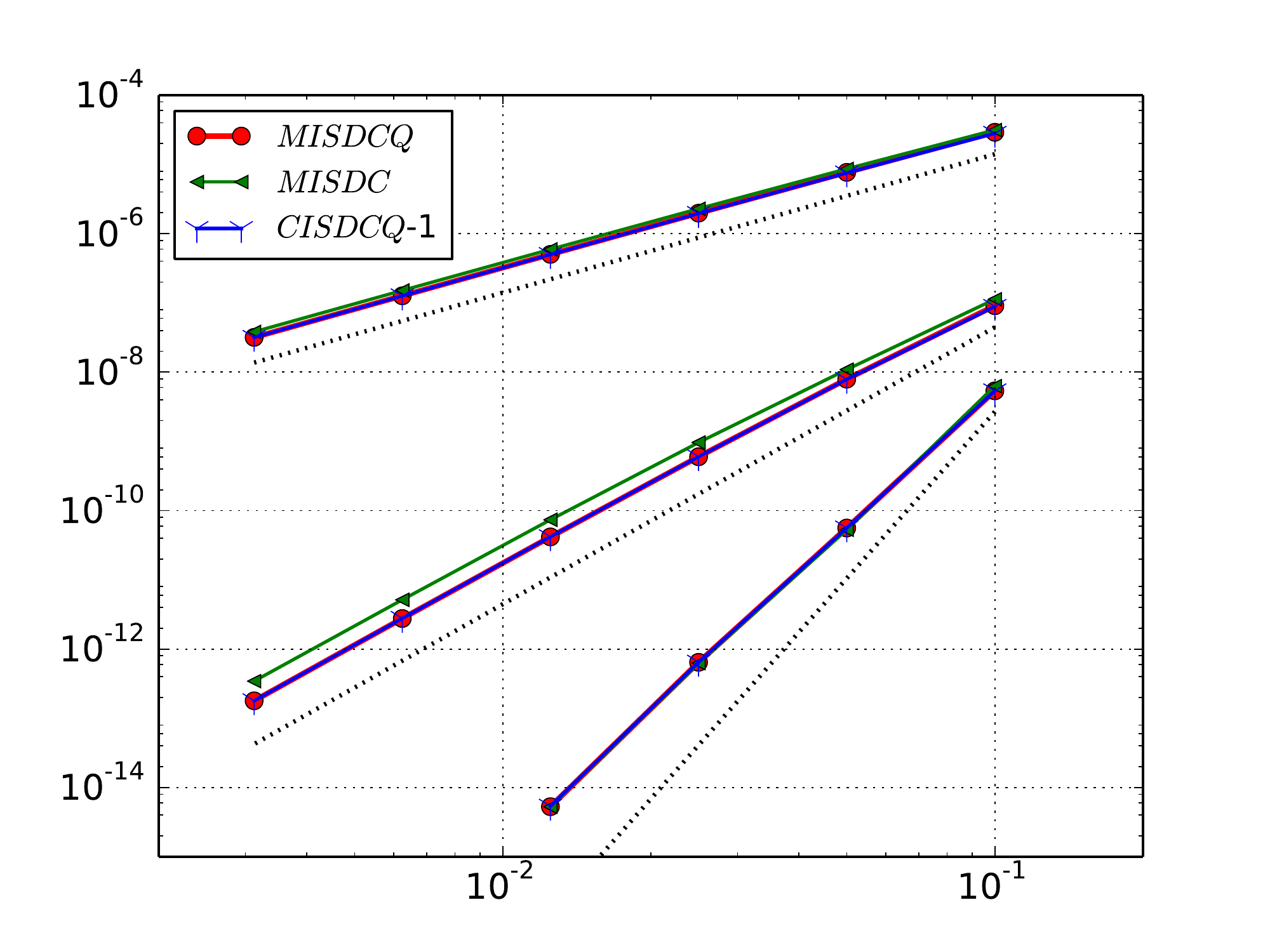}};
\node (ib_1) at (4.8,0.3) {\large $\Delta t$};
\node[rotate=90] (ib_1) at (0.175,3.4) {\large $|| \phi - \phi^{\textit{ref}} ||_1$};
\node (ib_1) at (4.43,5.5) {\scriptsize 2 sweeps};
\node (ib_1) at (4.43,3.9) {\scriptsize 4 sweeps};
\node (ib_1) at (5.73,1.45) {\scriptsize 8 sweeps};
\end{tikzpicture}
\vspace{-0.1cm}
\caption{\label{fig:nonlinear_problem_batstoi_convergence_refinement_in_time}
$L_1$-error of the MISDC, MISDCQ, and CISDCQ-1 solutions with respect to the reference solution as a function of time for the nonlinear test PDE 
\ref{nonlinear_test_pde} with $a = 1$, $d = 2$, $r = 4$. The dashed lines indicate the theoretical order of accuracy of the MISDC-type schemes. The solid 
lines corresponding to the MISDCQ and CISDCQ-1 schemes are on top of each other.
}
\end{figure}

We first perform a refinement study in time to assess the accuracy of CISDCQ. We verify that each CISDCQ sweep increases the temporal 
order of accuracy by one for mildly stiff problems. To this end, we refine the time step size while keeping the grid spacing fixed to 
$\Delta x = 0.1$ ($n_x = 200$). We set the coefficients in \ref{nonlinear_test_pde} to $a = 1$, $d = 2$, and $r = 4$.
Here, the error is quantified using the $L_1$-norm
\begin{equation}
|| \phi - \phi^{\textit{ref}} ||_1 = \frac{1}{n_{x}} \sum_{i = 1}^{n_{x}}  | \phi(x_i,t_{\textit{final}}) - \phi^{\textit{ref}}(x_i,t_{\textit{final}}) |,
\end{equation}
where $n_{x}$ is the number of cells and $\phi$ (respectively, $\phi^{\textit{ref}}$) denotes the approximate solution (respectively, the averaged reference solution). 
The results are compared to a reference solution generated with MISDCQ for $\Delta t_{\textit{ref}} = \Delta x / 32 = 0.003125$, which corresponds to an 
advective CFL number of $\sigma_{\textit{ref}} = 0.03125$. Fig.~\oldref{fig:nonlinear_problem_batstoi_convergence_refinement_in_time} shows 
the $L_1$-norm of the error with respect to the reference solution as a function of the time step size. The error obtained with 
CISDCQ-1 -- i.e., the CISDCQ scheme based on only one nested iteration on $\ell$ -- decreases at the same rate as with MISDC and MISDCQ 
when the time step size is reduced.  Specifically, CISDCQ-1 achieves second-order accuracy with two sweeps, and fourth-order 
accuracy with four sweeps. With eight sweeps, the asymptotic range is not reached and the three schemes only exhibit a seventh-order accuracy. 
The results obtained with CISDCQ-3 and CISDCQ-6, based on three and six nested iterations, respectively, are similar to those 
of CISDCQ-1 and are therefore not presented.

\begin{figure}[ht!]
\centering
\subfigure[]{
\begin{tikzpicture}
\node[anchor=south west,inner sep=0] at (0,0){\includegraphics[scale=0.26]{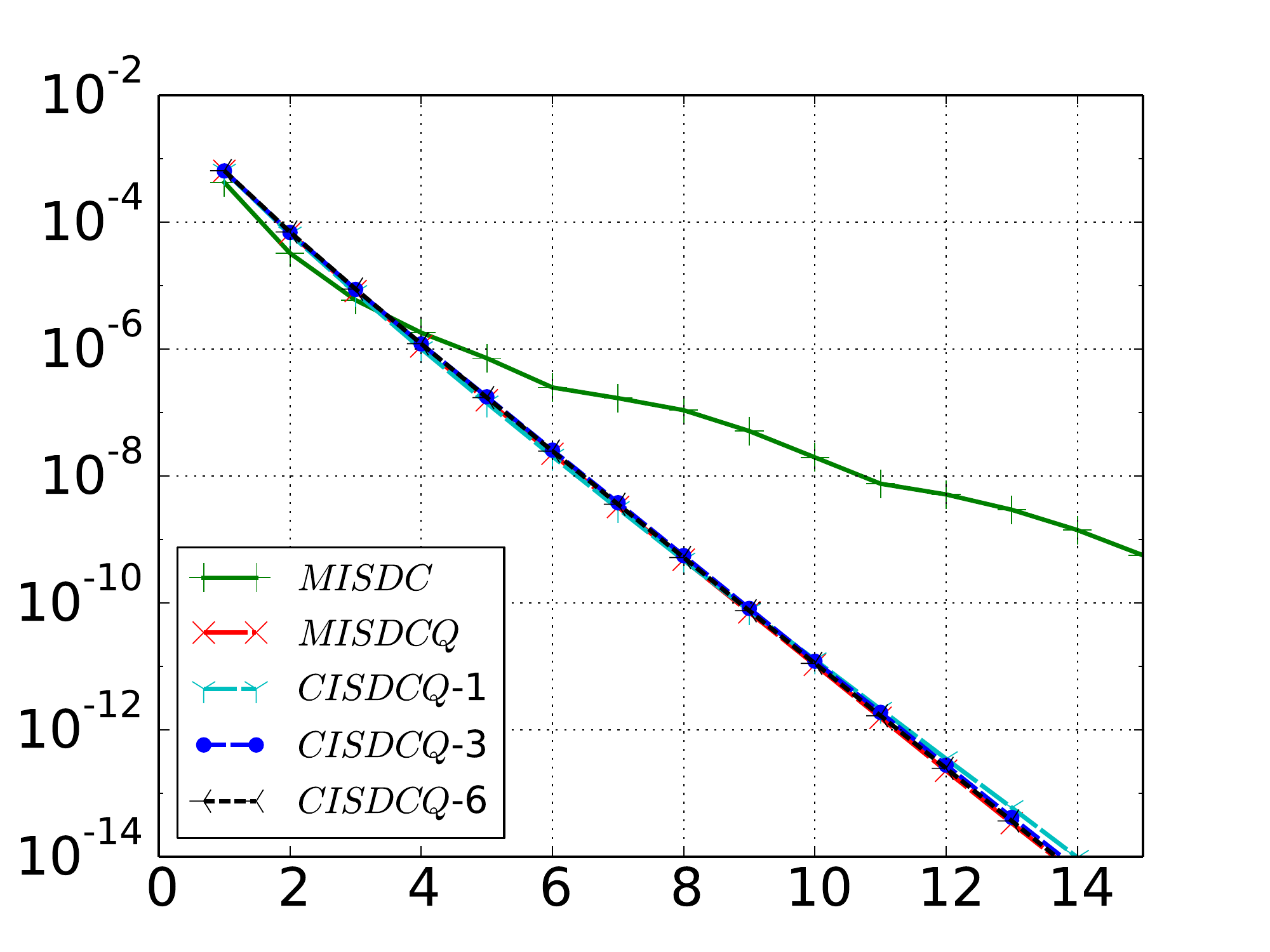}};
\node (ib_1) at (2.725,0.) {\scriptsize Number of sweeps};
 \node[rotate=90] (ib_1) at (-0.16,2.) {\scriptsize $|| \phi^{M,(k+1)} - \phi^{M,(k)} ||_1$};
\end{tikzpicture}
\label{fig:nonlinear_model_problem_fixed_ratio_d_r_1}
} 
\hspace{-0.7cm}
\subfigure[]{
\begin{tikzpicture}
\node[anchor=south west,inner sep=0] at (0,0){\includegraphics[scale=0.26]{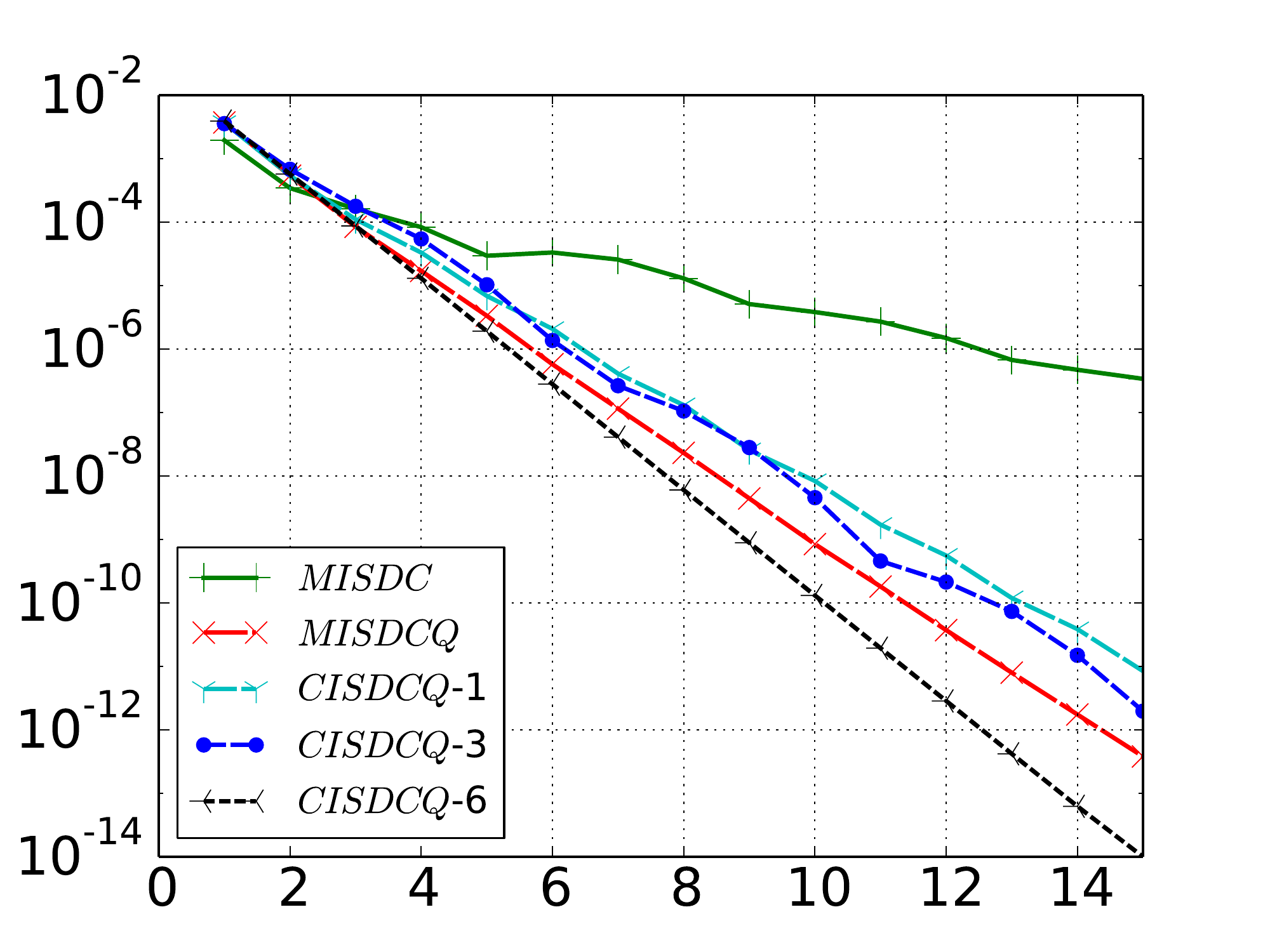}};
\node (ib_1) at (2.725,0.) {\scriptsize Number of sweeps};
\end{tikzpicture}
\label{fig:nonlinear_model_problem_fixed_ratio_d_r_2}
} 
\hspace{-0.7cm}
\subfigure[]{
\begin{tikzpicture}
\node[anchor=south west,inner sep=0] at (0,0){\includegraphics[scale=0.26]{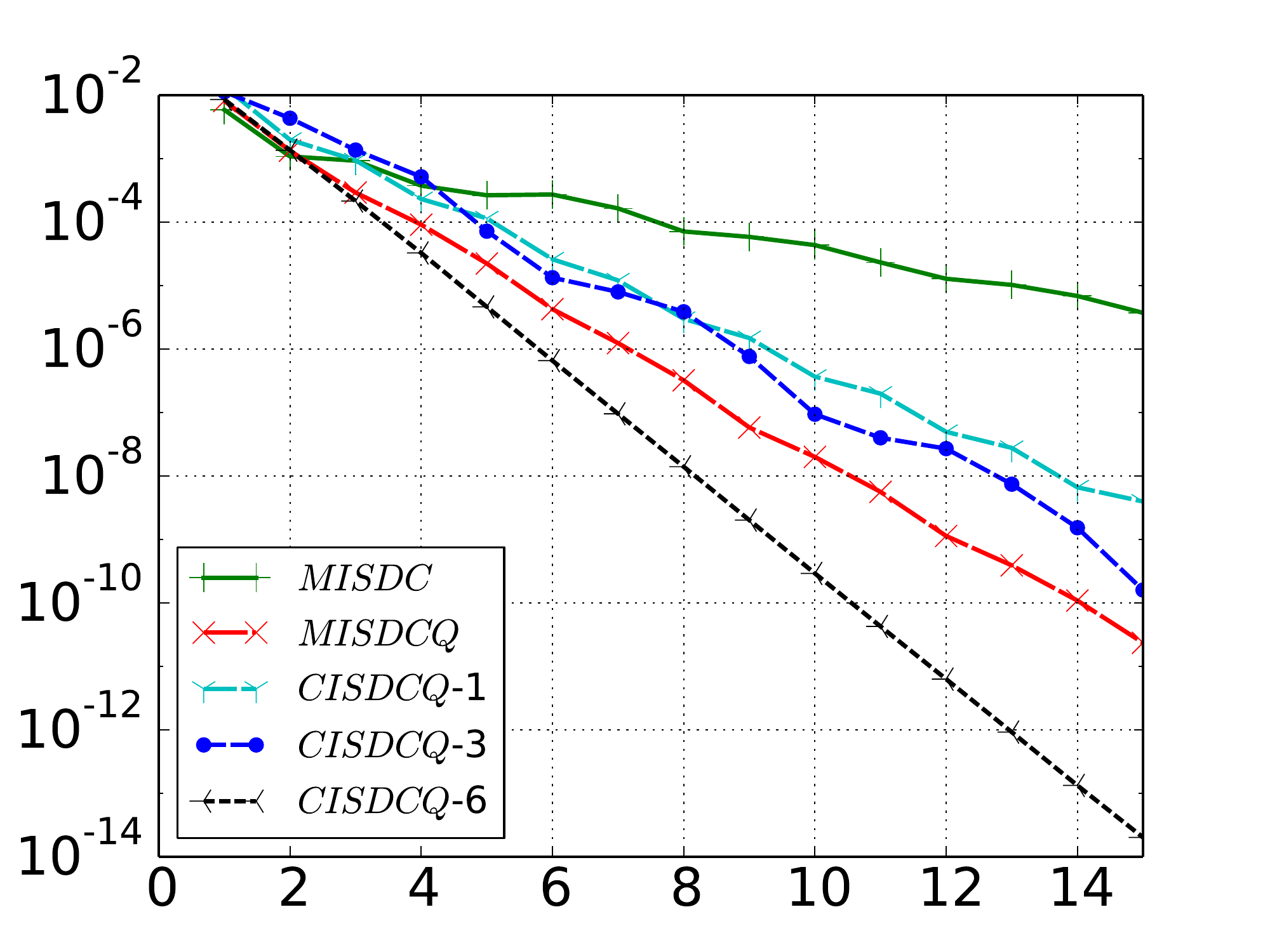}};
\node (ib_1) at (2.725,0.) {\scriptsize Number of sweeps};
\end{tikzpicture}
\label{fig:nonlinear_model_problem_fixed_ratio_d_r_3}
} 

\vspace{-0.4cm}
\caption{\label{fig:nonlinear_problem_batstoi_convergence_as_fct_number_of_sweeps}
Convergence of $|| \phi^{M,(k+1)} - \phi^{M,(k)} ||_1$ as a function of the number of sweeps for the nonlinear test PDE \ref{nonlinear_test_pde} 
with $a = 1$. 
The diffusion and reaction coefficients are set to $d = 2$, $r = 4$ in \oldref{fig:nonlinear_model_problem_fixed_ratio_d_r_1}, $d = 8$, $r = 16$ in 
\oldref{fig:nonlinear_model_problem_fixed_ratio_d_r_2}, and $d = 16$, $r = 32$ in \oldref{fig:nonlinear_model_problem_fixed_ratio_d_r_3}.
}
\end{figure}

Next, we consider a single time step and we investigate the convergence of the iterative spectral deferred correction process for different levels 
of stiffness as in Section \oldref{subsection_linear_model_problem}. The spatial discretization is still based on $n_x = 200$ cells and the time step 
is fixed to $\Delta t = \Delta x / 2 = 0.05$. The advection coefficient is still set to $a = 1$, which corresponds to an advective CFL number of 
$\sigma = 1/2$. We consider three cases with increasingly large  value of the diffusion and reaction coefficients $d$ and $r$. The convergence 
of the $L_1$-norm of the difference between two iterates, $|| \phi^{M,(k+1)} - \phi^{M,(k)} ||_1$, as a function of the number of sweeps for the first 
time step of the simulation is in 
Fig.~\oldref{fig:nonlinear_problem_batstoi_convergence_as_fct_number_of_sweeps}.

To evaluate the computational cost of CISDCQ-$\nu$ after 15 sweeps, we first compute the value $\epsilon \in \mathbb{R}^+$ such that 
$|| \phi^{M,(15)} - \phi^{M,(14)} ||_1 = \epsilon$, where $\phi^{M,(14)}$ and $\phi^{M,(15)}$ are the approximate solutions generated 
by MISDCQ after 14 and 15 sweeps, respectively. Then, we compute the number of sweeps $N_C$ necessary to achieve 
$|| \phi^{M,(N_C)} - \phi^{M,(N_C-1)} ||_1 \leq \epsilon$ with CISDCQ-$\nu$. Finally, using \ref{simple_parallel_speedup} 
and the assumption that $\Upsilon_{AD} = \Upsilon_R$, we compute the ratio $R_{\nu}$ for the CISDCQ-$\nu$ integration scheme 
as $R_{\nu} = 15 / N_C \times 2 M / ( 2 \nu  + M - 1 )$. The results can be found in Table \oldref{parallel_speedup_nonlinear_problem}.

In Fig.~\oldref{fig:nonlinear_model_problem_fixed_ratio_d_r_1}, the sweeps of the integration methods based on the LU decomposition of the integration
matrix, namely, MISDCQ and CISDCQ, converge at the same rate. The MISDC sweeps converge significantly slower and this scheme requires almost twice 
as many sweeps to achieve $|| \phi^{M,(k+1)} - \phi^{M,(k)} ||_1 \leq 10^{-11}$. For this mildly stiff problem, increasing the number of iterations 
on $\ell$ in the nested loop 
does not accelerate the convergence of the CISDCQ sweeps. Therefore, it is more advantageous to use CISDCQ-1 -- i.e., only one iteration in the nested 
loop --, which has a computational cost ratio of $R_{\nu = 1} = 1.6$ when $d = 2$ and $r = 4$. In CISDCQ-3 and CISDCQ-6, the overhead 
caused by the nested loop is not compensated for by an enhanced convergence of the sweeps, which is why we only obtain $R_{\nu = 3} \approx 0.9$ and 
$R_{\nu = 6} \approx 0.5$ for $d = 2$ and $r = 4$.

But, in Figs.~\oldref{fig:nonlinear_model_problem_fixed_ratio_d_r_2} and \oldref{fig:nonlinear_model_problem_fixed_ratio_d_r_3}, we clearly see that 
the nested iterations on $\ell$ have a positive impact on the convergence of the CISDCQ-3 and CISDCQ-6 sweeps when the problem becomes stiffer.  
For instance, the convergence of the CISDCQ-1 sweeps deteriorates significantly when $d$ and $r$ increase, whereas that of CISDCQ-6 is not affected 
by the stronger stiffness of the problem. However, this enhanced convergence of the CISDCQ-3 and CISDCQ-6 sweeps is not sufficient to yield a significant 
reduction in computational cost. In fact, in these two configurations, it is still more efficient to use CISDCQ-1 to reduce the 
cost, as shown by the values of $R_{\nu}$. We obtain $R_{\nu = 1} \approx 1.4$, $R_{\nu = 3} \approx 0.8$, and $R_{\nu = 6} \approx 0.6$ for $d = 8$ 
and $r = 16$. When $d = 16$ and $r = 32$, we achieve $R_{\nu = 1} \approx 1.2$, $R_{\nu = 3} \approx 0.7$, and $R_{\nu = 6} \approx 0.7$.
As in the previous example, an adaptive strategy to stop the nested iterations before $R_{\nu}$ starts decreasing would considerably benefit 
the CISDCQ algorithm.

\begin{table}[!ht]
\centering
 \scalebox {1}{
         \begin{tabular}{ccccc}
           \\ \toprule 
           $d$  & $r$   & $R_{\nu = 1}$ & $R_{\nu = 3}$ & $R_{\nu = 6}$   \\\toprule
           2    & 4     &    1.6      &    0.9      &   0.5       \\ 
           8    & 16    &    1.4      &    0.8      &   0.6       \\ 
           16   & 32    &    1.2      &    0.7      &   0.7       \\ 
           \bottomrule 
         \end{tabular}}
\caption{\label{parallel_speedup_nonlinear_problem}
Ratio of the computational cost of MISDCQ over that of CISDCQ-$\nu$ for the nonlinear test PDE \ref{nonlinear_test_pde} with $a = 1$ and 
multiple values of  $d$ and $r$.
}
\end{table}

\subsection{Dimethyl Ether Flame}
 
In this section, we assess the performance of CISDCQ using a one-dimensional unsteady simulation of a premixed flame based 
on a 39-species, 175-reaction dimethyl ether (DME) chemistry mechanism \citep{bansal2015direct}. This example is challenging for the SDC schemes discussed in this work because 
the DME chemistry mechanism is very stiff. The system is evolved on the relatively slow advection time scale to evaluate the 
ability of the schemes to capture the nonlinear coupling between the faster diffusion and reaction processes. The domain length is $d = 0.6 \, \text{cm}$. 
The inlet stream at $T = 298 \, \text{K}$, $p = 1 \, \text{atm}$, has composition, 
$\boldsymbol{Y}(\text{CH}_3\text{OCH}_3 \, : \, \text{O}_2 \, : \, \text{N}_2) = ( 0.0726 \, : \, 0.2160 \, : \, 0.7114 )$, 
obtained from the one-dimensional solution computed with the PREMIX package \citep{kee1985premix}. We arbitrarily set the inlet velocity to $5 \, \text{cm}.\text{s}^{-1}$.
The initial pressure in the domain is constant, equal to $1 \, \text{atm}$. The initial mixture composition is also constant, equal to the inlet stream composition given above. 
Finally, the initial temperature field is set to
\begin{equation}
T(x) = T_{\text{min}} + (T_{\text{max}} - T_{\text{min}}) \exp \bigg[- \frac{1}{2} \bigg( \displaystyle \frac{x - d/2}{\kappa} \bigg)^2 \bigg] \qquad x \in [0, 0.6],
\label{temperature}
\end{equation}
with $T_{\text{min}} = 298 \, \text{K}$, $T_{\text{max}} = 1615 \, \text{K}$, and $\kappa = 0.0275$. We simulate this test case 
using a constant time step size subdivided 
with three Gauss-Lobatto nodes. 
We note that with this small number of SDC nodes ($M = 2$), we use up to two processors 
to perform the $2\nu$ pipelined tasks employed by CISDCQ-$\nu$. That is, CISDCQ-2, CISDCQ-3, and CISDCQ-4 all allow to use two processors working
in parallel.

We first make sure that the CISDCQ scheme has the same order of accuracy as MISDCQ upon refinement in space and time. This is shown in Table 
\oldref{dme_problem_batstoi_convergence_accuracy_of_amisdcQ_2}, in which we perform a refinement study for a fixed advective CFL number of 
$\sigma \approx 0.65$, starting from $n_x = 128$ and $\Delta t = 4 \times 10^{-6} \, \text{s}$. The system is evolved for $8.8 \times 10^{-5} \, \text{s}$. 
The reference solution is obtained with MISDCQ with $n_x = 1024$ and $\Delta t = 5 \times 10^{-7} \, \text{s}$. The schemes are tested with eight 
sweeps per time step, as in \cite{pazner2016high}. Table \oldref{dme_problem_batstoi_convergence_accuracy_of_amisdcQ_2} demonstrates that MISDCQ 
and CISDCQ-2 have the same order of accuracy for the main thermodynamic variables as well as for the mass fraction of dioxygen. This shows that 
CISDCQ is as accurate as MISDCQ even when a relatively small number of nested iterations on $\ell$ is used. For this extremely stiff example, the 
MISDC sweeps fail to converge. The MISDC scheme is therefore significantly less accurate than MISDCQ and CISDCQ-2. 
The same is true for CISDCQ-1, whose results are therefore not included in Table \oldref{dme_problem_batstoi_convergence_accuracy_of_amisdcQ_2}.
We note that for the initial setup \ref{temperature} and the levels 
of spatio-temporal refinement used here, none of the schemes is in the asymptotic convergence regime, and more sweeps would be required 
to reach the formal order of accuracy of the schemes. This is why MISDCQ and CISDCQ-2 do not achieve fourth-order accuracy even though the spatial and 
temporal discretizations are formally fourth-order accurate. 

\begin{table}[!ht]
\centering
 \scalebox {1}{
         \begin{tabular}{ccccccc}\toprule 
                                    & Variable                        & $L^{128}_1$ & $r^{128/256}$  & $L^{256}_1$ & $r^{256/512}$ & $L^{512}_1$ \\\toprule  
        \multirow{4}{*}{MISDC }     & $Y_{\text{O}_2}$                   &  2.62E-03  &  1.05       & 1.26E-03   &   0.34      & 1.00E-03   \\
        \multirow{4}{*}{}           & $\rho$                          & 1.15E-05 &  0.87  & 6.29E-06  &  0.31  & 5.07E-06   \\
        \multirow{4}{*}{}           & $\rho h$                        & 6.99E+04 &  0.86  & 3.84E+04  &  0.31  & 3.11E+04   \\
        \multirow{4}{*}{}           & $T$                             & 2.86E+01 &  0.92  & 1.51E+01  &  0.31  & 1.21E+01   \\ \cmidrule(l){1-7}

        \multirow{4}{*}{MISDCQ}  & $Y_{\text{O}_2}$                   & 8.91E-04 &  2.42  & 1.67E-04  &   2.27  & 3.45E-05   \\
        \multirow{4}{*}{}           & $\rho$                          & 4.61E-06 &  2.33  & 9.19E-07  &   2.25 & 1.93E-07   \\
        \multirow{4}{*}{}           & $\rho h$                        & 2.78E+04 &  2.33  & 5.51E+03  &   2.26 & 1.15E+03   \\
        \multirow{4}{*}{}           & $T$                             & 1.07E+01 &  2.37  &  2.08E+00 &   2.27 & 4.32E-01   \\ \cmidrule(l){1-7}

        \multirow{4}{*}{CISDCQ-2}  & $Y_{\text{O}_2}$                   & 8.89E-04 &  2.42  & 1.66E-04  &   2.27  & 3.44E-05   \\
        \multirow{4}{*}{}           & $\rho$                          & 4.60E-06 &  2.32  & 9.18E-07  &   2.25 & 1.92E-07   \\
        \multirow{4}{*}{}           & $\rho h$                        & 2.77E+04 &  2.33  & 5.50E+03  &   2.25 & 1.15E+03   \\
        \multirow{4}{*}{}           & $T$                             & 1.07E+01 &  2.36  &  2.07E+00 &   2.26 & 4.32E-01   \\
           \bottomrule 
         \end{tabular}}
\caption{\label{dme_problem_batstoi_convergence_accuracy_of_amisdcQ_2}
$L_1$-norm of the error with respect to the reference solution as a function of the level of space-time refinement for the DME flame simulation. 
For brevity, this table only reports the mass fraction $Y_{\text{O}_2}$, but similar results are observed for the other species in the system.
The convergence rate is defined as $r^{c/f} = \log_2 ( L_1^{c} / L_1^{f} )$.
}
\end{table}

\begin{figure}[ht!]
\centering
\subfigure[]{
\begin{tikzpicture}
\node[anchor=south west,inner sep=0] at (0,0){\includegraphics[scale=0.45]{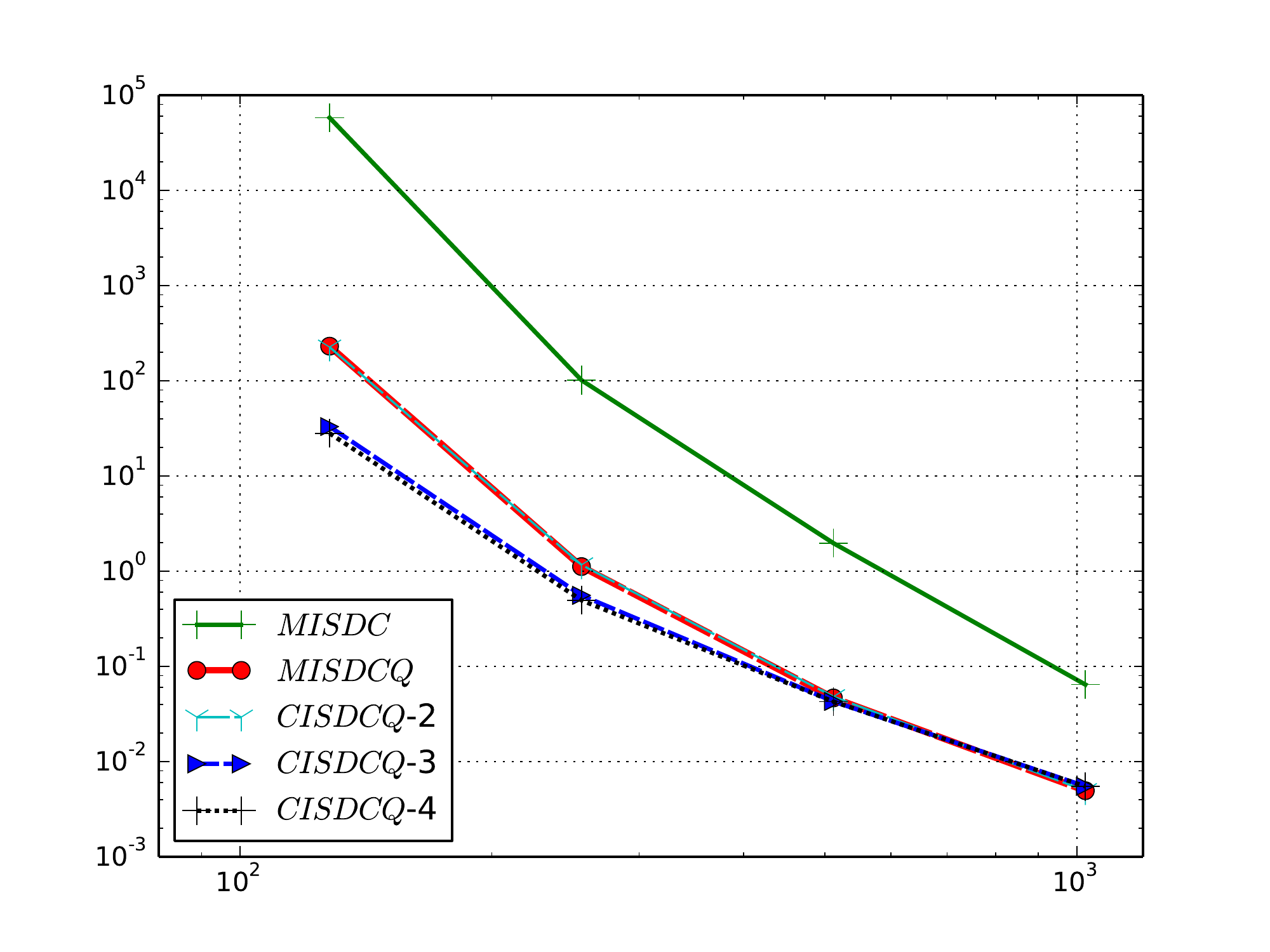}};
\node (ib_1) at (4.8,0.3) {\large $n$};
\node[rotate=90] (ib_1) at (0.175,3.4) {\large $|| p_{\textit{EOS}} - p_0 ||_{\infty}$};
\end{tikzpicture}
}
\vspace{-0.4cm}
\caption{\label{fig:dme_problem_batstoi_accuracy_convergence_for_thermodynamic_drift} 
Convergence of the $L_{\infty}$-norm of the thermodynamic drift as a function of the level of space-time refinement at the end of the DME flame simulation. 
}
\end{figure}

In Fig.~\oldref{fig:dme_problem_batstoi_accuracy_convergence_for_thermodynamic_drift}, we evaluate the magnitude of the thermodynamic drift 
as a function of 
the level of space-time refinement for the different schemes. To evaluate the magnitude of the drift, we use the $L_{\infty}$-norm 
\begin{equation}
|| p_{\textit{EOS}} - p_0 ||_{\infty} = \max_{i \in \{ 1, \dots, n_x \} }\big( | p_{\textit{EOS}}(x_i,t_{\textit{final}}) - p_0 | \big),
\end{equation}
where the thermodynamic pressure, $p_{\textit{EOS}}$, is defined in \ref{p_eos_definition}.
We still use a constant time step size corresponding to a CFL number $\sigma \approx 0.65$ and eight sweeps per time step. For this problem, 
MISDC yields a thermodynamic drift that is multiple orders of magnitude larger than that of the MISDCQ and CISDCQ schemes. 
MISDCQ and CISDCQ-2 achieve a similar thermodynamic drift for all levels of refinement, which is consistent with the results 
of Table \oldref{dme_problem_batstoi_convergence_accuracy_of_amisdcQ_2}. Increasing the number of nested iterations on $\ell$ yields a 
smaller norm of the drift. Specifically, CISDCQ-3 and CISDCQ-4 produce mass and energy fields that are more consistent with the 
EOS for the coarse levels of refinement ($n_x = 128, 256$).

\begin{figure}[ht!]
\centering
\subfigure[]{
\begin{tikzpicture}
\node[anchor=south west,inner sep=0] at (0,0){\includegraphics[scale=0.3]{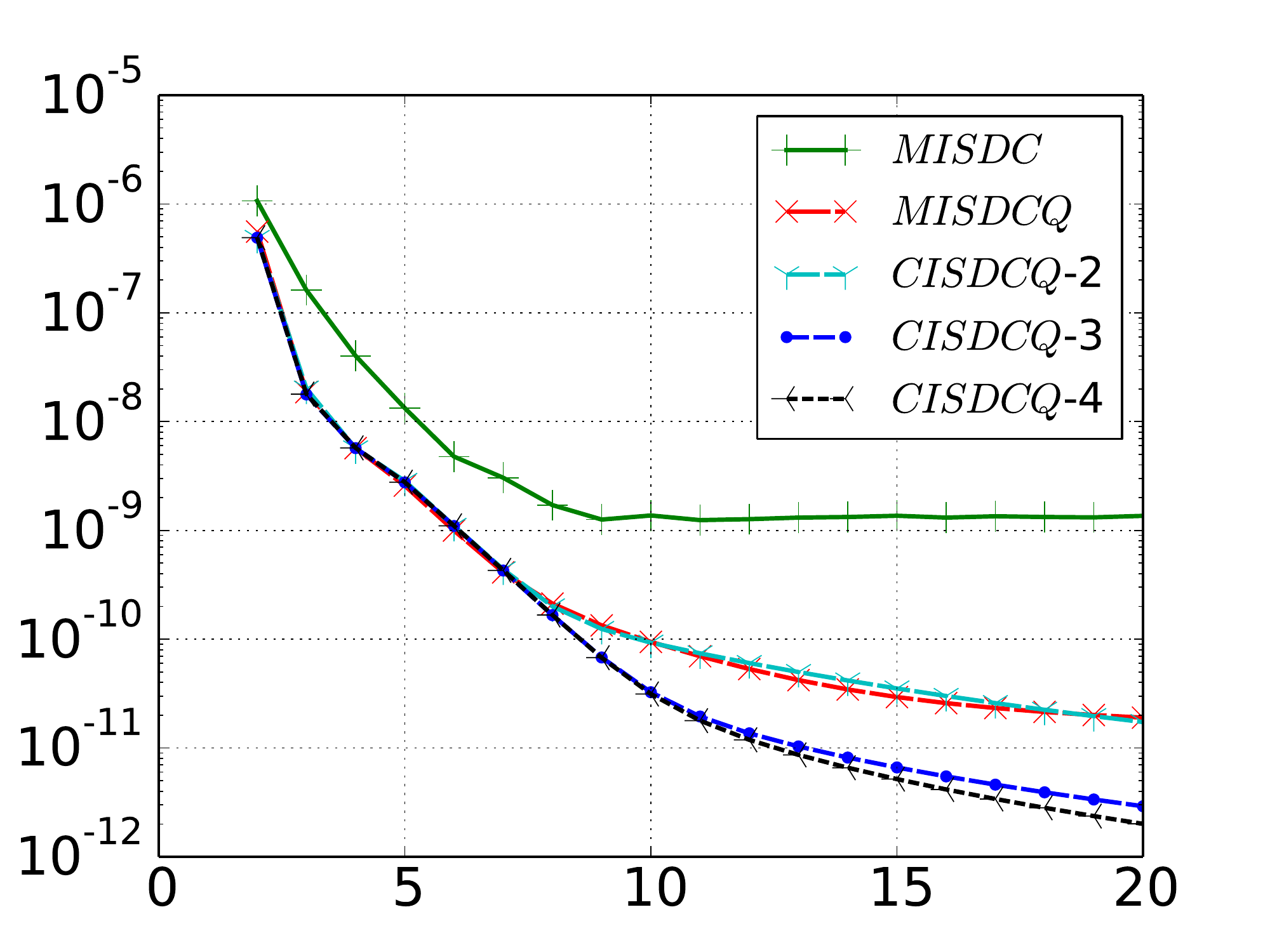}};
\node (ib_1) at (3.1,0.) {\scriptsize Number of sweeps};
\node[rotate=90] (ib_1) at (-0.16,2.25) {\scriptsize $|| \rho^{M,(k+1)} - \rho^{M,(k)} ||_1$};
\end{tikzpicture}
\label{fig:dme_problem_batstoi_convergence_as_fct_number_of_sweeps_density}
}
\subfigure[]{
\begin{tikzpicture}
\node[anchor=south west,inner sep=0] at (0,0){\includegraphics[scale=0.3]{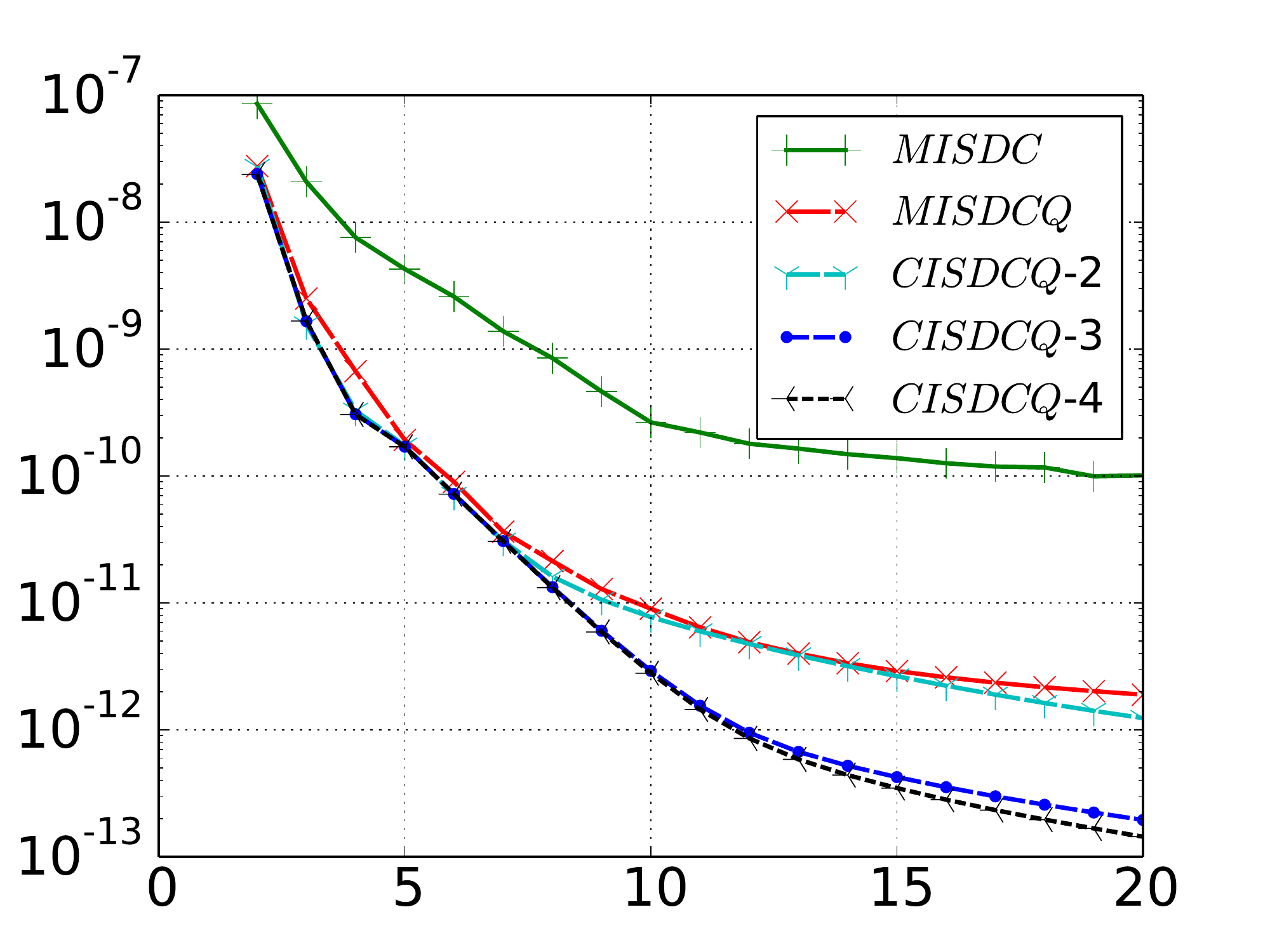}};
\node (ib_1) at (3.1,0.) {\scriptsize Number of sweeps};
\node[rotate=90] (ib_1) at (-0.16,2.25) {\scriptsize $|| T^{M,(k+1)} - T^{M,(k)} ||_1$};
\end{tikzpicture}
\label{fig:dme_problem_batstoi_convergence_as_fct_number_of_sweeps_temperature}
}
\subfigure[]{
\begin{tikzpicture}
\node[anchor=south west,inner sep=0] at (0,0){\includegraphics[scale=0.3]{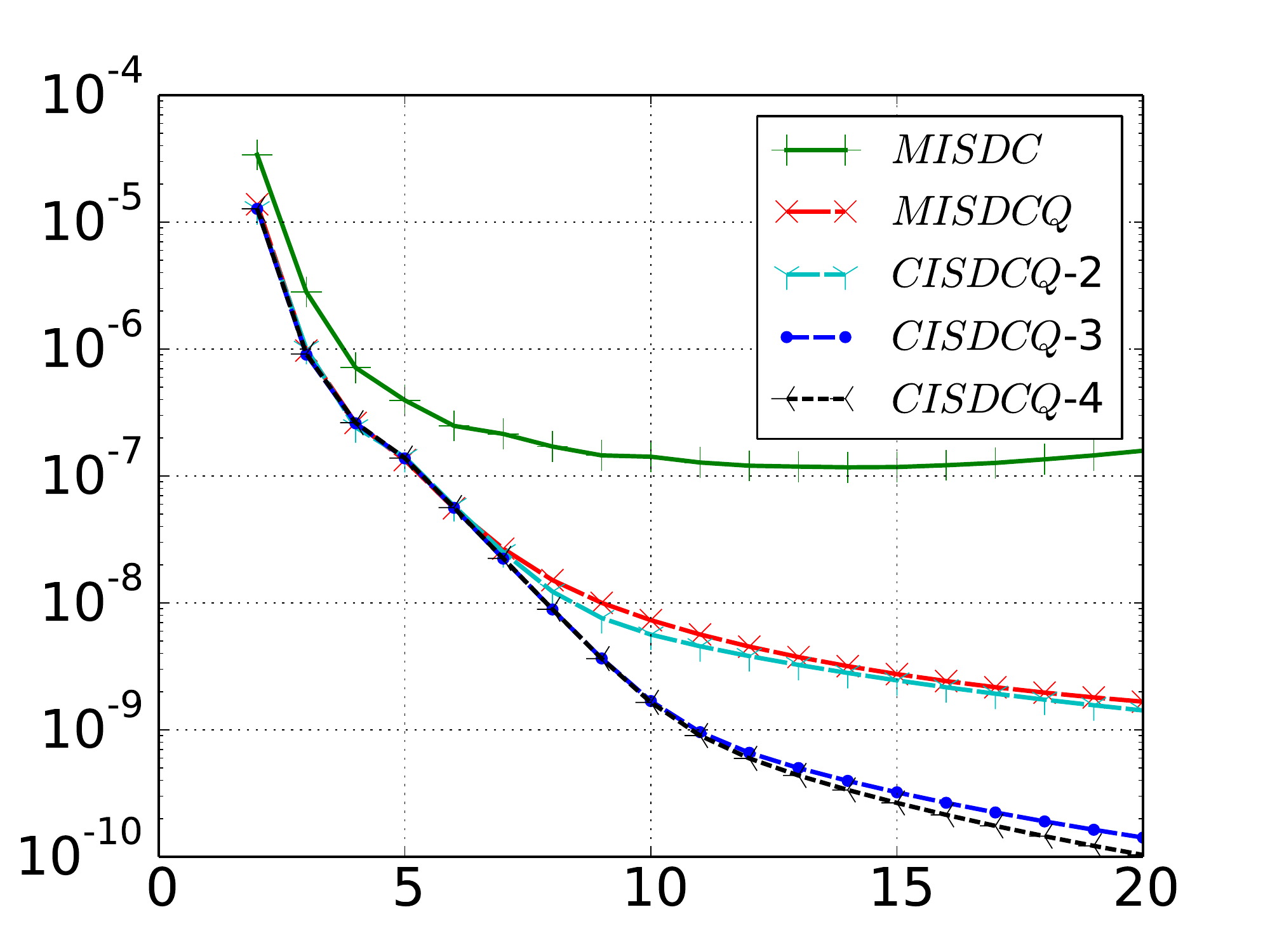}};
\node (ib_1) at (3.1,0.) {\scriptsize Number of sweeps};
\node[rotate=90] (ib_1) at (-0.16,2.25) {\scriptsize $|| Y_{\text{O}_2}^{M,(k+1)} - Y_{\text{O}_2}^{M,(k)} ||_1$};
\end{tikzpicture}
\label{fig:dme_problem_batstoi_convergence_as_fct_number_of_sweeps_yo2}
}
\subfigure[]{
\begin{tikzpicture}
\node[anchor=south west,inner sep=0] at (0,0){\includegraphics[scale=0.3]{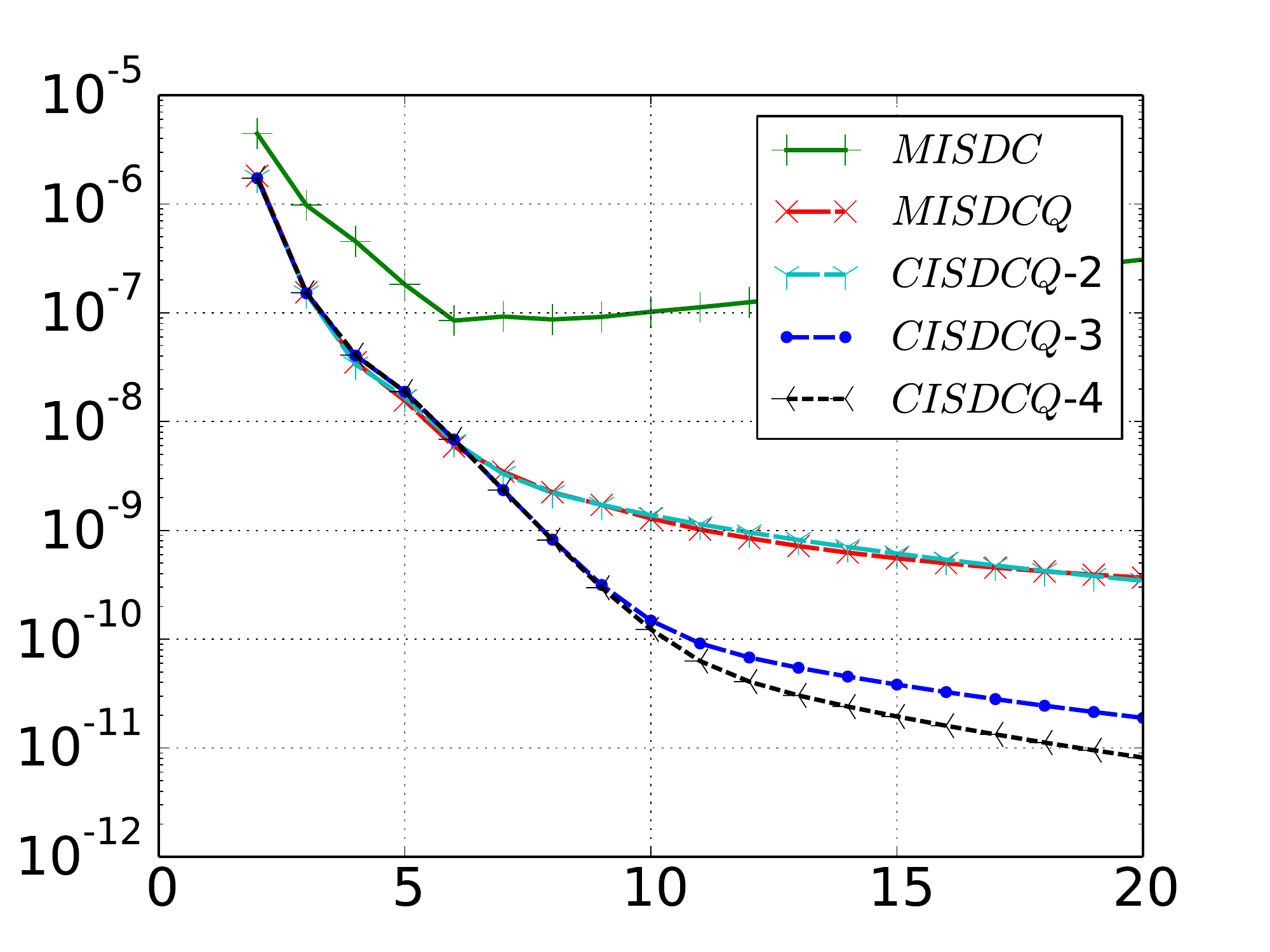}};
\node (ib_1) at (3.1,0.) {\scriptsize Number of sweeps};
\node[rotate=90] (ib_1) at (-0.16,2.25) {\scriptsize $|| Y_{\text{CH}_3\text{OCH}_2\text{O}_2}^{M,(k+1)} - Y_{\text{CH}_3\text{OCH}_2\text{O}_2}^{M,(k)} ||_1$};
\end{tikzpicture}
\label{fig:dme_problem_batstoi_convergence_as_fct_number_of_sweeps_ych3och2o2}
}
\vspace{-0.4cm}
\caption{\label{fig:dme_problem_batstoi_convergence_as_fct_number_of_sweeps} Convergence of the $L_1$-norm of the difference between two consecutive 
iterates for the last time step of the 
DME flame simulation performed with $n_x = 512$ control volumes. The density, temperature, and mass fractions $Y_{\text{O}_2}$ and 
$Y_{\text{CH}_3\text{OCH}_2\text{O}_2}$ are in Figs. 
\oldref{fig:dme_problem_batstoi_convergence_as_fct_number_of_sweeps_density}, 
\oldref{fig:dme_problem_batstoi_convergence_as_fct_number_of_sweeps_temperature}, 
\oldref{fig:dme_problem_batstoi_convergence_as_fct_number_of_sweeps_yo2}, 
and \oldref{fig:dme_problem_batstoi_convergence_as_fct_number_of_sweeps_ych3och2o2}, respectively.
}
\end{figure}

Finally, we assess the convergence of the correction process in the last time step of the DME flame simulation, for a level of 
refinement of $n_x = 512$ and $\Delta t = 10^{-6} \, \text{s}$, still yielding a CFL number $\sigma \approx 0.65$. The system is now evolved for 
$10^{-4} \, \text{s}$. We compute the $L_1$-norm of the difference between two consecutive iterates for the key variables as a function of the 
number of sweeps. The result, shown in Fig.~\oldref{fig:dme_problem_batstoi_convergence_as_fct_number_of_sweeps}, is consistent with those of 
Sections \oldref{subsection_convergence_analysis} and \oldref{subsubsection_nonlinear_test_pde}. 
The MISDC sweeps fail to converge to a fixed point solution because of the stiffness of the DME chemistry mechanism. This is 
particularly visible after sweep 15 in Figs.~\oldref{fig:dme_problem_batstoi_convergence_as_fct_number_of_sweeps_yo2} and 
\oldref{fig:dme_problem_batstoi_convergence_as_fct_number_of_sweeps_ych3och2o2}. The other schemes achieve convergence at different rates. 
CISDCQ-3 and CISDCQ-4 have practically converged at sweep 20. But, for CISDCQ-2 and MISDCQ, the convergence rate is relatively slower and more
than 40 sweeps are required to reach convergence (not shown in Fig.~\oldref{fig:dme_problem_batstoi_convergence_as_fct_number_of_sweeps}).

Using the same methodology as in Section \oldref{subsubsection_nonlinear_test_pde}, we compare the respective computational costs 
of the schemes after 20 and after 25 sweeps for the density (Fig.~\oldref{fig:dme_problem_batstoi_convergence_as_fct_number_of_sweeps}). The 
computational cost ratios $R_{\nu}$ are in Table \oldref{parallel_speedup_nonlinear_dme}. 
We find that the convergence of the sweeps to a fixed-point solution is too slow to achieve a computational cost reduction
when only two nested iterations are used ($R_{\nu = 2} = 0.8$ after 25 sweeps). But, with three iterations, the significantly 
faster convergence of the sweeps leads to a reduction in computational cost ($R_{\nu = 3} = 1.6$ after 25 sweeps). 
Finally, using four nested iterations does not improve the ratio $R_{\nu}$ due to the increase in 
the cost of the sweep ($S_{\nu = 4} = 1.2$ after 25 sweeps). As stated before, an adaptive method to vary the number of nested iterations 
on $\ell$ at each sweep would increase, up to a certain extent, the computational cost reduction obtained with CISDCQ.


\begin{table}[!ht]
\centering
 \scalebox {1}{
         \begin{tabular}{ccccc}
           \\ \toprule 
           $N_M$ & $R_{\nu = 2}$ & $R_{\nu = 3}$ & $R_{\nu = 4}$   \\\toprule
           20 sweeps &   0.8      &    1.4      &   1.1      \\ 
          25 sweeps &   0.8      &    1.6      &   1.2       \\ 
           \bottomrule 
         \end{tabular}}
\caption{\label{parallel_speedup_nonlinear_dme}
Ratio of the computational cost of MISDCQ over that of CISDCQ-$\nu$ for the DME flame simulation after 20 sweeps and after 25 sweeps.
}
\end{table}

\section{\label{section_conclusions}Conclusions}

In this work, we consider the time integration of the advection-diffusion-reaction systems arising from the spatial discretization of the equations 
governing low-Mach number combustion with complex chemistry. The problem is advanced in time on the scale of the relatively slow advection process, 
which makes it challenging to efficiently capture the nonlinear coupling with the faster diffusion and reaction processes. In the serial MISDC algorithm
of \cite{pazner2016high}, this entails solving expensive implicit diffusion and reaction systems at each substep. We first present an improved multi-implicit
scheme, MISDCQ, whose sweeps retain excellent convergence properties on stiff problems. Then, to reduce the time-to-solution, we modify MISDCQ to design 
a parallel-across-the-method integration scheme based on spectral deferred corrections. The scheme, referred to as Concurrent Implicit SDCQ (CISDCQ), 
is obtained by decoupling the diffusion step from the reaction step in serial MISDCQ. This allows 
concurrent implicit diffusion and reaction solves performed by different processors to reduce the time-to-solution. 

We study the accuracy, stability, convergence rate, and theoretical computational cost of the proposed parallel scheme. The linear convergence analysis 
demonstrates that CISDCQ is stable for stiff problems and that the CISDCQ iterative correction process can converge at a faster rate than that of the 
standard serial MISDC scheme. We propose a numerical methodology to apply CISDCQ to the discretized low-Mach number equation set with complex chemistry. 
The numerical results -- including a stiff premixed flame simulation -- confirm the findings of the convergence analysis and demonstrate the robustness 
and parallel efficiency of the new scheme. 

Future work includes extending the method presented here to multidimensional simulations. In multiple dimensions, the efficiency 
of the scheme could be improved by exploiting
local variations in temporal stiffness.  We plan to exploit local variabilities by combining SDC-based 
algorithms (CISDC, and Multi-Level SDC of \cite{speck2015multi}) with block-structured 
adaptive refinement techniques, which will be used to concentrate the computational effort
-- small time steps and finer spatial meshing -- near local structures in the flow.  This will
allow us to simultaneously deal with variability in temporal and spatial resolution requirements,
which typically are coincident in combustion problems.
%
A further improvement 
would consist in amortizing the cost of the fine sweeps by performing them in parallel, as in the Parallel Full 
Approximation Scheme in Space and in Time (PFASST) scheme \citep{emmett2012toward}.

As a final note, we mention an important but somewhat indirect benefit of our proposed concurrent update scheme. 
Linear solvers play a key role in the implementation of the low-Mach projection algorithm in multiple dimensions, both 
for the implicit diffusion and for the velocity projection operators.  However, linear solvers are notoriously difficult 
to scale up to the large processor counts typical of modern massively parallel computing hardware.  But because each of 
the algorithmic components in our scheme are managed by a partitioned subset of the processors working in an ``embarrassingly parallel'' 
mode with the other partitions, the requirements of the combined algorithm are reduced directly, by a factor equal
to the number of available concurrent tasks.  In the DME model discussed in this work for example, the implicit diffusion 
of 39 chemical species can be performed independently, {\em and}\ independent of the chemistry integration.  
While the parallelization of the diffusion step over the species can be implemented in the original MISDC scheme,
the results in this paper verify that we have developed an algorithm that will allow unique access to a new axis
of parallelism.  Whether it can be exploited effectively depends on the cost of communicating the subproblems across 
the machine, relative to the gains afforded by the smaller pool of compute processors.  Exploration of that trade-off will be the subject 
of future work.

\section*{Funding}

This research was supported by the Applied Mathematics Program of the U.S. DOE Office of Advanced Scientific Computing 
Research under contract number DE-AC02-05CH11231. Some computations used the resources of the National Energy
Research Scientific Computing Center, a DOE Office of Science User Facility supported by the Office of Science of the U.S.
DOE under Contract No. DE-AC02-05CH11231.

\appendix

\section{\label{section_derivation_of_the_amisdcq_iteration_matrix_g}Derivation of the CISDCQ iteration matrix}

We prove by induction on $\ell$ and $m$ that the matrix form of the CISDCQ update equations 
\ref{update_equation_amisdcq_discrete_form_D}-\ref{update_equation_amisdcq_discrete_form_R} leads to 
a relationship between $\boldsymbol{\Phi}^{(k+1,\ell+1)}_{m+1}$ and $\boldsymbol{\Phi}^{(k)}$ in the form of \ref{phi_kp1_as_fct_phi_k_phi_m_phi_l}.
We first introduce the matrices
\begin{align}
\boldsymbol{A}_{\kappa,m+1} &:= \boldsymbol{I} - \kappa \Delta t \boldsymbol{D}_{m+1} \boldsymbol{\tilde{Q}}^I, \qquad \kappa = \{ d, r \}
\end{align}
where $\boldsymbol{I} \in \mathbb{R}^{M \times M}$ is the identity matrix, and $\boldsymbol{D}_{m+1} \in \mathbb{R}^{M \times M}$  
contains only one non-zero entry,  equal to one, in slot $(m+1,m+1)$, with the convention that $\boldsymbol{D}_0 = \boldsymbol{0}$. 
We first consider the first nested iteration ($\ell = 0$). 
With the initialization procedure of \ref{initial_lagged_values_A} to \ref{initial_lagged_values_R}, the diffusion update 
\ref{update_equation_amisdcq_discrete_form_D} in matrix form reads
\begin{align}
\boldsymbol{D}_{m+1} \boldsymbol{\Phi}^{(k+1,1)}_{AD,m+1} &= \boldsymbol{D}_{m+1} \big[ \phi^0 \boldsymbol{1} + s \Delta t \phi^0 \boldsymbol{q} 
+ \Delta t  \big( s \boldsymbol{\tilde{Q}} - a \boldsymbol{\tilde{Q}}^E - (d+r) \boldsymbol{\tilde{Q}}^I \big) \boldsymbol{\Phi}^{(k)}
                                                                             \big] \nonumber \\
         &+ a \Delta t \boldsymbol{D}_{m+1} \boldsymbol{\tilde{Q}}^E \big( (\boldsymbol{I} - \boldsymbol{D}_m) \boldsymbol{\Phi}^{(k+1,1)}_m
                                                                                           + \boldsymbol{D}_m \boldsymbol{\Phi}^{(k+1,1)}_{AD,m+1}
                                                                                                                     \big) \nonumber \\                                  &+ d \Delta t \boldsymbol{D}_{m+1} \boldsymbol{\tilde{Q}}^I \big( (\boldsymbol{I} - \boldsymbol{D}_m) \boldsymbol{\Phi}^{(k+1,1)}_m                                                                                                             + \boldsymbol{D}_m \boldsymbol{\Phi}^{(k+1,1)}_{AD,m+1}                                                                                                       + \boldsymbol{D}_{m+1} \boldsymbol{\Phi}^{(k+1,1)}_{AD,m+1} 
                                                                                                                    \big) \nonumber \\
         &+ r \Delta t \boldsymbol{D}_{m+1} \boldsymbol{\tilde{Q}}^I \big( (\boldsymbol{I} - \boldsymbol{D}_m) \boldsymbol{\Phi}^{(k+1,1)}_m
                                                                                           + (\boldsymbol{D}_m + \boldsymbol{D}_{m+1}) \boldsymbol{\Phi}^{(k+1,0)}_{m+1} 
                                                                                                                    \big).
\label{matrix_form_of_amisdc_Q_update_D}
\end{align}
By construction of the vectors $\boldsymbol{\Phi}^{(k+1,\ell+1)}_{m+1}$ and $\boldsymbol{\Phi}^{(k+1,\ell+1)}_{AD,m+1}$, only the $(m+1)^{\text{th}}$ entry is modified at the
$(m+1)^{\text{th}}$ substep. Therefore, the following equalities hold
\begin{align}
(\boldsymbol{I} - \boldsymbol{D}_{m+1}) \boldsymbol{\Phi}^{(k+1,\ell+1)}_{AD,m+1} &= \boldsymbol{\Phi}^{(k+1,\ell+1)}_{AD,m}, \label{constraint_ad} \\
(\boldsymbol{I} - \boldsymbol{D}_{m+1}) \boldsymbol{\Phi}^{(k+1,\ell+1)}_{m+1}   &= \boldsymbol{\Phi}^{(k+1,\ell+1)}_{m}. \label{constraint}
\end{align}
Summing \ref{matrix_form_of_amisdc_Q_update_D} and \ref{constraint_ad}, and then multiplying the resulting equation by $\boldsymbol{A}^{-1}_{d,m+1}$,
we obtain
\begin{align}
\boldsymbol{\Phi}^{(k+1,1)}_{AD,m+1} 
&= \boldsymbol{A}^{-1}_{d,m+1} \boldsymbol{D}_{m+1}  \big[ \phi^0 \boldsymbol{1} + s \Delta t \phi^0 \boldsymbol{q} +
\Delta t \big( s \boldsymbol{\tilde{Q}} - a \boldsymbol{\tilde{Q}}^E - (d + r) \boldsymbol{\tilde{Q}}^I \big) \boldsymbol{\Phi}^{(k)} \big] \nonumber \\
&+ \boldsymbol{A}^{-1}_{d,m+1} \big( \boldsymbol{B}_{m+1} \boldsymbol{\Phi}^{(k+1,1)}_{m} + \boldsymbol{C}_{m+1} \boldsymbol{\Phi}^{(k+1,1)}_{AD,m} + \boldsymbol{E}_{m+1} \boldsymbol{\Phi}^{(k+1,0)}_{m+1} \big),
\label{expression_for_phi_ad}
\end{align}
where the matrices $\boldsymbol{B}_{m+1}$, $\boldsymbol{C}_{m+1}$, and $\boldsymbol{E}_{m+1}$ depend on the 
matrices $\boldsymbol{\tilde{Q}}$, $\boldsymbol{\tilde{Q}}^I$, and $\boldsymbol{\tilde{Q}}^E$, and on the scalars $m$, $\Delta t$, $a$, $d$, 
and $r$. The expression of these matrices is omitted for brevity. Next, we consider the reaction update.
In matrix form, the reaction update \ref{update_equation_amisdcq_discrete_form_R} is 
\begin{align}
\boldsymbol{D}_{m+1} \boldsymbol{\Phi}^{(k+1,1)}_{m+1}&= \boldsymbol{D}_{m+1} \boldsymbol{\Phi}^{(k+1, 1)}_{AD,m+1} 
+ r \Delta t \boldsymbol{D}_{m+1} \boldsymbol{\tilde{Q}}^I ( \boldsymbol{D}_m + \boldsymbol{D}_{m+1} ) \big( \boldsymbol{\Phi}^{(k+1,1)}_{m+1} - \boldsymbol{\Phi}^{(k+1, 0)}_{m+1} \big),
\end{align}
which gives, after using \ref{expression_for_phi_ad} to eliminate $\boldsymbol{\Phi}^{(k+1,1)}_{AD,m+1}$,
\begin{align}
\boldsymbol{D}_{m+1} \boldsymbol{\Phi}^{(k+1,1)}_{m+1} &= \boldsymbol{D}_{m+1} \boldsymbol{A}^{-1}_{d,m+1} \boldsymbol{D}_{m+1}  \big[ \phi^0 \boldsymbol{1} + s \Delta t \phi^0 \boldsymbol{q} +
\Delta t \big( s \boldsymbol{\tilde{Q}} - a \boldsymbol{\tilde{Q}}^E - (d + r) \boldsymbol{\tilde{Q}}^I \big) \boldsymbol{\Phi}^{(k)} \big] \nonumber \\
&+ \boldsymbol{D}_{m+1} \boldsymbol{A}^{-1}_{d,m+1} \big( \boldsymbol{B}_{m+1} \boldsymbol{\Phi}^{(k+1,1)}_{m} + \boldsymbol{C}_{m+1} \boldsymbol{\Phi}^{(k+1,1)}_{AD,m} + \boldsymbol{E}_{m+1} \boldsymbol{\Phi}^{(k+1,0)}_{m} \big) \nonumber \\
                                                  &+ r \Delta t \boldsymbol{D}_{m+1} \boldsymbol{\tilde{Q}}^I ( \boldsymbol{D}_m + \boldsymbol{D}_{m+1} ) \big( \boldsymbol{\Phi}^{(k+1,1)}_{m+1} - \boldsymbol{\Phi}^{(k+1, 0)}_{m+1} \big).
\label{matrix_form_of_amisdc_Q_update_R}
\end{align}
Summing \ref{constraint} and \ref{matrix_form_of_amisdc_Q_update_R}, and then multiplying the equation by $\boldsymbol{A}^{-1}_{r,m+1}$ leads to
\begin{align}
\boldsymbol{\Phi}^{(k+1,1)}_{m+1}&= \boldsymbol{A}^{-1}_{r,m+1} \boldsymbol{D}_{m+1} \boldsymbol{A}^{-1}_{d,m+1} \boldsymbol{D}_{m+1}  \big[ \phi^0 \boldsymbol{1} + s \Delta t \phi^0 \boldsymbol{q} +
\Delta t \big( s \boldsymbol{\tilde{Q}} - a \boldsymbol{\tilde{Q}}^E - (d + r) \boldsymbol{\tilde{Q}}^I \big) \boldsymbol{\Phi}^{(k)} \big] \nonumber \\
&+ \boldsymbol{A}^{-1}_{r,m+1} \big( \boldsymbol{R}_{m+1} \boldsymbol{\Phi}^{(k+1,1)}_m 
                                 + \boldsymbol{S}_{m+1} \boldsymbol{\Phi}^{(k+1, 1)}_{AD, m} 
                                 + \boldsymbol{T}_{m+1} \boldsymbol{\Phi}^{(k+1,0)}_{m+1} \big). 
\label{expression_for_phi}
\end{align}
The expressions of $\boldsymbol{R}_{m+1}$, $\boldsymbol{S}_{m+1}$, and $\boldsymbol{T}_{m+1}$ are obtained using \ref{matrix_form_of_amisdc_Q_update_R}. They depend on
the matrices 
$\boldsymbol{\tilde{Q}}$, $\boldsymbol{\tilde{Q}}^I$, and $\boldsymbol{\tilde{Q}}^E$, 
and on $m$, $\Delta t$, $a$, $d$, and $r$. 
Considering that $( \boldsymbol{D}_m + \boldsymbol{D}_{m+1} ) \boldsymbol{\Phi}^{(k+1,0)}_{m+1} = ( \boldsymbol{D}_m + \boldsymbol{D}_{m+1} ) \boldsymbol{\Phi}^{(k)}$, we can 
use \ref{expression_for_phi_ad}-\ref{expression_for_phi} to show by induction on $m$ that \ref{phi_kp1_as_fct_phi_k_phi_m_phi_l} holds for $\ell = 0$
and derive the expressions for $\boldsymbol{M}^{(1)}_{1,m+1}$, $\boldsymbol{M}^{(1)}_{2,m+1}$, $\boldsymbol{N}^{(1)}_{1,m+1}$, and $\boldsymbol{N}^{(1)}_{2,m+1}$.
The proof for $\ell \geq 1$ is analogous and is omitted here.



\bibliography{biblio}

\begin{thebibliography}{}

\bibitem[Bansal et~al., 2015]{bansal2015direct}
Bansal, G., Mascarenhas, A., and Chen, J.~H. (2015).
\newblock Direct numerical simulations of autoignition in stratified
  dimethyl-ether ({D}{M}{E})/air turbulent mixtures.
\newblock {\em Combustion and Flame}, 162(3):688--702.

\bibitem[Bourlioux et~al., 2003]{bourlioux2003high}
Bourlioux, A., Layton, A.~T., and Minion, M.~L. (2003).
\newblock High-order multi-implicit spectral deferred correction methods for
  problems of reactive flow.
\newblock {\em Journal of Computational Physics}, 189(2):651--675.

\bibitem[Burrage, 1997]{burrage1997parallel}
Burrage, K. (1997).
\newblock Parallel methods for {O}{D}{E}s.
\newblock {\em Advances in Computational Mathematics}, 7(1-2):1--31.

\bibitem[Butcher, 1997]{butcher1997order}
Butcher, J.~C. (1997).
\newblock Order and stability of parallel methods for stiff problems.
\newblock {\em Advances in Computational Mathematics}, 7(1):79--96.

\bibitem[Christlieb et~al., 2012]{christlieb2012parallel}
Christlieb, A.~J., Haynes, R.~D., and Ong, B.~W. (2012).
\newblock A parallel space-time algorithm.
\newblock {\em SIAM Journal on Scientific Computing}, 34(5):C233--C248.

\bibitem[Christlieb et~al., 2010]{christlieb2010parallel}
Christlieb, A.~J., Macdonald, C.~B., and Ong, B.~W. (2010).
\newblock Parallel high-order integrators.
\newblock {\em SIAM Journal on Scientific Computing}, 32(2):818--835.

\bibitem[Christlieb et~al., 2009]{christlieb2009comments}
Christlieb, A.~J., Ong, B.~W., and Qiu, J.-M. (2009).
\newblock Comments on high-order integrators embedded within integral deferred
  correction methods.
\newblock {\em Communications in Applied Mathematics and Computational
  Science}, 4(1):27--56.

\bibitem[Day and Bell, 2000]{day2000numerical}
Day, M.~S. and Bell, J.~B. (2000).
\newblock Numerical simulation of laminar reacting flows with complex
  chemistry.
\newblock {\em Combustion Theory and Modelling}, 4(4):535--556.

\bibitem[Duarte et~al., 2013]{duarte2013time}
Duarte, M., Descombes, S., Tenaud, C., Candel, S., and Massot, M. (2013).
\newblock Time--space adaptive numerical methods for the simulation of
  combustion fronts.
\newblock {\em Combustion and Flame}, 160(6):1083--1101.

\bibitem[Dutt et~al., 2000]{dutt2000spectral}
Dutt, A., Greengard, L., and Rokhlin, V. (2000).
\newblock Spectral deferred correction methods for ordinary differential
  equations.
\newblock {\em BIT Numerical Mathematics}, 40(2):241--266.

\bibitem[Emmett and Minion, 2012]{emmett2012toward}
Emmett, M. and Minion, M.~L. (2012).
\newblock Toward an efficient parallel in time method for partial differential
  equations.
\newblock {\em Communications in Applied Mathematics and Computational
  Science}, 7(1):105--132.

\bibitem[Falgout et~al., 2014]{falgout2014parallel}
Falgout, R.~D., Friedhoff, S., Kolev, T.~V., MacLachlan, S.~P., and Schroder,
  J.~B. (2014).
\newblock Parallel time integration with multigrid.
\newblock {\em SIAM Journal on Scientific Computing}, 36(6):C635--C661.

\bibitem[Gander, 1999]{gander1999waveform}
Gander, M.~J. (1999).
\newblock A waveform relaxation algorithm with overlapping splitting for
  reaction diffusion equations.
\newblock {\em Numerical Linear Algebra with Applications}, 6(2):125--145.

\bibitem[Hagstrom and Zhou, 2007]{hagstrom2007spectral}
Hagstrom, T. and Zhou, R. (2007).
\newblock On the spectral deferred correction of splitting methods for initial
  value problems.
\newblock {\em Communications in Applied Mathematics and Computational
  Science}, 1(1):169--205.

\bibitem[Iserles and N{\o}rsett, 1990]{iserles1990theory}
Iserles, A. and N{\o}rsett, S.~P. (1990).
\newblock On the theory of parallel {R}unge-{K}utta methods.
\newblock {\em IMA Journal of Numerical Analysis}, 10(4):463--488.

\bibitem[Kee et~al., 1985]{kee1985premix}
Kee, R.~J., Grcar, J.~F., Smooke, M.~D., Miller, J.~A., and Meeks, E. (1985).
\newblock P{R}{E}{M}{I}{X}: a {F}ortran program for modeling steady laminar
  one-dimensional premixed flames.
\newblock Technical report, Sandia National Laboratory, USA.

\bibitem[Kee et~al., 1983]{kee1983fortran}
Kee, R.~J., Warnatz, J., and Miller, J.~A. (1983).
\newblock Fortran computer-code package for the evaluation of gas-phase
  viscosities, conductivities, and diffusion coefficients.
\newblock Technical report, Sandia National Laboratory, USA.

\bibitem[Knio et~al., 1999]{knio1999semi}
Knio, O.~M., Najm, H.~N., and Wyckoff, P.~S. (1999).
\newblock A semi-implicit numerical scheme for reacting flow: Ii. stiff,
  operator-split formulation.
\newblock {\em Journal of Computational Physics}, 154(2):428--467.

\bibitem[Layton and Minion, 2004]{layton2004conservative}
Layton, A.~T. and Minion, M.~L. (2004).
\newblock Conservative multi-implicit spectral deferred correction methods for
  reacting gas dynamics.
\newblock {\em Journal of Computational Physics}, 194(2):697--715.

\bibitem[Lelarasmee et~al., 1982]{1270004}
Lelarasmee, E., Ruehli, A.~E., and Sangiovanni-Vincentelli, A.~L. (1982).
\newblock The waveform relaxation method for time-domain analysis of large
  scale integrated circuits.
\newblock {\em IEEE Transactions on Computer-Aided Design of Integrated
  Circuits and Systems}, 1(3):131--145.

\bibitem[Lions et~al., 2001]{lions2001resolution}
Lions, J.-L., Maday, Y., and Turinici, G. (2001).
\newblock R{\'e}solution d'{E}{D}{P} par un sch{\'e}ma en temps parar{\'e}el.
\newblock {\em Comptes Rendus de l'Acad{\'e}mie des Sciences-Series
  I-Mathematics}, 332(7):661--668.

\bibitem[Majda and Sethian, 1985]{majda1985derivation}
Majda, A. and Sethian, J. (1985).
\newblock The derivation and numerical solution of the equations for zero
  {M}ach number combustion.
\newblock {\em Combustion Science and Technology}, 42(3-4):185--205.

\bibitem[McCorquodale and Colella, 2011]{mccorquodale2011high}
McCorquodale, P. and Colella, P. (2011).
\newblock A high-order finite-volume method for conservation laws on locally
  refined grids.
\newblock {\em Communications in Applied Mathematics and Computational
  Science}, 6(1):1--25.

\bibitem[Minion, 2003]{minion2003semi}
Minion, M.~L. (2003).
\newblock Semi-implicit spectral deferred correction methods for ordinary
  differential equations.
\newblock {\em Communications in Mathematical Sciences}, 1(3):471--500.

\bibitem[Miranker and Liniger, 1967]{miranker1967parallel}
Miranker, W.~L. and Liniger, W. (1967).
\newblock Parallel methods for the numerical integration of ordinary
  differential equations.
\newblock {\em Mathematics of Computation}, 21(99):303--320.

\bibitem[Motheau and Abraham, 2016]{motheau2016high}
Motheau, E. and Abraham, J. (2016).
\newblock A high-order numerical algorithm for {D}{N}{S} of low-mach-number
  reactive flows with detailed chemistry and quasi-spectral accuracy.
\newblock {\em Journal of Computational Physics}, 313:430--454.

\bibitem[Najm and Knio, 2005]{najm2005modeling}
Najm, H.~N. and Knio, O.~M. (2005).
\newblock Modeling low mach number reacting flow with detailed chemistry and
  transport.
\newblock {\em Journal of Scientific Computing}, 25(1-2):263.

\bibitem[Najm et~al., 1998]{najm1998semi}
Najm, H.~N., Wyckoff, P.~S., and Knio, O.~M. (1998).
\newblock A semi-implicit numerical scheme for reacting flow: I. stiff
  chemistry.
\newblock {\em Journal of Computational Physics}, 143(2):381--402.

\bibitem[Nievergelt, 1964]{nievergelt1964parallel}
Nievergelt, J. (1964).
\newblock Parallel methods for integrating ordinary differential equations.
\newblock {\em Communications of the ACM}, 7(12):731--733.

\bibitem[Nonaka et~al., 2012]{nonaka2012deferred}
Nonaka, A., Bell, J.~B., Day, M.~S., Gilet, C., Almgren, A.~S., and Minion,
  M.~L. (2012).
\newblock A deferred correction coupling strategy for low {M}ach number flow
  with complex chemistry.
\newblock {\em Combustion Theory and Modelling}, 16(6):1053--1088.

\bibitem[Pazner et~al., 2016]{pazner2016high}
Pazner, W.~E., Nonaka, A., Bell, J.~B., Day, M.~S., and Minion, M.~L. (2016).
\newblock A high-order spectral deferred correction strategy for low {M}ach
  number flow with complex chemistry.
\newblock {\em Combustion Theory and Modelling}, 20(3):521--547.

\bibitem[Pember et~al., 1998]{pember1998adaptive}
Pember, R.~B., Howell, L.~H., Bell, J.~B., Colella, P., Crutchfield, W.~Y.,
  Fiveland, W.~A., and Jessee, J.~P. (1998).
\newblock An adaptive projection method for unsteady, low-mach number
  combustion.
\newblock {\em Combustion Science and Technology}, 140(1-6):123--168.

\bibitem[Perini et~al., 2012]{perini2012analytical}
Perini, F., Galligani, E., and Reitz, R.~D. (2012).
\newblock An analytical {J}acobian approach to sparse reaction kinetics for
  computationally efficient combustion modeling with large reaction mechanisms.
\newblock {\em Energy \& Fuels}, 26(8):4804--4822.

\bibitem[Rehm and Baum, 1978]{baum1978equations}
Rehm, R.~G. and Baum, H.~R. (1978).
\newblock The equations of motion for thermally driven buoyant flows.
\newblock {\em Journal of Research of the National Bureau of Standards},
  83:297--308.

\bibitem[Safta et~al., 2010]{safta2010high}
Safta, C., Ray, J., and Najm, H.~N. (2010).
\newblock A high-order low-mach number {A}{M}{R} construction for chemically
  reacting flows.
\newblock {\em Journal of Computational Physics}, 229(24):9299--9322.

\bibitem[Speck, 2017]{speck2017parallelizing}
Speck, R. (2017).
\newblock Parallelizing spectral deferred corrections across the method.
\newblock {\em arXiv preprint arXiv:1703.08079}.

\bibitem[Speck et~al., 2015]{speck2015multi}
Speck, R., Ruprecht, D., Emmett, M., Minion, M.~L., Bolten, M., and Krause, R.
  (2015).
\newblock A multi-level spectral deferred correction method.
\newblock {\em BIT Numerical Mathematics}, 55(3):843--867.

\bibitem[Warnatz and Peters, 1982]{warnatz1982numerical}
Warnatz, J. and Peters, N. (1982).
\newblock Numerical methods in laminar flame propagation.
\newblock {\em Vieweg, Braunchweig}, pages 87--111.

\bibitem[Weiser, 2015]{weiser2015faster}
Weiser, M. (2015).
\newblock Faster {S}{D}{C} convergence on non-equidistant grids by {D}{I}{R}{K}
  sweeps.
\newblock {\em BIT Numerical Mathematics}, 55(4):1219--1241.

\bibitem[Xia et~al., 2007]{xia2007efficient}
Xia, Y., Xu, Y., and Shu, C. (2007).
\newblock Efficient time discretization for local discontinuous {G}alerkin
  methods.
\newblock {\em Discrete and Continuous Dynamical Systems Series B}, 8(3):677.

\bibitem[Zhang et~al., 2012]{zhang2012fourth}
Zhang, Q., Johansen, H., and Colella, P. (2012).
\newblock A fourth-order accurate finite-volume method with structured adaptive
  mesh refinement for solving the advection-diffusion equation.
\newblock {\em SIAM Journal on Scientific Computing}, 34(2):B179--B201.

\end{thebibliography}

\end{document}